\documentclass[10pt]{article}
\usepackage{amsmath,amssymb,latexsym,epsfig,subfigure,multirow}
\usepackage{amsthm}
\usepackage{verbatim}
\usepackage{fancyhdr}
\newtheorem{theorem}{Theorem}[section]
\newtheorem{lemma}{Lemma}[section]
\newtheorem{definition}{Definition}[section]

\newtheorem{remark}{Remark}

\newcommand{\pa}{\partial}
\newcommand{\f}{\frac}

\begin{document}

%\pagenumbering{roman}
%\begin{verbatim}

\title{ Second order Method for Solving 3D Elasticity Equations with  Complex Interfaces  }

\author{Bao Wang$^{1}$, Kelin Xia$^{1}$ and Guo-Wei Wei$^{1,2,3}$
\footnote{
Corresponding author. Email: wei@math.msu.edu}  \\
\\
\small \it     $^1$Department of Mathematics, Michigan State University, East Lansing, MI 48824, USA \\
\small \it     $^2$Department of Electrical and Computer Engineering, \\
\small \it          Michigan State University, East Lansing, MI 48824, USA  \\
\small \it     $^3$Center for Mathematical Molecular Biosciences, \\
\small \it          Michigan State University, East Lansing, MI 48824, USA  }

\date{\today}

\maketitle

\begin{abstract}
Elastic materials are ubiquitous in nature and indispensable components in man-made devices and equipments. When a device or equipment involves composite or multiple elastic materials, elasticity interface problems come into play.  The solution of  three dimensional (3D) elasticity interface problems is significantly more difficult than that of elliptic counterparts due to the coupled vector components and cross derivatives in the governing elasticity equation. This  work introduces  the matched interface and boundary (MIB) method for solving  3D elasticity interface problems. The proposed MIB  method utilizes  fictitious values on irregular grid points near the material interface to replace function values in the discretization so that the elasticity equation can be discretized using the standard finite difference schemes as if there were no material interface. The interface jump conditions are rigorously enforced on the intersecting points between the interface and the mesh lines. Such an enforcement determines the fictitious values. A number of new technique are developed to construct efficient MIB schemes for dealing with cross derivative in coupled governing equations. The proposed method is extensively validated over both weak and strong discontinuity of the solution, both piecewise constant and position-dependent material parameters, both smooth and nonsmooth interface geometries, and both small and large contrasts in the Poisson's ratio and shear modulus across the interface. Numerical experiments indicate that the present MIB method is of second order convergence in both $L_\infty$ and $L_2$ error norms.
\end{abstract}
%\vspace*{1cm}
{\it Keywords:}~
Elasticity Interface Problem;
Complex  interface;
Matched interface and boundary.

\newpage
{\setcounter{tocdepth}{4} \tableofcontents}

\newpage

\section{Introduction}
%-----------------------------------------------------------------------------------------
%%%%%%%%%% problem setting
Although materials, such as solids, are composed of atoms or molecules, which are discrete in nature,  continuum models based on the continuum mechanics are highly accurate and applicable to  length scales much greater than that of inter-atomic distances \cite{ZhanChen:2010a}. One of the most widely applied continuum models is elasticity theory, which describes how  solid materials return to their original shapes once being deformed by applied forces.  Linear elasticity theory is often employed when the deformation is relatively small. In such a case,  the stress-strain relation is governed by the constitutive equation. One class of elastic materials is isotropic homogeneous, whose constitutive equations can be uniquely determined with any two terms of six moduli, namely,  bulk modulus, Young's modulus, Lam$\acute{e}$'s first parameter, shear modulus, Poisson's ratio and p-wave modulus \cite{Anandarajah:2010}. For isotropic inhomogeneous materials, the inhomogeneity is often modeled by position-dependent  moduli in their constitutive equations. For example, in seismic wave equations, inhomogeneity is accounted by  position-dependent Lam$\acute{e}$'s  parameters \cite{shearer:1999}. Similar  models have also been employed in the elasticity   analysis of  biomolecules \cite{Wei:2009,Wei:2013,KLXia:2013d}.

Interface description in the elasticity modeling  is  indispensable whenever elastic materials encounter rapid changes or discontinuities in material properties due to    voids, pores, inclusions, dislocations, cracks or composite structures \cite{Dvorak:2013,Fries:2010,Sukumar:2001,Stolarska:2001}. The resulting problem is called  an elasticity interface problem, which is of considerable importance in man-made materials, devices, equipments, tissue engineering, biomedical science and biophysics \cite{Sukumar:2001,Stolarska:2001,Wei:2009,Wei:2013,KLXia:2013d}. Mathematically, discontinuities in  elasticity interface problems can be classified into two types, namely, strong ones and weak ones. Strong discontinuities are referred to situations where the displacement has jumps across the interface, while weak discontinuities have a continuous displacement but with jumps in the gradient of the displacement.   In general, analytical solution to elasticity interface problems is difficult to obtain, except for simple interface geometries. In 1950s, Eshelby found that under a uniformly applied stress, an infinite and elastically isotropic system with an ellipsoidal inhomogeneity has a uniform eigenstrain distribution inside the ellipsoidal domain \cite{Eshelby:1956,Eshelby:1957}. For arbitrarily shaped inhomogeneity,  semianalytic approaches  have been proposed  for finding stress tensors \cite{Mathiesen:2008}.

%%Numerical scheme for elasticity interface problems
Numerical approaches, such as finite element methods (FEMs),  boundary element methods (BEMs)  and finite difference methods (FDMs), are the main workhorse for elasticity interface problems arising from practical applications. Based on computational meshes used, these methods can be classified  into two types, i.e., schemes  utilize body-fitting meshes and algorithms based on special interface schemes. Body-fitting meshes are generated according to geometry of the interface so that no mesh lines cut through the interface. In this type of methods, locally adaptive meshes are frequently employed based on local refinement techniques \cite{XuZL:2003}. In the second type of algorithms, regular meshes that may   cut through the interface are used. Consequently,  sophisticated numerical schemes are needed   to incorporate the interface conditions into element shape functions or operator discretizations. Immersed interface method (IIM) originally proposed for elliptic interface problems  \cite{LeVeque:1994} has been developed to solve  two-dimensional (2D) elasticity interface problems with isotropic homogenous media \cite{YangXZ:2003}. This  finite difference based approach achieves second order accuracy. A second-order sharp numerical method  has been developed for linear elasticity equations \cite{Theillard:2013}.
 Many finite element based methods have also been proposed for elasticity interface problems. Among them, partition of unity method (PUM), the generalized finite element method (GFEM) and extended finite element method (XFEM)    are designed to capture the non-smooth property of the solution over the interface   \cite{Sukumar:2001,Stolarska:2001,Fries:2010}.  Enrichment functions are utilized to handle the material interface.  Discontinuous Galerkin based methods have also been constructed to deal with strong and weak discontinuities \cite{Hansbo:2002,Becker:2009,Mergheim:2006}. Recently,  immerse finite element (IFM) method has been developed to solve elasticity problems with interface jump conditions \cite{LiZL:2005,XieH:2011}. This approach  locally modifies  finite element basis functions  to enforce the jump conditions across the interface. Most recently, a Nitsche type method has been proposed for elasticity interface problems \cite{Michaeli:2013}.

There are few numerical issues in the  solution of  elasticity interface problems. One issue is to deal with complex interface geometry. It is easy to construct a numerical method for some special designed simple interface shapes. However, it is a challenge to automatically  deal with  complex interface geometries. Another issue is to develop robust numerical schemes for handling   interfaces of Lipschitz continuity or geometric singularities, such as cusps, sharp edges,   tips and self-intersecting surfaces \cite{HouSM:2012}. It is still a major challenge to develop second order accurate schemes for arbitrarily complex interface geometries in a three-dimensional (3D) setting. One example of arbitrarily complex interface geometries is the protein molecular surfaces \cite{Yu:2007,Yu:2007a,Geng:2007a}. The other issue is position dependent material parameters. It is necessary for numerical methods to be able to treat spatially varying coefficients. Additionally,  taking care of  strong discontinuities, handling large contrast between material parameters across the interfaces and treatment of the Poisson's ratio near the incompressible limit are also valid numerical issues in elasticity interface problems \cite{LinT:2012,LinT:2013}.
% A locking-free bilinear IFM was proposed to   handle nearly incompressible media and  achieve the second order accuracy  \cite{LinT:2012,LinT:2013}.
Finally, it is a challenge to develop second order accurate schemes for arbitrarily complex interface geometries in three-dimensional (3D) setting. As a vector equation, the existence of three deformation components gives rise to an extraordinary requirement for numerical schemes to be unusually efficient. Although many elegant and efficient algorithms have been  developed for 2D and 3D elasticity interface problems, to our best knowledge, there is little literature about second order convergent schemes for arbitrarily complex interface geometries in 3D, including  interfaces of Lipschitz continuity.

%%MIB algorithms
The matched interface and boundary (MIB) method was originally constructed  for dealing with material interfaces in  Maxwell's equations  \cite{Zhao:2004} and Poisson equation  \cite{Yu:2007,Yu:2007a,Zhou:2006c,Zhou:2006d,Geng:2007a}. The essential idea of the MIB method is to extend the solution beyond the interface so that the derivatives near the interface can be discretized as if there were no interface. The extension along the interface is carried out by iteratively incorporating the  lowest order of interface conditions so that  in principle, arbitrarily high order accuracy can be achieved. Sixteenth order accuracy was achieved for simple interface geometries  \cite{Zhao:2004,Zhou:2006c} and up to sixth order convergence was realized for 3D complex interface shapes i \cite{Yu:2007a}.  For arbitrarily complex interfaces with geometric singularities, robust second order numerical convergence was observed \cite{Yu:2007,Yu:2007a,Geng:2007a}.   In the past decade, the MIB method has been successfully applied to a variety of problems. For example, in computational biophysics, an MIB based Poisson-Boltzmann solver, MIBPB \cite{DuanChen:2011a}, has been developed   for the analysis of the electrostatic potential of biomolecules \cite{Yu:2007, Geng:2007a,Zhou:2008b}, molecular dynamics \cite{Geng:2011} and charge transport phenomenon \cite{QZheng:2011a, QZheng:2011b}.  Zhao has constructed second order and fourth order  MIB schemes for the Helmholtz problems \cite{SZhao:2010a,SZhao:2008a}. A second order MIB method has been  developed by Zhou and coworkers to solve the Navier-Stokes equations with discontinuous viscosity and density \cite{YCZhou:2012a}. Recently, the MIB method has been extended to solve elliptic equations with multi-material interfaces \cite{KLXia:2011}.
%Most recently, the second order MIB scheme has been constructed for  2D elasticity interface problems \cite{BWang:2014a}.  Second order  convergence was observed for both piecewise constant material parameters and spatial dependent material parameters via several different kind of interfaces. The robustness of the MIB elasticity interface algorithm was tested through both large and small contrasts of the two Lam\'{e}'s parameters \cite{BWang:2014a}.

%%%%%%%%%%% MIB for elasticity in 3D
The object of the present work is to develop MIB schemes for solving 3D elasticity interface problems. We consider both smooth and sharp interfaces  for isotropic homogeneous and inhomogeneous elastic materials. First, we extend our earlier MIB method  for elliptic interface problems to elasticity counterparts. To this end, we take care of both central derivatives and cross derivatives in the elasticity equation. Several numerical techniques
namely, disassociation,
% which will employ the fictitious values for the central derivatives to the discretization of the cross derivatives;
extrapolation, % which use function values or fictitious values at three grid points to extrapolate the secondary fictitious values for the cross derivatives;
and neighbor combination, % which is powerful to handle the limiting case that there is only one grid point from one side of the interface.
are proposed to compute the fictitious values for the discretization of the cross derivatives.
Additionally, to make the present MIB method efficient for dealing with three coupled vector components, we carefully optimize our algorithms so that the resulting discretization matrix is as symmetric and diagonally dominant as possible.   Moreover, to handle geometric singularities, we develop a technique to simultaneously employ two sets of interface conditions from two intersecting points where the interface meets mesh lines.
%more class of interface conditions from another intersection point, proposed to handle the geometric singularities in elliptical interface problems \cite{Yu:2007}, to find the fictitious values for the discretization of the central derivatives of the elasticity interface problems.
Finally, we validate the proposed MIB for wide variety of  elasticity interface problems, including large contrast in material parameters across the interface, strong interface discontinuity, sharp-edged interface and variable material coefficients.

The rest of this paper is organized as follows.  The formulation of 3D elasticity interface problems is presented in Section \ref{theory}. Section \ref{algorithm} is devoted to the construction of MIB algorithms for elasticity interface problems. Methods for determining  fictitious values are proposed for both  central derivatives and cross derivatives in the elasticity equation.  The present methods are  extensively validated by analytical tests with complex interface geometries, including interfaces of Lipschitz continuity in Section \ref{validation}. We demonstrate that the second order accuracy is achieved by the proposed  MIB method.  This paper ends with a conclusion.

%=================================================================================================================================
\section{Formulation of the elasticity interface problem}
\label{theory}

The 3D linear elasticity motion considered in the present work is governed by the following linear elasticity equations
\begin{equation}
\label{elas}
\nabla \cdot \mathbb{T}+\mathbf{F}=\frac{d^2\mathbf{u}}{dt^2},
\end{equation}
where $\mathbb{T}$ is the strain tensor, $\mathbf{F}=\left(F_1(\mathbf{x}), F_2(\mathbf{x}), F_3(\mathbf{x})\right)^T:=(F_1, F_2, F_3)^T$ is the external force on the elasticity body, $\mathbf{u}=\left(u_1(\mathbf{x}), u_2(\mathbf{x}), u_3(\mathbf{x})\right)^T$ is a displacement vector, $\mathbf{x}=\left(x, y, z\right)^T$ is a position vector, and $*^T$ is the transpose of quantity $*$.
\\

For  isotropic homogeneous media, the strain tensor $\mathbb{T}$ is a 3 by 3 symmetric matrix which has the form
\begin{equation}
\mathbb{T}=\lambda \rm{tr}(\sigma)I+2\mu\sigma,
\end{equation}
where $\lambda$ is the Lam\'{e}'s parameter, $\mu$ is the shear modulus, $I$ is a 3 by 3 identity matrix, and $\sigma$ is a stress tensor which can be further written as
\begin{equation}
\sigma=\frac{1}{2}\left(\nabla \mathbf{u}+(\nabla\mathbf{u})^T\right).
\end{equation}
The static state elasticity equation is given by
\begin{equation} \label{elas_s}
\nabla \cdot \mathbb{T}+\mathbf{F}=\mathbf{0}.
\end{equation}
In the present work, we focus on the static state problem (\ref{elas_s}).

\subsection{Interface jump conditions}
\begin{figure}[!ht]
\centering
\includegraphics[width=6cm,height=4cm]{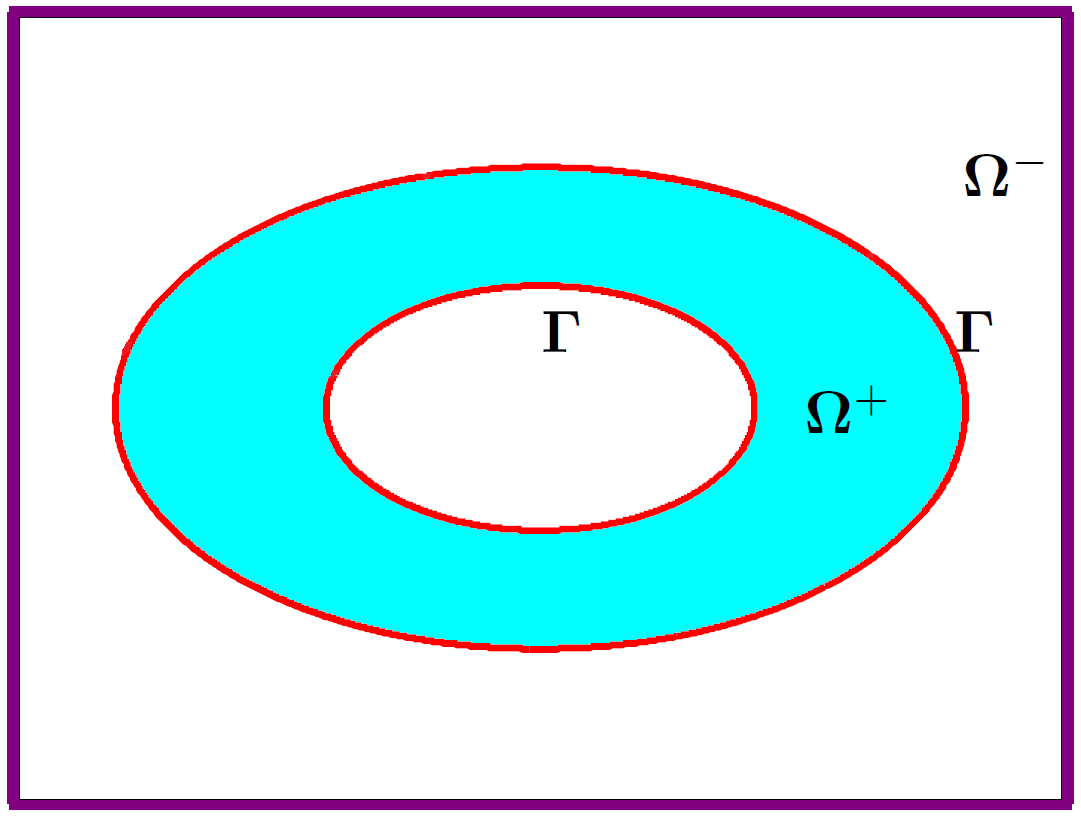}
\caption{Illustration of the elasticity interface problem at a cross section $(x=x_i)$. The whole domain consists of two subdomains $\Omega^+$ and $\Omega^-$ by the interface $\Gamma$.}
\label{inter_pro}
\end{figure}

Consider a two-phase elastic body having two different elastic materials in a rectangular prism domain $\Omega \subset \mathbb{R}^3$.  The two phase elastic motion is separated by an arbitrarily complex interface $\Gamma$, which splits the whole domain $\Omega$ into $\Omega^+$ and $\Omega^-$, i.e., $\Omega=\Omega^+\cup\Gamma\cup\Omega^-$, as illustrated in Fig. \ref{inter_pro}.

\begin{lemma}
For the 3D elasticity equations of the static state (\ref{elas_s}), if the source term $\mathbf{F}$ has a potential function representation $U$, i.e., $\nabla U=\mathbf{F}$, then the static state elasticity equations can be written as a homogeneous equation. More precisely, there exist another 3 by 3 matrix $\tilde{\mathbb{T}}$ such that
$$
\nabla\cdot\tilde{\mathbb{T}}=\mathbf{0},
$$
where $\mathbf{0}$ is the 3D zero vector.
\end{lemma}

\begin{definition}
\textbf{Weak Solution:} $\tilde{\mathbb{T}}$ is said to be the weak solution of the homogeneous equation $\nabla\cdot\tilde{\mathbb{T}}=\mathbf{0}$ provided
$$
\int_{\Omega}\nabla\phi\cdot \tilde{\mathbb{T}}d{\bf r}=0,
$$
holds for  $\forall \phi\in C^\infty_0(\Omega)$, where $C^\infty_0(\Omega)$ is the space of smooth functions with compact support on $\Omega$, and $d{\bf r}=dxdydz$ is the volume integral element.
\end{definition}

\begin{theorem}
Let $\mathbb{T}$ be a second order tensor in $\mathbb{R}^3$ which can be written as an $3$ by $3$ matrix. For the elasticity equations
\begin{equation}
\nabla \cdot \mathbb{T}+\mathbf{F}=\mathbf{0},
\end{equation}
where $\mathbf{F}$ is an $3$-dimensional vector valued function  and $\mathbf{0}\in \mathbb{R}^3$. If the force term has a potential function $U$, i.e., $\nabla U=\mathbf{F}$, then across the interface, the weak solution satisfies the following interface conditions
\begin{equation}
\label{inter2}
[\mathbb{T}\cdot \mathbf{n}]=\mathbf{T},
\end{equation}
where $\mathbf{T}$ is an $3$-dimensional vector-valued function, $[*]$ is the jump of   quantity $*$ across the interface, and
$\mathbf{n}$ is the normal direction of the interface.
\end{theorem}

\begin{remark}
If the material has no fracture which corresponds to the weak discontinuity in the linear elasticity interface problem, the following interface condition is enforced in traditionally elasticity interface problems
$$
[\mathbf{u}]=\mathbf{0}.
$$
However, fractures often occur in realistic materials which corresponds to the strong discontinuity.  In this work, our numerical scheme is designed for both strong and weak discontinuity of elasticity interface problems.
\end{remark}

\subsection{Linear elasticity interface problem}
In this work, we only consider the static state elasticity equation (\ref{elas_s}). As discussed above, the 3D elasticity interface problem can be formulated as
\begin{eqnarray}
     \nabla\cdot\mathbb{T}+\mathbf{F} &=& \mathbf{0}, \quad \mbox{in} \quad  \Omega^+\cup\Omega^-,\\
    \left[\mathbf{u}\right] \mid_\Gamma&=& \mathbf{b}, \quad \mbox{on} \quad \Gamma,\\
    \left[\mathbb{T}\cdot\mathbf{n}\right]\mid_\Gamma &=& \mathbf{T}, \quad \mbox{on} \quad \Gamma,\\
    \mathbf{u} &=& \mathbf{u}^0, \quad  \mbox{on} \quad \partial\Omega.
\end{eqnarray}
where $\mathbf{u}=(u_1, u_2, u_3)^T: \Omega\rightarrow \mathbb{R}^3$ is the displacement field and $\mathbf{n}=(n_1, n_2, n_3)^T$ is the unit outer normal vector to the interface $\Gamma$. Function $\mathbf{F}$, as stated above, is a 3D vector-valued function of body force field. Vector $\mathbf{u}^0=(u_1^0, u_2^0, u_3^0)^T$ is the Dirichlet boundary conditions. For elasticity interface problem, generally, if vector $\mathbf{b}=(b_1, b_2, b_3)$ does not equal 0, it is called strong discontinuity, otherwise weak discontinuity. Here vector valued function $\mathbf{T}=(\phi, \psi, \eta)^T$ are the jump of the traction $\mathbb{T}\cdot\mathbf{n}$ across the interface $\Gamma$.

In material science, the stress-strain relation is usually described by constitutive equation, which in terms of Lam\'{e}'s parameters can be expressed as,
\begin{equation}
\label{strain_stress}
\mathbb{T}=\lambda \rm{tr}(\sigma)I+2\mu\sigma.
\end{equation}
Here stain tensor $\sigma$ is defined as,
$$
\sigma=\frac{1}{2}\left(\nabla \mathbf{u}+(\nabla \mathbf{u})^T\right).
$$

Dramatic different elasticity behaviors can be observed between inhomogeneous media and homogeneous media. To take this property into consideration, we elaborate the elasticity interface problem in both situations.
For inhomogeneous material, Lam\'{e}'s parameters are position dependent, i.e., $\lambda=\lambda(x,y,z))$ and $\mu=\mu(x,y,z))$. Using the constitutive equation in Eq. (\ref{strain_stress}), the governing equation of elasticity interface problem can be expressed as,
\begin{eqnarray}
\label{el11}
\nabla\lambda(\nabla \cdot \mathbf{u})+\nabla\mu \cdot \left[\nabla\mathbf{u}+(\nabla\mathbf{u})^T \right] +(\lambda+\mu)\nabla\nabla \cdot \mathbf{u} +\mu \nabla^2\mathbf{u}=-\mathbf{F}.
\end{eqnarray}
We can spell out all the terms as following,
\begin{equation}
\label{el11}
\resizebox{0.92\hsize}{!}{$%
(\lambda+2\mu)\frac{\pa^2 u_1}{\pa x^2}+\mu\frac{\pa^2 u_1}{\pa y^2}+\mu\frac{\pa^2 u_1}{\pa z^2}+(\lambda+\mu)\frac{\pa^2 u_2}{\pa x\pa y}+(\lambda+\mu)\frac{\pa^2 u_3}{\pa x\pa z}+
\lambda_x\left(\frac{\pa u_1}{\pa x}+\frac{\pa u_2}{\pa y}+\frac{u_3}{\pa z}\right)+2\mu_x\frac{\pa u_1}{\pa x}+\mu_y\left(\frac{\pa u_1}{\pa y}+\frac{\pa u_2}{\pa x}\right)+\mu_z\left(\frac{u_1}{\pa z}+\frac{\pa u_3}{\pa x}\right)=-F_1,
$%
}%
\end{equation}
\begin{eqnarray}
\label{el12}
\resizebox{0.92\hsize}{!}{$%
\mu\frac{\pa^2 u_2}{\pa x^2}+(\lambda+2\mu)\frac{\pa^2 u_2}{\pa y^2}+\mu\frac{\pa^2 u_2}{\pa z^2}+(\lambda+\mu)\frac{\pa^2 u_1}{\pa x\pa y}+(\lambda+\mu)\frac{\pa^2 u_3}{\pa y\pa z}+
\mu_x\left(\frac{\pa u_2}{\pa x}+\frac{\pa u_1}{\pa y}\right)+\lambda_y\left(\frac{\pa u_1}{\pa x}+\frac{\pa u_2}{\pa y}+\frac{\pa u_3}{\pa z}\right)+2\mu_y\frac{\pa u_2}{\pa y}+\mu_z\left(\frac{\pa u_2}{\pa z}+\frac{\pa u_3}{\pa y}\right)=-F_2,
$%
}%
\end{eqnarray}
\begin{eqnarray}
\label{el13}
\resizebox{0.92\hsize}{!}{$%
\mu\frac{\pa^2 u_3}{\pa x^2}+\mu\frac{\pa^2 u_3}{\pa y^2}+(\lambda+2\mu)\frac{\pa^2 u_3}{\pa z^2}+(\lambda+\mu)\frac{\pa^2 u_1}{\pa x\pa z}+(\lambda+\mu)\frac{\pa^2u_2}{\pa y\pa z}+
\mu_x\left(\frac{\pa u_3}{\pa x}+\frac{\pa u_1}{\pa z}\right)+\mu_y\left(\frac{\pa u_3}{\pa y}+\frac{\pa u_2}{\pa z}\right)+\lambda_z\left(\frac{\pa u_1}{\pa x}+\frac{\pa u_2}{\pa y}+\frac{\pa u_3}{\pa z}\right)+
2\mu_z\frac{\pa u_3}{\pa z}=-F_3.
$%
}%
\end{eqnarray}
With the constitutive equations, the jump conditions regarding to the stain tensor can be represented as,
\begin{eqnarray}
\label{el17}
\left[\left(\lambda\left(\frac{\pa u_1}{\pa x}+\frac{\pa u_2}{\pa y}+\frac{\pa u_3}{\pa z}\right)+2\mu \frac{\pa u_1}{\pa x} \right)n_1+\mu\left(\frac{\pa u_2}{\pa x}+\frac{\pa u_1}{\pa y}\right)n_2+\mu\left(\frac{\pa u_3}{\pa x}+\frac{\pa u_1}{\pa z}\right)n_3\right]\mid_\Gamma=\phi, \ \mbox{on}\ \Gamma, \\
\left[\mu\left(\frac{\pa u_1}{\pa y}+\frac{\pa u_2}{\pa x}\right)n_1+\left(\lambda\left(\frac{\pa u_1}{\pa x}+\frac{\pa u_2}{\pa y}+\frac{\pa u_3}{\pa z}\right)+2\mu\frac{\pa u_2}{\pa y} \right)n_2+\mu\left(\frac{\pa u_3}{\pa y}+\frac{\pa u_2}{\pa z}\right)n_3\right]\mid_\Gamma=\psi, \ \mbox{on}\ \Gamma,  \\
\left[\mu\left(\frac{\pa u_1}{\pa z}+\frac{\pa u_3}{\pa x}\right)n_1+\mu\left(\frac{\pa u_2}{\pa z}+\frac{\pa u_3}{\pa y}\right)n_2+\left(\lambda\left(\frac{\pa u_1}{\pa x}+\frac{\pa u_2}{\pa y}+\frac{\pa u_3}{\pa z}\right)+2\mu\frac{\pa u_3}{\pa z} \right)n_3\right]\mid_\Gamma=\eta, \ \mbox{on}\ \Gamma.
\end{eqnarray}
Together with the Dirichlet boundary conditions and the jump conditions,
%regarding to displacement as following
%\begin{eqnarray}\label{el14}
%&&\left[u_1\right]\mid_\Gamma=b_1,  \ \mbox{on}\ \Gamma, \\
%&&\left[u_2\right]\mid_\Gamma=b_2,  \ \mbox{on}\ \Gamma, \\
%&&\left[u_3\right]\mid_\Gamma=b_3,  \ \mbox{on}\ \Gamma,
%\end{eqnarray}
we set up the general formulation for linear elasticity interface problem with inhomogeneous media.

For homogeneous material, algebraic relations exist between different elasticity moduli, i.e., Bulk modulus $K$, Young's modulus $E$, Lam \'{e}'s first parameter $\lambda$, Shear modulus $\mu$, Poisson's ratio $\nu$  and P-wave modulus $M$. For instance, Lam\'{e}'s parameters can be represented by Young's modulus $E$ and Poisson's ratio $\nu$ as,
$$
\mu=\f{E}{2(1+\nu)},\lambda=\f{E\nu}{(1+\nu)(1-2\nu)}.
$$
Due to the constant moduli, the governing equation can be simplified as
\begin{eqnarray}
\label{el11_simplify}
(\lambda+\mu)\nabla\nabla \cdot \mathbf{u} +\mu \nabla^2\mathbf{u}=-\mathbf{F}.
\end{eqnarray}

With the above algebraic relations of elasticity moduli, the governing equation can be further written as,
\begin{eqnarray}\label{elas_eq1}
2(1-\nu)\f{\pa^2 u_1}{\pa x^2}+(1-2\nu)\f{\pa^2 u_1}{\pa y^2}+(1-2\nu)\f{\pa^2 u_1}{\pa z^2}+\f{\pa^2 u_2}{\pa x\pa y}+\f{\pa^2 u_3}{\pa x\pa z}=f_1,\\
(1-2\nu)\f{\pa^2 u_2}{\pa x^2}+2(1-\nu)\f{\pa^2 u_2}{\pa y^2}+(1-2\nu)\f{\pa^2 u_2}{\pa z^2}+\f{\pa^2 u_1}{\pa x\pa y}+\f{\pa^2 u_3}{\pa y\pa z}=f_2,\\
(1-2\nu)\f{\pa^2 u_3}{\pa x^2}+(1-2\nu)\f{\pa^2 u_3}{\pa y^2}+2(1-\nu)\f{\pa^2 u_3}{\pa z^2}+\f{\pa^2 u_1}{\pa x\pa z}+\f{\pa^2 u_2}{\pa y\pa z}=f_3.
\end{eqnarray}
Here $(f_1,f_2,f_3)$ are prerequisite terms, and they can be related to the body force by $(f_1,f_2,f_3)=(-\frac{F_1}{\mu+\lambda},-\frac{F_2}{\mu+\lambda},-\frac{F_3}{\mu+\lambda})$.

Also the second set of jump conditions can be rewritten as,
\begin{eqnarray}
\label{deri_jump1}
\left[\f{2\mu}{1-2\nu}\left((1-\nu)\f{\pa u_1}{\pa x}+\nu\f{\pa u_2}{\pa y}+\nu\f{\pa u_3}{\pa z}\right)n_1+\mu\left(\f{\pa u_1}{\pa y}+\f{\pa u_2}{\pa x}\right)n_2+\mu\left(\f{\pa u_1}{\pa z}+\f{\pa u_3}{\pa x}\right)n_3\right]|_\Gamma=\phi, \ \mbox{on}\ \Gamma, \\
\left[\mu\left(\f{\pa u_1}{\pa y}+\f{\pa u_2}{\pa x}\right)n_1+\f{2\mu}{1-2\nu}\left(\nu\f{\pa u_1}{\pa x}+(1-\nu)\f{\pa u_2}{\pa y}+\nu\f{\pa u_3}{\pa z}\right)n_2+\mu\left(\f{\pa u_3}{\pa y}+\f{\pa u_2}{\pa z}\right)n_3\right]|_\Gamma=\psi, \ \mbox{on}\ \Gamma, \\
\left[\mu\left(\f{\pa u_1}{\pa z}+\f{\pa u_3}{\pa x}\right)n_1+\mu\left(\f{\pa u_2}{\pa z}+\f{\pa u_3}{\pa y}\right)n_2+\f{2\mu}{1-2\nu}\left(\nu\f{\pa u_1}{\pa x}+\nu\f{\pa u_2}{\pa y}+(1-\nu)\f{\pa u_3}{\pa z}\right)n_3\right]|_\Gamma=\rho, \ \mbox{on}\ \Gamma,
\end{eqnarray}

The above two sets of equations, together with the Dirichlet boundary conditions and jump conditions regarding to displacement, constitute general formulation for linear elasticity interface problem in homogeneous media.

%=================================================================================================================================
\section{Method and algorithm}\label{algorithm}

In this section, the MIB method for elliptical interface problems is extended to solve elasticity interface problems.
%Let the computational domain be $\Omega=[a, b]\times[c, d]\times[e, f]$ and denote $h$  the grid size in griding the computational domain.  We set  the grid points to be
%$$
%x_i=a+(i-1)h,\ y_j=c+(j-1)h, \ z_k=e+(k-1)h,\ i=1, 2, \cdots, n_x;\ j=1, 2, \cdots, n_y;\ k=1, 2, \cdots, n_z,
%$$
%where $n_x, n_y$ and $n_z$ are the total number of grid points along the $x$- $y$- and $z$-directions, respectively. We denote grid point $(i, j, k)=:(x_i,y_j,z_k)$.
Due to the existence of the interface,  the direct application of the standard central finite difference (CFD) schemes leads to a dramatic decrease in the accuracy
and convergence of the numerical solution. To maintain the designed order of accuracy, MIB method extends function values across the interface. The resulting extended function values are called fictitious values, which are employed, together with function values on the other side of the interface, for the CFD discretization of the PDE across the interface. For example, at a grid point $(i, j, k)$ near the interface, if its finite difference scheme refers to some grid points that are in the other side of the interface, fictitious values from other side of the interface  are utilized in the finite difference  discretization.
To extend  the function values to the other side of the interface and enable the MIB discretization of second order convergence,
the interface conditions on both function values and normal derivatives are utilized and enforced.
% which yields the continuity of both function values and derivatives, and thus  the extension is $C^1$.

The location of a fictitious value is called an irregular grid point, while those grid points where no fictitious value is required are called regular grid points.  Loosely speaking, irregular grid points form extended domains on both sides of the interface. The extended domains ensure that the standard central finite difference scheme can be uniformly applied without the  loss of numerical accuracy.

Additionally, derivatives involved in the elasticity equation are classified into central derivatives and cross derivatives. Central derivatives involve only one direction, while cross derivatives refer to more than one direction. These two situations are to be handled in different manners in the present method. Additional care is needed for discretizing cross derivatives to the second order accuracy.

Moreover, interfaces are classified into smooth ones and nonsmooth ones. The nonsmooth interfaces are Lipschitz continuous with geometric singularities, such as cusps, tips and/or sharp edges. To maintain designed order of accuracy, nonsmooth interfaces are much more difficult to deal with.

\subsection{General MIB algorithms for Laplace operator}\label{smoothinterface}

\subsubsection{Simplification of interface jump conditions}

As the interface normal direction varies along the interface, which is very troublesome from a computational perspective. It is necessary to define a set of local coordinates at each intersection point of the interface and the Cartesian mesh, so that different interface geometries can be treated in a systematical manner. In this section, we present the local coordinate transformation formula. At a specific intersection point,  the local coordinate system is chosen to be $(\xi, \eta, \zeta)$, where $\xi$ is along the normal direction and $\eta$ is in the $xy$ plane. This local coordinate system can be obtained from the Cartesian coordinate system via the following transformation
\begin{equation}
\label{trans}
\left(
  \begin{array}{c}
    \xi \\
    \eta\\
    \zeta
  \end{array}
\right)=
\mathbf{P}\cdot
\left(
  \begin{array}{c}
    x \\
    y\\
    z
  \end{array}
\right),
\end{equation}
where $\mathbf{P}\doteq \{P(i, j)\}_{i, j=1, 2, 3}$ is a transformation matrix
\begin{equation}
\label{trans_mat}
\mathbf{P}=
\left(
  \begin{array}{ccc}
    \sin\phi\cos\theta & \sin\phi\sin\theta &\cos\phi \\
    -\sin\theta & \cos\theta &0\\
    -\cos\phi\cos\theta &-\cos\phi\sin\theta &\sin\phi
  \end{array}
\right),
\end{equation}
where $\theta$ and $\phi$ are the azimuth and zenith angle with respect to the normal direction, respectively.

In the new local coordinate system, the interface conditions on function values and normal derivatives become (here for simplicity, we only discuss the constant material parameter case, the case of spatially dependent material parameters  can be treated similarly)
\begin{equation}
\label{jump_val}
[\mathbf{u}]|_\Gamma=\mathbf{b},
\end{equation}
and
\begin{equation}
\label{deri_jump_normal}
[\mathbb{T}\cdot\xi]|_\Gamma=\mathbf{T}.
\end{equation}

To achieve better stability and higher efficiency, which is essential for the present 3D vector equation, only the lowest order jump conditions are utilized in the MIB method.  Therefore, we avoid generating high order (derivative) jump conditions, even if in arbitrarily high order MIB methods \cite{Zhao:2004,Zhou:2006c,Yu:2007}. However, we hope to have as many low order jump conditions as possible so as to gain flexibility in dealing with complex interface geometries. To this end, we differentiate the jump condition of the vector function to derive two additional sets of interface jump conditions along
$\eta$ and  $\zeta$ directions, respectively
\begin{equation}
\label{jump1}
[\mathbf{u}_\eta]|_\Gamma=\left(-\sin\theta \frac{\pa \mathbf{u^+}}{\pa x}+\cos\theta \frac{\pa \mathbf{u^+}}{\pa y}\right)-\left(-\sin\theta \frac{\pa \mathbf{u^-}}{\pa x}+\cos\theta \frac{\pa \mathbf{u^-}}{\pa y}\right),
\end{equation}
and
\begin{equation}
\label{jump2}
[\mathbf{u}_\zeta]|_\Gamma=\left(-\cos\theta\cos\theta \frac{\pa \mathbf{u^+}}{\pa x}-\cos\phi\sin\theta \frac{\pa \mathbf{u^+}}{\pa y}+\sin\phi \frac{\pa \mathbf{u^+}}{\pa z}\right)-\left(-\cos\theta\cos\theta \frac{\pa \mathbf{u^-}}{\pa x}-\cos\phi\sin\theta \frac{\pa \mathbf{u^-}}{\pa y}+\sin\phi \frac{\pa \mathbf{u^-}}{\pa z}\right).
\end{equation}
where $\mathbf{u}=(u_1, u_2, u_3)^T$.

In summary, at a specific intersection point of the interface and the mesh line, there are four sets of interface conditions (\ref{jump_val})-(\ref{jump2}), which only refer to the function values and lowest order derivatives. This property is crucial to endow the MIB method with high efficiency and stability in handling complex interface geometries since no higher order derivative is referred in determining fictitious values. Additionally, lowest order derivatives lead to a more banded matrix and a smaller condition number, which are crucial in solving 3D vector interface problems.

%\subsubsection{Elimination of interfacial derivatives}

In the MIB method, the function values near the interface are extended across the interface by introducing fictitious values. The extension is done along one mesh line at a time, so that it is locally a 1D-like scheme for a higher dimensional interface. Fictitious values can be determined by the aforementioned interface conditions (\ref{jump_val})-(\ref{jump2}).
These conditions involve eighteen derivatives $\frac{\pa \mathbf{u}^{\pm}}{\pa x_j}$, where $x_j=x,y,z$ and ${\bf u}=(u_1,u_2,u_3)^T$. These derivatives are to be evaluated on the interface and thus are called interfacial derivatives.  Due to the geometric complexity, some of these eighteen inetrfacial derivatives can be very difficult to compute numerically.
In general, these interfacial derivatives are grouped into six sets because $u_1$ $u_2$ and $u_3$ can be treated in a similar manner in most situations.

In a second order scheme, we typically have two (sets of) fictitious values along one specific mesh line at one time. However, there are four sets of interface conditions.  Therefore,
 two sets of interface conditions are  redundant. This redundancy gives two more degrees of freedom for us to design efficient and robust second order schemes in a complex interface geometry. Our basic idea is to  algebraically eliminate two sets of interfacial derivatives that are the most difficult to compute by using two sets of redundant interface conditions.   Therefore,  at each intersection point we only need to evaluate four sets of derivatives that are relatively easy to approximate numerically.

The  two sets of derivatives that are to be eliminated are selected by the following principles.
\begin{itemize}
\item Two sets of  fictitious values along the mesh line that intersects with the interface are determined at one time.

\item Two sets of derivatives along the mesh line that intersect with the interface must be kept.

\item In the remaining four sets of derivatives, select two sets that are most  difficult to evaluate due  to the local geometry and eliminate them by two sets of jump conditions.
\end{itemize}

Denote $\mathbf{T}:=(T_1, T_2, T_3)^T$ in interface conditions (\ref{deri_jump_normal}), and    by further introducing the matrix notation, the interface conditions (\ref{deri_jump_normal}-\ref{jump2}) can be rewritten as:
\begin{eqnarray*}
\label{jump}
&\left(T_1, T_2, T_3, [\frac{\pa u_1}{\pa \eta}], [\frac{\pa u_2}{\pa \eta}], [\frac{\pa u_3}{\pa \eta}], [\frac{\pa u_1}{\pa \zeta}], [\frac{\pa u_2}{\pa \zeta}], [\frac{\pa u_3}{\pa \zeta}] \right)^T
\\&=
\mathbf{C}
\left(\frac{\pa u^+_1}{\pa x}, \frac{\pa u^-_1}{\pa x}, \frac{\pa u^+_1}{\pa y}, \frac{\pa u^-_1}{\pa y}, \frac{\pa u^+_1}{\pa z}, \frac{\pa u^-_1}{\pa z}, \frac{\pa u^+_2}{\pa x}, \frac{\pa u^-_2}{\pa x}, \frac{\pa u^+_2}{\pa y}, \frac{\pa u^-_2}{\pa y}, \frac{\pa u^+_2}{\pa z}, \frac{\pa u^-_2}{\pa z},
\frac{\pa u^+_3}{\pa x}, \frac{\pa u^-_3}{\pa x}, \frac{\pa u^+_3}{\pa y}, \frac{\pa u^-_3}{\pa y}, \frac{\pa u^+_3}{\pa z}, \frac{\pa u^-_3}{\pa z}\right)^T
\end{eqnarray*}
where
$$
\mathbf{C}=
\left(
  \begin{array}{ccc}
    C_{1, 1} & C_{1, 2} &C_{1, 3} \\
    C_{2, 1} & C_{2, 2} &C_{2, 3}\\
    C_{3, 1} &C_{3, 2} &C_{3, 3}
  \end{array}
\right),
$$

$$
C_{1, 1}=
\left(
\begin{array}{cccccc}
    M^+P(1, 1) & -M^-P(1, 1) &\mu^+P(1, 2) &-\mu^-P(1, 2) &\mu^+P(1, 3) &-\mu^-P(1, 3)\\
    \lambda^+P(1, 2) & -\lambda^-P(1, 2) &\mu^+P(1, 1) &-\mu^-P(1, 1) &0 &0\\
    \lambda^+P(1, 3) & -\lambda^-P(1, 3) &0 &0 &\mu^+P(1, 1) &-\mu^-P(1, 1)
  \end{array}
\right)
$$

$$
C_{1, 2}=
\left(
\begin{array}{cccccc}
    \mu^+P(1, 2) & -\mu^-P(1, 2) &\lambda^+P(1, 1) &-\lambda^-P(1, 1) &0 &0\\
    \mu^+P(1, 1) & -\mu^-P(1, 1) &M^+P(1, 2) &-M^-P(1, 2) &\mu^+P(1, 3) &-\mu^-P(1, 3)\\
    0 & 0 &\lambda^+P(1, 3) &-\lambda^-P(1, 3) &\mu^+P(1, 2) &-\mu^-P(1, 2)
  \end{array}
\right)
$$

$$
C_{1, 3}=
\left(
\begin{array}{cccccc}
    \mu^+P(1, 3) & -\mu^-P(1, 3) &0 &0 &\lambda^+P(1, 1) &\lambda^-P(1, 1)\\
    0 &0&  \mu^+P(1, 3) &-\mu^-P(1, 3) &\lambda^+P(1, 2) & -\lambda^-P(1, 2)\\
    \mu^+P(1, 1) & -\mu^-P(1, 1) &\mu^+P(1, 2) &-\mu^-P(1, 2) &M^+P(1, 3) &-M^-P(1, 3)
  \end{array}
\right)
$$

$$
C_{2, 1}=
\left(
\begin{array}{cccccc}
    P(2, 1) & -P(2, 1) &P(2, 2) &-P(2, 2) &0 &0\\
    0 &0& 0 &0 &0 & 0\\
    0 &0& 0 &0 &0 & 0
  \end{array}
\right)
$$

$$
C_{2, 2}=
\left(
\begin{array}{cccccc}
    0 &0& 0 &0 &0 & 0\\
    P(2, 1) & -P(2, 1) &P(2, 2) &-P(2, 2) &0 &0\\
    0 &0& 0 &0 &0 & 0
  \end{array}
\right)
$$

$$
C_{2, 3}=
\left(
\begin{array}{cccccc}
    0 &0& 0 &0 &0 & 0\\
    0 &0& 0 &0 &0 & 0\\
    P(2, 1) & -P(2, 1) &P(2, 2) &-P(2, 2) &0 &0
  \end{array}
\right)
$$

$$
C_{3, 1}=
\left(
\begin{array}{cccccc}
    P(3, 1) & -P(3, 1) &P(3, 2) &-P(3, 2) &P(3, 3) &-P(3, 3)\\
    0 &0& 0 &0 &0 & 0\\
    0 &0& 0 &0 &0 & 0
  \end{array}
\right)
$$

$$
C_{3, 2}=
\left(
\begin{array}{cccccc}
    0 &0& 0 &0 &0 & 0\\
    P(3, 1) & -P(3, 1) &P(3, 2) &-P(3, 2) &P(3, 3) &-P(3, 3)\\
    0 &0& 0 &0 &0 & 0
  \end{array}
\right)
$$
$$
C_{3, 3}=
\left(
\begin{array}{cccccc}
    0 &0& 0 &0 &0 & 0\\
    0 &0& 0 &0 &0 & 0\\
    P(3, 1) & -P(3, 1) &P(3, 2) &-P(3, 2) &P(3, 3) &-P(3, 3)
  \end{array}
\right).
$$

In the above expressions, $M=\frac{2\mu(1-\nu)}{1-2\nu}$, $\lambda$ and $\mu$ are the $p$-wave module, bulk modulus and shear modulus, respectively.
Here $*^+$ and  $*^-$ are the limiting values of the quantity $*$ inside and outside the interface, respectively.

\begin{lemma}
Consider the matrix:
$$
A\doteq
\left(
\begin{array}{cccccc}
    M^+P(1, 1) & -M^-P(1, 1) &\mu^+ P(1, 2) &-\mu^- P(1, 2) &\mu^+ P(1, 3) &-\mu^-P(1, 3)\\
    P(2, 1) &-P(2, 1)&P(2, 2) &-P(2, 2) &0 &0\\
    P(3, 1) &-P(3, 1) &P(3, 2) &-P(3, 2) &P(3, 3) &-P(3, 3)
  \end{array}
\right)
$$
where $M^+, M^-, \mu^+, \mu^-, P(i, j), i, j=1, 2, 3$ are the same as above. Then $\forall 1\leq l, m\leq 6, l\neq m$, there
exists constants $a, b, c$ such that the $l$-th and $m$-th column of the vector $aA(1, :)+bA(2, :)+cA(3, :)$ are both zero, where
$A(1, :), A(2, :), A(3, :)$ are the first, second and the last row of the matrix $A$.
\end{lemma}

\begin{proof}
If $l=5, m=6$ or $l=6, m=5$ we simply let $a=0, b=1, c=0$ then it is obvious that the $5$-th and $6$-th column of the vector
$aA(1, :)+bA(2, :)+cA(3, :)$ are both zero.
\\

Otherwise, we let:
$$
a=A(2, l)A(3, m)-A(3, l)A(2, m),
$$
$$
b=A(3, l)A(1, m)-A(1, l)A(3, m),
$$
$$
c=A(1, l)A(2, m)-A(2, l)A(1, m),
$$
then we have the $l$-th and $m$-th column of the vector
$aA(1, :)+bA(2, :)+cA(3, :)$ are both zero.
\end{proof}

Now suppose that according to the local geometry the $l$-th and $m$-th elements of the array:
\begin{equation}
\label{derivs}
\left(\frac{\pa \mathbf{u}^+}{\pa x}, \frac{\pa \mathbf{u}^-}{\pa x}, \frac{\pa \mathbf{u}^+}{\pa y},
\frac{\pa \mathbf{u}^-}{\pa y}, \frac{\pa \mathbf{u}^+}{\pa z}, \frac{\pa \mathbf{u}^-}{\pa z}\right),
\end{equation}
are to be eliminated, where $1\leq l, m \leq 6, l\neq m$. We are going to seek the combined interface conditions
for computing the two pairs of fictitious values at the two irregular grid points.
\\

First, if $l=5, m=6$ or $l=6, m=5$ then we simply employ the interface conditions Eqs. (\ref{jump_val}) and (\ref{jump1}) for computing the fictitious values.
Otherwise, we have the following results.

\begin{lemma}
For given $1\leq l, m\leq 6$ $l\neq 5$ or $6$, or, $m\neq 5$ or $6$, then
there exists constants $a_i, b_i, c_i, d_i, e_i, f_i, g_i, i=1, 2, 3$ such that the $l$-th, $m$-th $(l+6)$-th $(m+6)$-th $(l+12)$-th and $(m+12)$-th elements of the following vectors are all zero:
$$
a_1C(1, :)+b_1C(4, :)+c_1C(7, :)+d_1C(5, :)+e_1C(8, :)+f_1C(6, :)+g_1C(9, :),
$$
$$
a_2C(2, :)+b_2C(4, :)+c_2C(7, :)+d_2C(5, :)+e_2C(8, :)+f_2C(6, :)+g_2C(9, :),
$$
$$
a_3C(3, :)+b_3C(4, :)+c_3C(7, :)+d_3C(5, :)+e_3C(8, :)+f_3C(6, :)+g_3C(9, :),
$$
where $C(i, :), i=1,...,9$ is the $i$-th column of the above matrix $C$.
\end{lemma}

\begin{proof}
We only show that there exists constants $a_1, b_1, c_1, d_1, e_1, f_1, g_1$ such that the results stated in the lemma are true for the first vector, and the other two are the same.
\\

Consider the following three matrices:
$$
A_1\doteq
\left(
\begin{array}{cccccc}
    M^+P(1, 1) & -M^-P(1, 1) &\mu^+ P(1, 2) &-\mu^- P(1, 2) &\mu^+ P(1, 3) &-\mu^-P(1, 3)\\
    P(2, 1) &-P(2, 1)&P(2, 2) &-P(2, 2) &0 &0\\
    P(3, 1) &-P(3, 1) &P(3, 2) &-P(3, 2) &P(3, 3) &-P(3, 3)
  \end{array}
\right)
$$

$$
A_2\doteq
\left(
\begin{array}{cccccc}
    \mu^+P(1, 2) & -\mu^-P(1, 2) &\lambda^+ P(1, 1) &-\lambda^- P(1, 1) &0 &0\\
    P(2, 1) &-P(2, 1)&P(2, 2) &-P(2, 2) &0 &0\\
    P(3, 1) &-P(3, 1) &P(3, 2) &-P(3, 2) &P(3, 3) &-P(3, 3)
  \end{array}
\right)
$$

$$
A_3\doteq
\left(
\begin{array}{cccccc}
    \mu^+P(1, 3) & -\mu^-P(1, 3) &0 &0 &\lambda^+ P(1, 1) &-\lambda^-P(1, 1)\\
    P(2, 1) &-P(2, 1)&P(2, 2) &-P(2, 2) &0 &0\\
    P(3, 1) &-P(3, 1) &P(3, 2) &-P(3, 2) &P(3, 3) &-P(3, 3)
  \end{array}
\right)
$$
According the previous lemma, let
\begin{eqnarray*}
a_1=A_1(2, l)A_1(3, m)-A_1(3, l)A_1(2, m)=A_2(2, l)A_2(3, m)-A_2(3, l)A_2(2, m)&\\
=A_3(2, l)A_3(3, m)-A_3(3, l)A_3(2, m)=C(4, l)C(7, m)-C(7, l)C(4, m),
\end{eqnarray*}
$$
b_1=A_1(3, l)A_1(1, m)-A_1(1, l)A_1(3, m)=C(7, l)C(1, m)-C(1, l)C(7, m),
$$
$$
c_1=A_1(1, l)A_1(2, m)-A_1(2, l)A_1(1, m)=C(1, l)C(4, m)-C(4, l)C(1, m),
$$
$$
d_1=A_2(3, l)A_2(1, m)-A_2(1, l)A_2(3, m)=C(8, l+6)C(1, m+6)-C(1, l+6)C(8, m+6),
$$
$$
e_1=A_2(1, l)A_2(2, m)-A_2(2, l)A_2(1, m)=C(1, l+6)C(5, m+6)-C(5, l+6)C(1, m+6),
$$
$$
f_1=A_3(3, l)A_3(1, m)-A_3(1, l)A_3(3, m)=C(9, 1+12)C(1, m+12)-C(1, l+12)C(9, m+12),
$$
$$
g_1=A_3(1, l)A_3(2, m)-A_3(2, l)A_3(1, m)=C(1, l+12)C(6, m+12)-C(6, l+12)C(1, m+12),
$$
then we have the $l$-th and $m$-th column of the following vectors are all zero:
$$
a_1A_1(1, :)+b_1A_1(2, :)+c_1A_1(3, :),
$$
$$
a_1A_2(1, :)+d_1A_2(2, :)+e_1A_2(3, :),
$$
$$
a_1A_3(1, :)+f_1A_3(2, :)+g_1A_3(3, :).
$$
Note the relationship of the matrix $C$ and the matrices $A_1, A_2$ and $A_3$, it ends up that the
$l$-th, $m$-th $(l+6)$-th $(m+6)$-th $(l+12)$-th and $(m+12)$-th elements of the following vector are all zero:
$$
a_1C(1, :)+b_1C(4, :)+c_1C(7, :)+d_1C(5, :)+e_1C(8, :)+f_1C(6, :)+g_1C(9, :).
$$
\end{proof}

According to the above lemma, if $l\neq 5$ or $6$, or, $m\neq 5$ or $6$. The following two sets of interface conditions can be employed to compute the fictitious values, which do not contains the $l$-th and $m$-th columns' elements of (\ref{derivs}).

The first set of interface condition is due to the jump of function values, i.e.,
\begin{eqnarray}
\label{inter1}
[u_1]|_\Gamma=b_1,
\end{eqnarray}
\begin{eqnarray}
\label{inter22}
[u_2]|_\Gamma=b_2,
\end{eqnarray}
\begin{eqnarray}
\label{inter3}
[u_3]|_\Gamma=b_3.
\end{eqnarray}
The other set is due to derivatives,
\begin{eqnarray}
\label{inter4}
a_1T_1+b_1[\frac{\pa u_1}{\pa \eta}]+c_1\left[\frac{\pa u_1}{\pa \zeta}\right]+d_1\left[\frac{\pa u_2}{\pa \eta}\right]+e_1\left[\frac{\pa u_2}{\pa \zeta}\right]+
f_1\left[\frac{\pa u_3}{\pa \eta}\right]+g_1\left[\frac{\pa u_3}{\pa \zeta}\right]=\\\nonumber
\left(a_1C(1, :)+b_1C(4, :)+c_1C(7, :)+d_1C(5, :)+e_1C(8, :)+f_1C(6, :)+g_1C(9, :)\right)\cdot \mathbf{\alpha},
\end{eqnarray}

\begin{eqnarray}
\label{inter5}
a_2T_2+b_2\left[\frac{\pa u_1}{\pa \eta}\right]+c_2\left[\frac{\pa u_1}{\pa \zeta}\right]+d_2\left[\frac{\pa u_2}{\pa \eta}\right]+e_2\left[\frac{\pa u_2}{\pa \zeta}\right]+
f_2\left[\frac{\pa u_3}{\pa \eta}\right]+g_2\left[\frac{\pa u_3}{\pa \zeta}\right]=\\\nonumber
\left(a_2C(2, :)+b_2C(4, :)+c_2C(7, :)+d_2C(5, :)+e_2C(8, :)+f_2C(6, :)+g_2C(9, :)\right)\cdot \mathbf{\alpha},
\end{eqnarray}

\begin{eqnarray}
\label{inter6}
a_3T_3+b_3\left[\frac{\pa u_1}{\pa \eta}\right]+c_3\left[\frac{\pa u_1}{\pa \zeta}\right]+d_3\left[\frac{\pa u_2}{\pa \eta}\right]+e_3\left[\frac{\pa u_2}{\pa \zeta}\right]+
f_3\left[\frac{\pa u_3}{\pa \eta}\right]+g_3\left[\frac{\pa u_3}{\pa \zeta}\right]=\\\nonumber
\left(a_3C(3, :)+b_3C(4, :)+c_3C(7, :)+d_3C(5, :)+e_3C(8, :)+f_3C(6, :)+g_3C(9, :)\right)\cdot \mathbf{\alpha},
\end{eqnarray}

where
$a_1=a_2=a_3=C(4, l)C(7, m)-C(7, l)C(4, m)$,

$b_1=C(7, l)C(1, m)-C(1, l)C(7, m)$,

$c_1=C(1, l)C(4, m)-C(4, l)C(1, m)$,

$d_1=C(8, l+6)C(1, m+6)-C(1, l+6)C(8, m+6)$,

$e_1=C(1, l+6)C(5, m+6)-C(5, l+6)C(1, m+6)$,

$f_1=C(9, 1+12)C(1, m+12)-C(1, l+12)C(9, m+12)$,

$g_1=C(1, l+12)C(6, m_12)-C(6, l+12)C(1, m+12)$,

$b_2=C(7, l)C(2, m)-C(2, l)C(7, m)$,

$c_2=C(2, l)C(4, m)-C(4, l)C(2, m)$,

$d_2=C(8, l+6)C(2, m+6)-C(2, l+6)C(8, m+6)$,

$e_2=C(2, l+6)C(5, m+6)-C(5, l+6)C(2, m+6)$,

$f_2=C(9, 1+12)C(2, m+12)-C(2, l+12)C(9, m+12)$,

$g_2=C(2, l+12)C(6, m+12)-C(6, l+12)C(2, m+12)$,

$b_3=C(7, l)C(3, m)-C(3, l)C(7, m)$,

$c_3=C(3, l)C(4, m)-C(4, l)C(3, m)$,

$d_3=C(8, l+6)C(3, m+6)-C(3, l+6)C(8, m+6)$,

$e_3=C(3, l+6)C(5, m+6)-C(5, l+6)C(3, m+6)$,

$f_3=C(9, 1+12)C(3, m+12)-C(3, l+12)C(9, m+12)$,

$g_3=C(3, l+12)C(6, m+12)-C(6, l+12)C(3, m+12)$.

In the following, we omit the discussion for the case that $l=5, m=6$ or $l=6, m=5$, which is essential the same as the other cases.

\subsubsection{General fictitious scheme}

\begin{figure}[!ht]
\centering
\includegraphics[width=8cm,height=8cm]{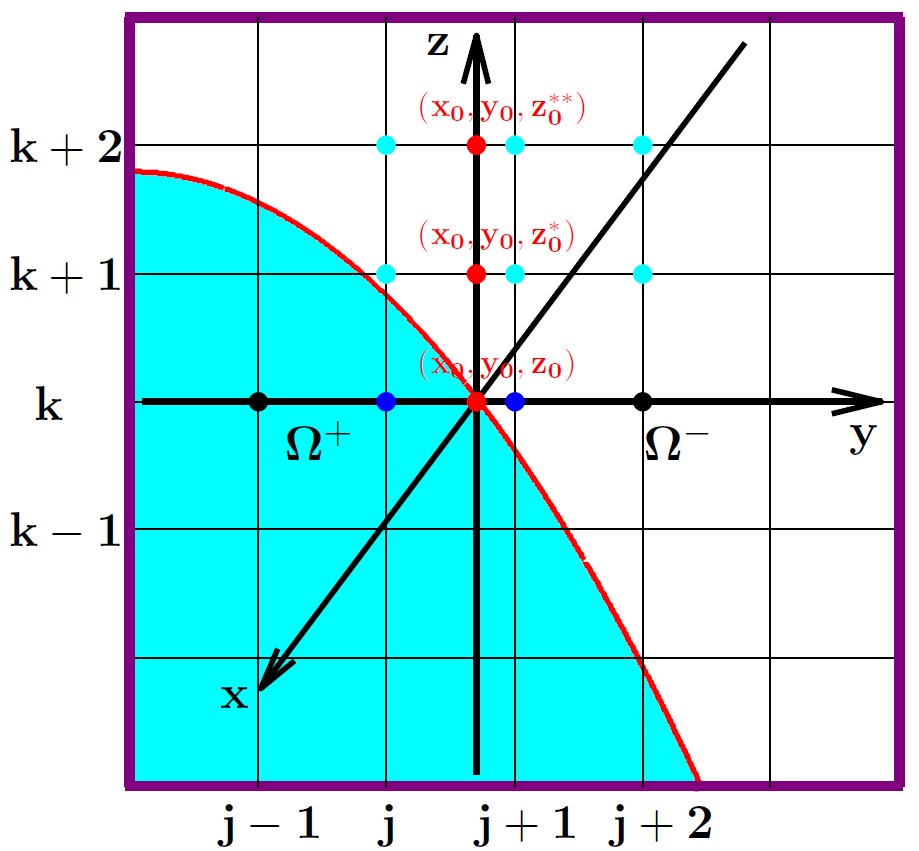}
\caption{Illustration of a smooth interface at cross section $(x=x_i)$. The $k$-th mesh line along the $y$-direction intersects the interface at point $(x_0, y_0, z_0)$. A pair of fictitious values at irregular grid points $(i, j, k)$ and $(i, j+1, k)$ in blue color is to be determined. To this end function values and derivatives at $(x_0, y_0, z_0)$ approximated from $\Omega^+$ and $\Omega^-$ are to be matched. Four points along the $y$-direction are utilized to approximate quantities at $(x_0, y_0, z_0)$, and six cyan points are adopted to approximate $\frac{\partial \mathbf{u}^-(x_0, y_0, z_0)}{\partial z}$.}
\label{3D_inter}
\end{figure}

Consider the geometry illustrated in Fig. \ref{3D_inter}, two sets of fictitious values  $\mathbf{f}(i, j, k):=(f_1^c(i, j, k), f_2^c(i, j, k), f_3^c(i, j, k))^T$ and
$\mathbf{f}(i+1, j, k):=(f_1^c(i+1, j, k), f_2^c(i+1, j, k), f_3^c(i+1, j, k))^T$ are to be determined on the irregular grid points $(i, j, k)$ and $(i+1, j, k)$  for  discretizing central derivatives.

We denote the left domain as $\Omega^+$ and the right one as $\Omega^-$, in this case, $\mathbf{u}^+$, $\mathbf{u}^-$,
$\mathbf{u}_y^+$, and  $\mathbf{u}_y^-$ at $(x_0, y_0, z_0)$ can be easily approximated through interpolations and standard finite difference (FD) schemes from information in $\Omega^+$ and $\Omega^-$, respectively:
\begin{equation}
\mathbf{u}^+=(\omega_{0, j-1}, \omega_{0, j}, \omega_{0, j+1})\cdot\left(\mathbf{u}(i, j-1, k), \mathbf{u}(i, j, k), \mathbf{f}(i, j+1, k)\right)^T,
\end{equation}
\begin{equation}
\mathbf{u}^-=(\tilde{\omega}_{0, j}, \tilde{\omega}_{0, j+1}, \tilde{\omega}_{0, j+2})\cdot\left(\mathbf{f}(i, j, k), \mathbf{u}(i, j+1, k), \mathbf{u}(i, j+2, k)\right)^T,
\end{equation}
\begin{equation}
\mathbf{u}_y^+=(\omega_{1, j-1}, \omega_{1, j}, \omega_{1, j+1})\cdot\left(\mathbf{u}(i, j-1, k), \mathbf{u}(i, j, k), \mathbf{f}(i, j+1, k)\right)^T,
\end{equation}
\begin{equation}
\mathbf{u}_y^-=(\tilde{\omega}_{1, j}, \tilde{\omega}_{1, j+1}, \tilde{\omega}_{1, j+2})\cdot\left(\mathbf{f}(i, j, k), \mathbf{u}(i, j+1, k), \mathbf{u}(i, j+2, k)\right)^T,
\end{equation}
where $\omega_{m, n}$, $\tilde{\omega}_{m, n}$ denote the interpolation or finite difference weights, the first subscript $n$ represents either the interpolation ($n=0$) or the first order derivatives ($n=1$) at interface point $(x_0, y_0, z_0)$, while the second subscript denotes the node index. All the coefficients/weights are generated from  standard Lagrange polynomials \cite{fonberg:1998}.

We only need to compute two of the remaining four vector valued interface quantities. If $\mathbf{u}_x^-$ and $\mathbf{u}_z^-$ can be conveniently computed, then $\mathbf{u}_x^+$ and $\mathbf{u}_z^+$ are eliminated by using the above elimination process with setting $l=1$ and $m=5$.

Here we provide a detailed scheme to approximate $\frac{\partial \mathbf{u}^-}{\partial z}$. Other derivatives can be approximated in the same manner. Without the loss of generality, we only demonstrate how to approximate the first component $\frac{\partial u_1^-}{\partial z}$.

As  shown in Fig. \ref{3D_inter}, to approximate $\frac{\partial u_1^-}{\partial z}$, we need $u_1$ values along the auxiliary line $y=y_0$ on the $yz$-plane. However, these values are unavailable on the grid and have to be approximated by the interpolation schemes along the $y$-direction. Therefore six more auxiliary grid points are involved. In this situation, $\frac{\partial u_1^-}{\partial z}|_{(x_0, y_0, z_0)}$ can be approximated as
\begin{eqnarray}
\frac{\partial u_1^-}{\partial z}=
\left(
\begin{array}{c}
 w_{1, k} \\
 w_{1, k+1}  \\
 w_{1, k+2}
\end{array}
\right)^T\cdot
\left(
\begin{array}{ccccccccc}
    \omega_{0, j} &\omega_{0, j+1}& \omega_{0, j+2} &0 &0 & 0 &0 &0 &0\\
    0 &0& 0 &\omega_{0, j}^{'} &\omega_{0, j+1}^{'} & \omega_{0, j+2}^{'} &0 &0 &0\\
    0 &0& 0 &0 &0 & 0 &\omega_{0, j}^* &\omega_{0, j+1}^* &\omega_{0, j+2}^*
\end{array}
\right)
\cdot \mathbf{U}.
\end{eqnarray}
Here $\mathbf{U}=(f_1^c(i, j, k), u_1(i, j+1, k), u_1(i, j+2, k), u_1(i, j, k+1), u_1(i, j+1, k+1),u_1(i, j+2, k+1), u_1(i, j, k+2), u_1(i, j+1, k+2), u_1(i, j+2, k+2))^T$.
%\begin{eqnarray}
%\tiny{\mathbf{U}=\left(
%\begin{array}{ccccccccc}
%    f_1^c(i, j, k) &u_1(i, j+1, k)& u_1(i, j+2, k) &u_1(i, j, k+1) &u_1(i, j+1, k+1) &u_1(i, j+2, k+1) &u_1(i, j, k+2)& u_1(i, j+1, k+2) &u_1(i, j+2, k+2)
%\end{array}
%\right)}.
%\end{eqnarray}
%where
%$$
%T=\left(f_1^c(i, j, k), u_1(i, j+1, k), u_1(i, j+2, k), u_1(i, j, k+1), u_1(i, j+1, k+1),\\
%u_1(i, j+2, k+1), u_1(i, j, k+2), u_1(i, j+1, k+2), u_1(i, j+2, k+2) \right)
%$$
By solving the above six interface conditions Eq.(\ref{inter1})-(\ref{inter6}) together, six fictitious values  $\mathbf{f}(i, j, k)$ and $\mathbf{f}(i, j+1, k)$ can be easily represented in terms of 48 function values and 12 interface jump conditions around them.

\subsubsection{Matrix optimization}
The MIB matrix is banded due to the reason that the interfaces are 2D surfaces and typically there is only one fictitious value on each side of the interface in a second order MIB scheme. However, to determine each pair of fictitious values, 12 auxiliary grid points are involved and their distribution affects  the convergence property of the resulting MIB matrix.   In most cases, the choice of these 12 auxiliary grid points is not unique. In general, it is very important to make the MIB matrix optimally symmetric and diagonally dominated so as to accelerate the speed of the convergence of the resulting linear algebraic solver. This aspect becomes more important in elasticity interface problems than in elliptic interface problems because the matrix size is much larger.  We therefore select 12 auxiliary grid points as close to the interface as possible. This strategy has  been employed in our earlier MIBPB II software package \cite{Yu:2007,Yu:2007a} for solving  elliptical interface problems.  A more detailed description can be found elsewhere \cite{Yu:2007a}. In the present work, we utilize the same strategy to construct the MIB matrix for elasticity interface problems.

\subsubsection{Fictitious scheme for interface with large curvatures}

The key assumption in the above scheme is that there should be at least two grid points on each mesh line inside a subdomain so that fictitious values on the mesh line can be determined. However, when the curvature of the interface is very large, the above requirement cannot be guaranteed on all mesh sizes.

\begin{figure}[!ht]
\small
\centering
\includegraphics[width=6cm,height=6cm]{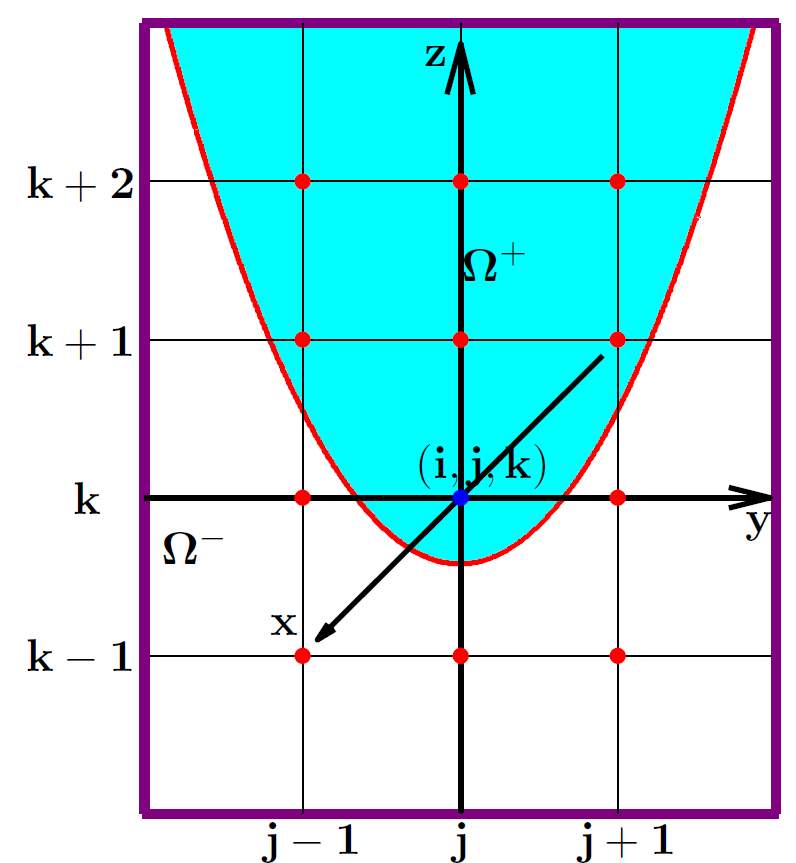}
\caption{An illustration of the disassociation type of irregular grid points  at cross section $(x=x_i)$. Fictitious values $\mathbf{f}(i, j, k)$ cannot be computed from $y$-direction by the aforementioned scheme. Nevertheless, they can be computed from the $z$-direction. In the discretization schemes, the fictitious value at $(i, j, k)$ found from the vertical direction is utilized for both vertical and horizontal discretizations of
the derivatives in the governing equation.
}
\label{dis}
\end{figure}

As shown in  Fig. \ref{dis}, the above scheme is not applicable for finding the fictitious values at grid point $(i, j, k)$ along the
$y$-direction, since there is only one point inside the interface along the $y$-direction, which is not enough for the interpolation.
Nevertheless, there is no problem to find the fictitious values at  grid point $(i, j, k)$ along the $z$-direction. Hence, it is possible to replace the fictitious values to be found along $y$-direction with the fictitious values found along the $z$-direction or the $x$-direction.

Note that this replacement does not reduce the numerical accuracy in general, since if  fictitious values found along $z$-direction has the numerical accuracy $O(h^m)$ for some integer $m$, this estimate holds for the fictitious values at $(i, j, k)$ no matter how they
are determined, where $h$ is the grid size of the uniform mesh. %Hence, simply replace the fictitious values found along $z$-direction by that along $x$-direction is applicable.

\begin{remark}
In principle, to make the final algebraic linear system as symmetric and  banded as possible, if the fictitious values can be found along the given direction, one should avoid  using the disassociation technique.
\end{remark}

%\textcolor{red}{Discretization formulation}

\subsubsection{Fictitious scheme for interface with sharp edge}\label{sharpinterface}

Geometric singularities, such as tips, cusps and self-interesting surfaces, are ubiquitous in science and engineering problems. Due to the existence of geometric singularities, the schemes proposed above may not work because fictitious values can be determined. Therefore, it is crucial to develop some special schemes for determining fictitious values near geometric singularities.

\begin{figure}
\begin{center}
\begin{tabular}{ccc}
\includegraphics[width=0.5\textwidth]{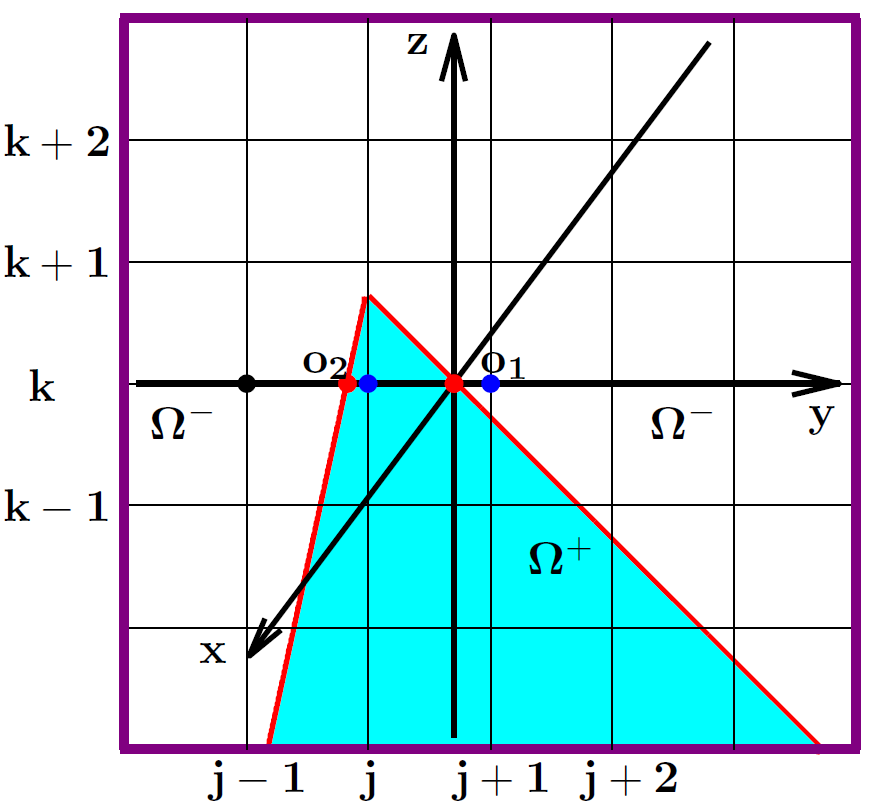}
\includegraphics[width=0.5\textwidth]{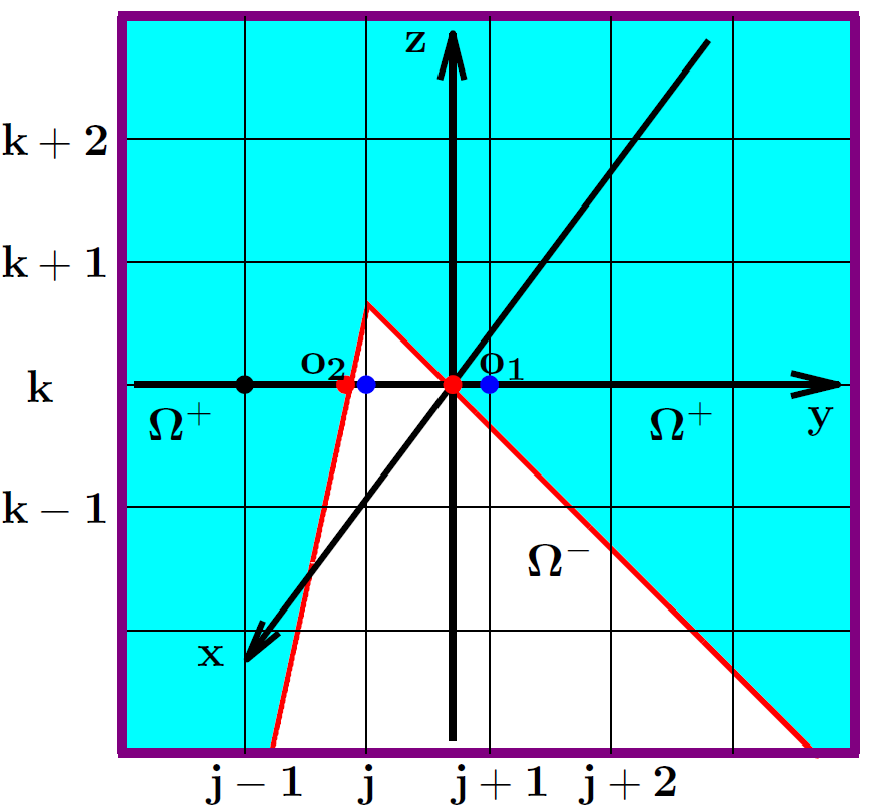}
\end{tabular}
\end{center}
\caption{Illustration of two types of sharp-edged interfaces  at cross section $(x=x_i)$. The jump conditions of the function values at the point $o_1$ and $  o_2$ (two red points) and the jump of derivatives at $o_1$ are employed to compute fictitious values at two blue irregular grid points.
}
\label{sharp}
\end{figure}

According to the local interface geometry, the sharp-edged interface can be classified into two classes, one is locally convex, the other is locally concave, as shown in Fig. \ref{sharp}.

Let us discuss how to determine fictitious values at grid point $(i, j, k)$ for the convex interface case, and the other case can be treated in the same manner.

In the MIB scheme, a pair of fictitious values on a mesh is determined at one time. Suppose that the fictitious values at grid point $(i, j, k)$ is going to be determined along the positive $y$-direction and the interface intersects the mesh line at point $o_1$, the MIB scheme determines fictitious values at $(i, j, k)$ and $(i, j+1, k)$ simultaneously.

In the left chart of Fig. \ref{sharp}, the point $(i, j-1, k)$ will be referred in the discretization of the interface conditions (\ref{inter1})-(\ref{inter6}).  Due to the sharp-edged interface, $(i, j-1, k)$ is not in the same subdomain with $(i, j, k)$, and fictitious values
at grid point $(i, j, k)$ cannot be calculated directly from the interface conditions (\ref{inter1})-(\ref{inter6}). In this case, one more set of fictitious values at grid point $(i, j-1, k)$ will be involved, so that there are nine fictitious values to be determined while there are only
six interface conditions available.

Note that the jump of the function values at point $o_2$, which is another  intersection point of the interface with the mesh line, can be utilized to compute fictitious values. Now there are  nine interface conditions, namely, three jumps of function values at $o_1$, three jumps of function values at $o_2$ and three jumps of derivatives at $o_1$.

The discretization of interface conditions (\ref{inter1})-(\ref{inter6}) in this sharp-edged interface situation can be obtained simply by replacing $\mathbf{u}(i, j-1, k)$ with fictitious values $\mathbf{f}(i, j-1, k)$, where $\mathbf{f}(i, j-1, k):=(f_1^c(i, j-1, k), f_2^c(i, j-1, k), f_3^c(i, j-1, k))^T$ is the fictitious values at node $(i, j-1, k)$. Three more interface conditions at $o_2$ can be discretized as
\begin{equation}
\label{inter7}
[\mathbf{u}]|_{o_2}=\left(\omega'_{0, j-1}, \omega'_{0, j}, \omega'_{0, j+1}\right)\cdot\left((\mathbf{f}(i, j-1, k), \mathbf{u}(i, j, k), \mathbf{u}(i, j+1, k))^T-(\mathbf{u}(i, j-1, k), \mathbf{f}(i, j, k), \mathbf{u}(i, j+1, k))^T\right).
\end{equation}

Fictitious values $\mathbf{f}(i, j-1, k), \mathbf{f}(i, j, k)$ and $\mathbf{f}(i, j+1, k)$ can be calculated from the modified discretization of interface conditions (\ref{inter1})-(\ref{inter6}) and Eqs. (\ref{inter7}).

%\textcolor{red}{Discretization formulation}
\subsubsection{Second Order MIB Finite Difference for Central Derivatives}
All the fictitious values referred in the MIB discretization of the central derivatives can be obtained by the above schemes. At any grid point the second order MIB method applies to all the central derivatives referred in the governing equations of the elasticity interface problem. At an irregular grid point if the CFD scheme refers to some grid points in the other side of the interface, the MIB scheme simply replace the function values at that point by its fictitious values. For instance, the second order MIB finite difference for $\frac{\partial^2 u_1}{\partial y^2}$ at grid point $(i, j, k)$ and $(i, j+1, k)$ in the left chart of Fig. \ref{sharp} are given, respectively, by:
\begin{eqnarray*}
\frac{\partial^2 u_1}{\partial y^2}(i, j, k)=\frac{1}{h^2}\left(f^c_1(i, j-1, k)-2u_1(i, j, k)+f^c_1(i, j+1, k)\right),
\end{eqnarray*}
and
\begin{eqnarray*}
\frac{\partial^2 u_1}{\partial y^2}(i, j+1, k)=\frac{1}{h^2}\left(f^c_1(i, j, k)-2u_1(i, j+1, k)+u_1(i, j+2, k)\right).
\end{eqnarray*}

\subsection{General MIB algorithms for cross derivatives}

The cross derivatives in the elasticity equations make the second order CFD scheme more complicated as the points referred in the CFD schemes are restricted not only to the nearest neighbor points, but also the next nearest neighbor points. This situation does not occur to the elliptic interface problems.

A  critical idea of the MIB method is to reduce  high dimensional problems to locally lower dimensional problems. As such in determining fictitious values for the elliptical interface problems, the MIB scheme carries out  1D-like extensions, which  makes the MIB highly efficient for versatile interface geometries and geometric singularities. Similar idea is applied in the present elasticity interface problem in  determining  fictitious values for both central derivatives and cross derivatives. Based on local interface geometric information, different schemes are designed, including, disassociation type, extrapolation type and neighbor combination type.

\subsubsection{Disassociation scheme}

First we define the disassociation type of fictitious values.
\begin{definition}
An irregular grid point associated with  cross derivatives is called a disassociation type provided that the irregular grid point
is also an irregular grid point associated with  central derivatives.
\end{definition}

The fictitious values on the disassociation type of irregular grid points for cross derivatives can be replaced by fictitious values
found for the central derivatives. Their order of approximation was analyzed in an earlier  paper \cite{Zhou:2006d}.

As illustrated in Fig. \ref{dis}, grid point $(i, j, k)$ is not only irregular in central derivatives, but also irregular in cross derivatives. In this case, fictitious values for the central derivatives at grid point $(i, j, k)$ are obtained based on the numerical scheme proposed for central derivatives. %  In the current scenario, the fictitious values at the grid point $(i, j, k)$ for cross derivatives is simply replaced with that for central derivatives.

%\paragraph{Discretization formulation}

\subsubsection{Extrapolation scheme}

\begin{figure}
\begin{center}
\begin{tabular}{ccc}
\includegraphics[width=0.333\textwidth]{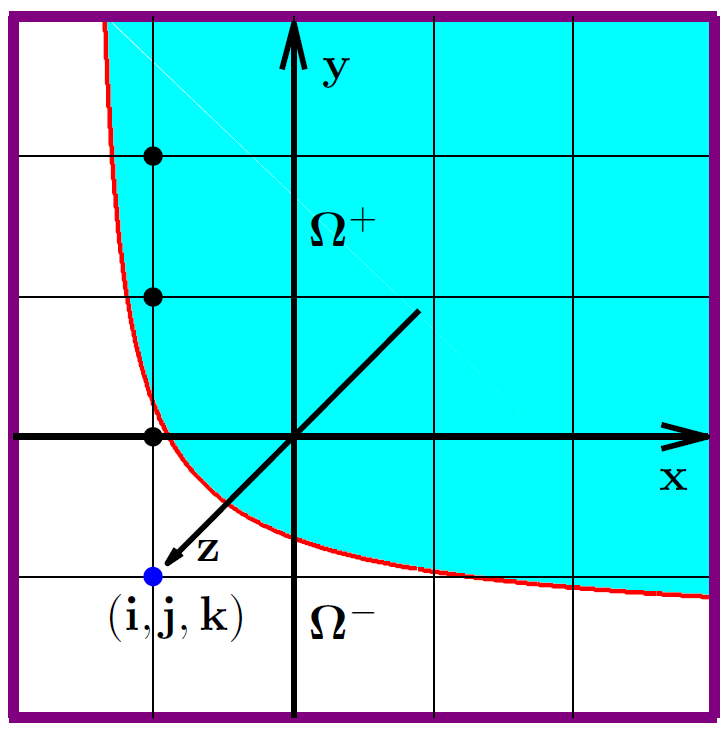}
\includegraphics[width=0.333\textwidth]{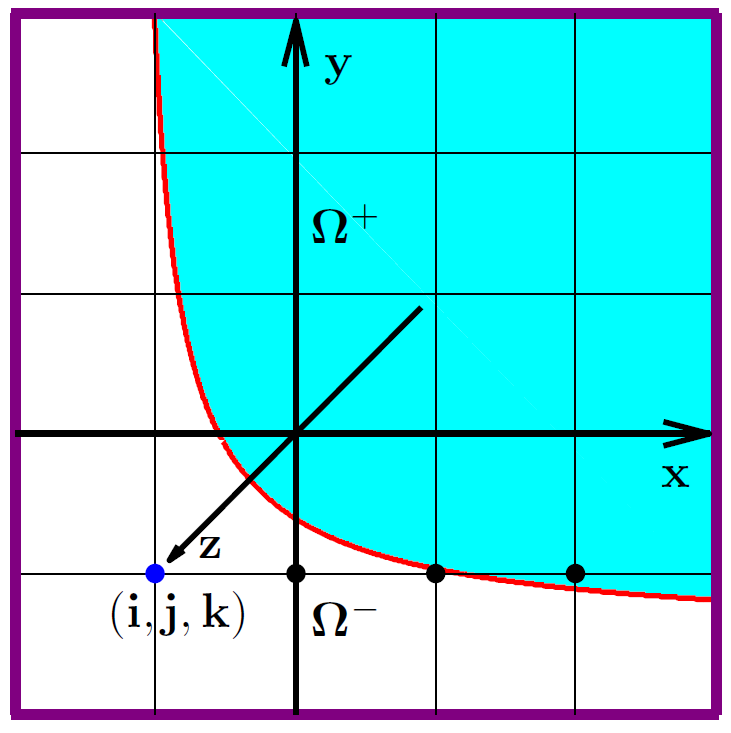}
\includegraphics[width=0.333\textwidth]{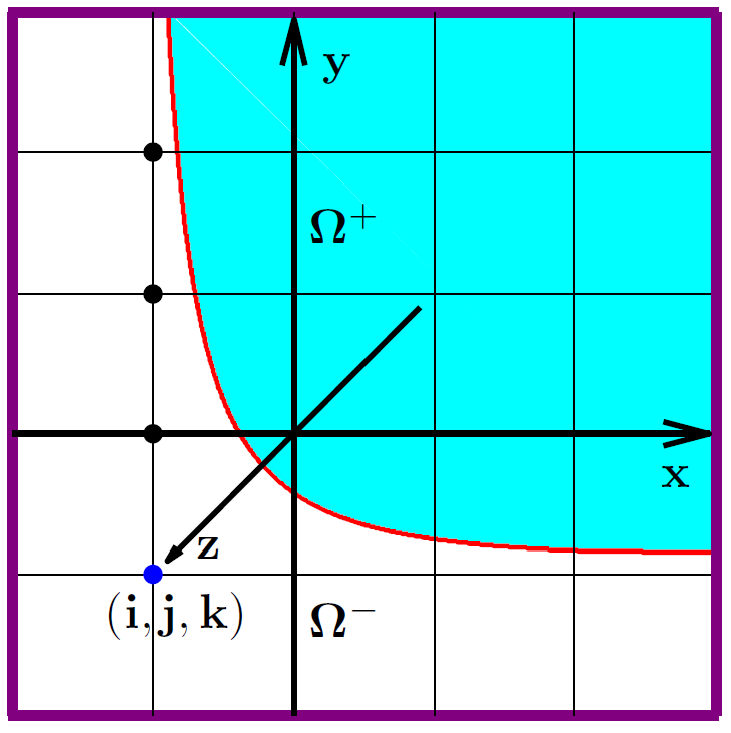}
\end{tabular}
\end{center}
\caption{Illustration of extrapolation type of irregular grid points used in cross derivatives  at cross section $(z=z_k)$. In the left case, fictitious value at the bottom red point and function values at other two red points are employed to approximate fictitious value at $(i, j, k)$.  For the middle case, function values at the right most red point and fictitious values at other two red points are utilized to extrapolate fictitious value at $(i, j, k)$. For the right case, fictitious values at three red points are used to approximate fictitious value at $(i, j, k)$.
}
\label{ext1}
\end{figure}

If a grid point is irregular in the CFD scheme of the cross derivatives while regular for that of the central derivatives, the aforementioned disassociation technique fails. Further, if there exists a direction along which three values are available (function value or fictitious value), then the
extrapolation method is applied.   Suppose that we are seeking the fictitious values at grid point $(i, j, k)$ and project the problem into $xy$-plane, according to the local geometry, the MIB scheme can be classified into three cases.
\begin{itemize}
\item Scheme I. Two function values and one fictitious value are used for the extrapolation. Function values at grid point
$(i, j+2, k)$ and $(i, j+3, k)$, fictitious values at $(i, j+1, k)$ are available and used to extrapolate fictitious values for the cross
derivatives at grid point $(i, j, k)$, see the left chart of Fig. \ref{ext1}.

\item Scheme II. One function value and two fictitious values are used for the extrapolation. Function values at grid point
$(i+3, j, k)$, fictitious values at $(i+1, j, k)$ and $(i+2, j, k)$ are available and used to extrapolate fictitious values for the cross
derivatives at grid point $(i, j, k)$, see the middle  chart of  Fig. \ref{ext1}.

\item Scheme III. Three fictitious values at grid point $(i, j+1, k)$, $(i, j+2, k)$ and $(i, j+3, k)$ are applied to extrapolate the
fictitious value for the cross derivative at the grid point, see the right chart of Fig. \ref{ext1}.
\end{itemize}

%\textcolor{red}{Discretization formulation}

%\subsubsection{Combination scheme}
Now we consider a very special case, in which fictitious values for cross derivatives cannot be obtained with the above schemes.

\begin{figure}[!ht]
\small
\centering
\includegraphics[width=9cm,height=8cm]{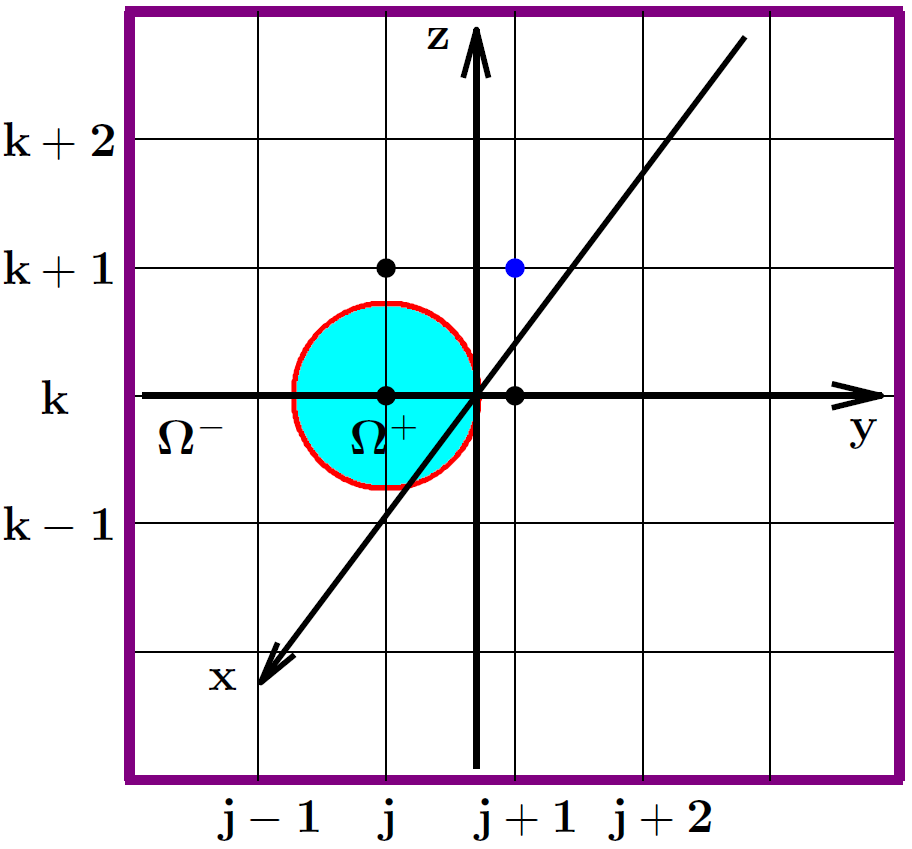}
\caption{Illustration of a single  point situation  at cross section $(x=x_i)$. All the nearest and next nearest neighbor grid points are referred in the second order CFD scheme at grid point $(i, j, k)$. First, fictitious values at its two nearest neighbor grid point, $(i, j+1, k)$ and $(i, j, k+1)$   can be determined by the fictitious scheme for  sharp-edged interfaces. Second, by the neighbor combination scheme, the fictitious values at the blue point  $(i, j+1, k+1)$ can be approximated by the function or fictitious values at three black points.}
\label{1pt}
\end{figure}

As illustrated in Fig. \ref{1pt}, the interface is a sphere centered at $(i, j, k)$ with radius less than grid size $h$.
The CFD scheme at  grid point $(i, j, k)$ refers to all its distance one and two neighbor points. Note that all these points are in the other
side of the interface. To attain a convergent discretization at grid point $(i, j, k)$, all the fictitious values at its neighbor points should be found. Here without the loss of generality, let us only consider the way to find fictitious values for $u_1$ at grid points $(i, j+1, k), (i, j, k+1)$ and $(i, j+1, k+1)$. Due to the symmetry, other fictitious values can be obtained in the same manner.

First, the fictitious values at  grid points $(i, j+1, k)$ and $(i, j, k+1)$ can be found by the sharp interface scheme for central derivatives. Denote the obtained fictitious values to be $f_1^c(i, j+1, k)$ and $f_1^c(i, j, k+1)$, respectively. Further let the analytic extension of the exact solution at these grid points to be $\hat{u}_1(i, j+1, k)$ and $\hat{u}_1(i, j, k+1)$. The numerical extension based on the above MIB scheme satisfies: $f_1^c(i, j+1, k)=\hat{u}_1(i, j+1, k)+O(h^3)$ and $f_1^c(i, j, k+1)=\hat{u}_1(i, j, k+1)+O(h^3)$.
\\

Now the only fictitious value  to be determined is $f_1^c(i, j+1, k+1)$. Based on the Taylor expansion and the above MIB extension estimates,  following equations hold for the uniform Cartesian mesh with grid size $h$.
$$
u_1(i, j+1, k+1)=u_1(i, j, k)+\frac{\partial u_1 }{\partial y}(i, j, k)h+\frac{\partial u_1 }{\partial z}(i, j, k)h+O(h^2)
$$

$$
u_1(i, j+1, k)=u_1(i, j, k)+\frac{\partial u_1 }{\partial y}(i, j, k)h+O(h^2)
$$

$$
u_1(i, j, k+1)=u_1(i, j, k)+\frac{\partial u_1 }{\partial z}(i, j, k)h+O(h^2)
$$

$$
f_1^c(i, j+1, k)=\hat{u}_1(i, j+1, k)+O(h^3),
$$
$$
f_1^c(i, j, k+1)=\hat{u}_1(i, j, k+1)+O(h^3),
$$
where $h$ is the size of the Cartesian mesh.

Therefore, let  fictitious value at grid point $(i, j+1, k+1)$ to be:
$$
f_1^c(i, j+1, k+1)=f_1^c(i, j+1, k)+f_1^c(i, j, k+1)-u_1 (i, j, k).
$$
By direct calculation, the following estimate holds
$$
f_1^c(i, j+1, k+1)=\hat{u}_1(i, j+1, k+1)+O(h^2),
$$
where $\hat{u}_1(i, j+1, k+1)$ is the analytic extension of the exact solution at grid point $(i, j+1, k+1)$.

\begin{remark}
The proposed scheme for finding fictitious values at $(i, j+1, k+1)$ may reduce the numerical accuracy, while based on  numerous numerical
tests, the proposed scheme is still of second order convergence globally.
\end{remark}

%\textcolor{red}{Discretization formulation}
\subsubsection{Second Order MIB Finite Difference for Cross Derivatives}
It is obviously that all the fictitious values at irregular grid points are guaranteed to be found by the above extension and combination schemes. The local combination scheme may lead to some numerical accuracy reduction, however, in most case, this scheme is used quite seldom. Based on our numerous numerical tests examples, the MIB scheme still has the second order numerical accuracy for both $L_\infty$ and $L_2$ error for the elasticity interface problem.

Similar to the MIB discretization of the central derivatives, in the discretization of cross derivatives, when grid point from the other subdomain referred, fictitious values at that point are adopted to replace the function values in the CFD discretization. For instance, the MIB discretization of the $\frac{\partial^2 u_1}{\partial y\partial z}$ at grid point $(i, j, k)$ in Fig. \ref{1pt} is given by:

\begin{eqnarray*}
\frac{\partial^2 u_1}{\partial y\partial z}(i, j, k)=\frac{f_1^c(i, j+1, k+1)+f_1^c(i, j-1, k-1)-f_1^c(i, j+1, k-1)-f_1^c(i, j-1, k+1)}{4h^2}.
\end{eqnarray*}

\section{Numerical experiments} \label{validation}

Numerous numerical tests are designed in this section to investigate the accuracy, efficiency and robustness of the proposed MIB method for solving  3D  elasticity interface problems with both smooth and non-smooth material interfaces. We consider a large number of complex geometric shapes, including  sphere, hemisphere, ellipsoid, cylinder, torus, acorn-like, apple-shaped, flower-like, and pentagon-star shapes in our tests.  Both piecewise constant material parameters and position-dependent material parameters are tested in our investigation. Furthermore, problems with small and large contrast in Poisson's ratio and shear modulus across the interface are also examined.

The standard bi-conjugate gradient  method is employed to solve the linear algebraic system generated by the MIB discretization of the governing equation of the elasticity interface problems. Numerical solutions are compared to the designed analytical solution. Both $L_2$ and $L_\infty$ error measurements are employed in examining the accuracy and convergence of the MIB algorithm for 3D elasticity interface problems
$$
L_\infty(u_k):=\max{|u_k(m, n, l)-\hat{u}_k(m, n, l)|}, k=1, 2, 3; \forall m=1, 2,\cdots n_x; \forall n=1, 2,\cdots n_y; \forall l=1, 2,\cdots n_z
$$
and
$$
L_2:=\sqrt{\frac{1}{n_x*n_y*n_z}\sum^{n_x}_{m=1}\sum^{n_y}_{n=1}\sum^{n_z}_{l=1}(u_k(m, n, l)-\hat{u}_k(m, n, l))^2},
$$
where  $u_k$ $\hat{u}_k$ are the numerical and exact solutions, respectively. Here $L_\infty$ is the maximum error over all the grid points in the computational domain.

\subsection{Smooth interface}

\subsubsection{Piecewise constant shear modulus}

In this section, the proposed MIB method is tested for the piecewise constant material parameters associated with  smooth material interfaces. Problems with both large and small contrasts of Poisson's ratio and shear modulus across the interface are considered in our investigation.

%*******************************************************************************************************************************************
\textbf{Example~1.}
In this example, the computational domain is set to  $[-3, 3]\times[-3, 3]\times[-3, 3]$ and the interface is a sphere which
is defined by $x^2+y^2+z^2=4$.  A sphere is the simplest irregular or complex interface in 3D.
The exact solution is designed to be
$$
u_1(x, y)=
\left\{\begin{array}{ll}
x^2+y^2+z^2-4+\cos(x)\cos(y)\cos(z), &\ \  \mbox{in}\ \ \Omega^+,\\
\cos(x)\cos(y)\cos(z), &\ \  \mbox{in}\ \ \Omega^-,
\end{array}\right.
$$
$$
u_2(x, y)=
\left\{\begin{array}{ll}
x^2+y^2+z^2-4+xy+\cos(x)\cos(y)\cos(z), &\ \  \mbox{in}\ \ \Omega^+,\\
xy+\cos(x)\cos(y)\cos(z), &\ \  \mbox{in}\ \ \Omega^-,
\end{array}\right.
$$
and
$$
u_3(x, y)=
\left\{\begin{array}{ll}
x^2+y^2+z^2-4+yz+\cos(x)\cos(y)\cos(z), &\ \  \mbox{in}\ \ \Omega^+,\\
yz+\cos(x)\cos(y)\cos(z), &\ \  \mbox{in}\ \ \Omega^-.
\end{array}\right.
$$

Note that the above solution guarantees the continuity of the solution across the interface. The Dirichlet boundary conditions and interface jump conditions can be derived from the above exact solution. We consider a series of three cases to test the robustness of the proposed MIB method for large contrasts in material parameters across the interface.

\textbf{Case~1.}
First, let the piecewise constant type of Poisson's ratio and shear modulus  to be
$$
\nu=
\left\{\begin{array}{ll}
\nu^+=0.20, &\ \  \mbox{in}\ \ \Omega^+,\\
\nu^-=0.24, &\ \  \mbox{in}\ \ \Omega^-,
\end{array}\right.
$$
and
$$
\mu=
\left\{\begin{array}{ll}
\mu^+=1500000, &\ \  \mbox{in}\ \ \Omega^+,\\
\mu^-=2000000, &\ \  \mbox{in}\ \ \Omega^-.
\end{array}\right.
$$

Table \ref{Ex1_Case1_linf} lists the grid refinement analysis for the $L_\infty$ error of  the Case 1 of Example 1. We obtain a quite robust second order accuracy in the $L_\infty$ error norm.
It is also interesting to examine the convergence in the $L_2$ error  norm as well. Table \ref{Ex1_Case1_2} presents the grid refinement analysis for the $L_2$ error of the Case 1 of Example 1. We again found highly accurate solutions.

\begin{table}
\caption{ The $L_\infty$ errors for the Case 1 of Example 1.}
\label{Ex1_Case1_linf}
\centering
\begin{tabular}{lllllllll}
\hline
%\multicolumn{5}{c}{The grid refinement analysis of $L_\infty$ error for case1 of example1 } \\
\cline{1-7}
$n_x \times n_y \times n_z$       &  $\  L_\infty(u_1) $ & Order              &$ L_\infty(u_2)$ & Order                &$  L_\infty(u_3)$   &Order  \\
\hline
$10\times 10\times10$ &$6.70\times 10^{-2}$     &           &$6.31\times 10^{-2}$      &   &$5.68\times 10^{-2}$     &        \\
$20\times 20\times20$ &$1.39\times 10^{-2}$     &2.27       &$1.36\times 10^{-2}$     &2.21 &$1.31\times 10^{-2}$     & 2.12   \\
$40\times40\times40$ &$2.72\times 10^{-3}$     &2.35       &$2.94\times 10^{-3}$      &2.21  &$2.69\times 10^{-3}$     &2.28  \\
$80\times80\times80$ &$7.58\times 10^{-4}$     &1.84       &$7.28\times 10^{-4}$      &2.01  &$7.17\times 10^{-4}$     &1.91  \\
\hline
\end{tabular}
\end{table}

\begin{table}
\caption{ The $L_2$ errors for the  Case 1 of Example 1.}
\label{Ex1_Case1_2}
\centering
\begin{tabular}{lllllllll}
\hline
%\multicolumn{5}{c}{The grid refinement analysis of $L_2$ error for case1 of example 1} \\
\cline{1-7}
$n_x \times n_y \times n_z$       & $L_2(u_1)$ &Order              &$L_2(u_2)$      &Order                &$L_2(u_3)$   &Order  \\
\hline
$10\times 10\times10$ &$1.30\times 10^{-2}$     &           &$1.29\times 10^{-2}$      &     &$1.30\times 10^{-2}$     &        \\
$20\times 20\times20$ &$3.20\times 10^{-3}$     &2.02        &$3.20\times 10^{-3}$     &2.01  &$3.15\times 10^{-3}$     &2.05        \\
$40\times40\times40$ &$8.34\times 10^{-4}$     &1.94       &$8.39\times 10^{-4}$      &1.93  &$8.27\times 10^{-4}$     &1.93  \\
$80\times80\times80$ &$2.25\times 10^{-4}$     &1.89       &$2.25\times 10^{-4}$      &1.90  &$2.23\times 10^{-4}$     &1.89  \\
\hline
\end{tabular}
%\caption{The grid refinement analysis of $L_2$ error for Example 1 case 1}
\end{table}

\begin{figure}
\begin{center}
\begin{tabular}{ccc}
\includegraphics[width=0.333\textwidth]{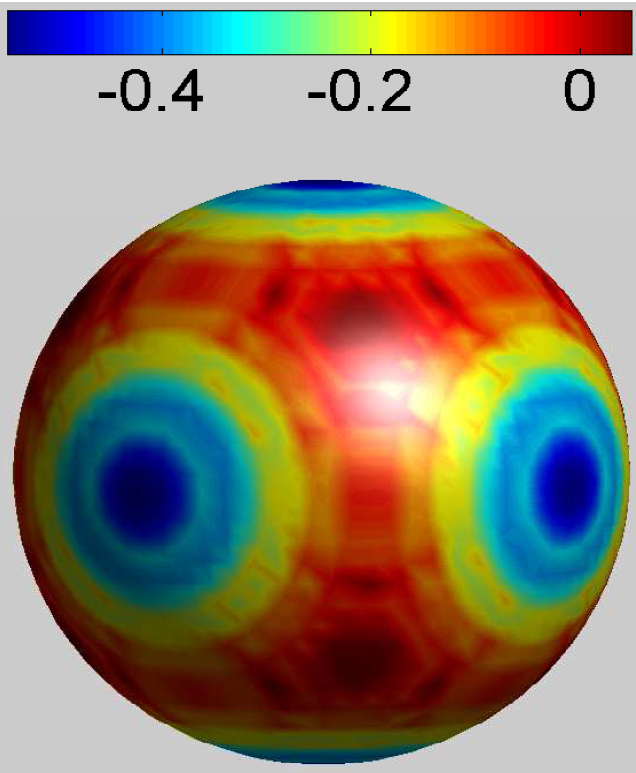}
\includegraphics[width=0.333\textwidth]{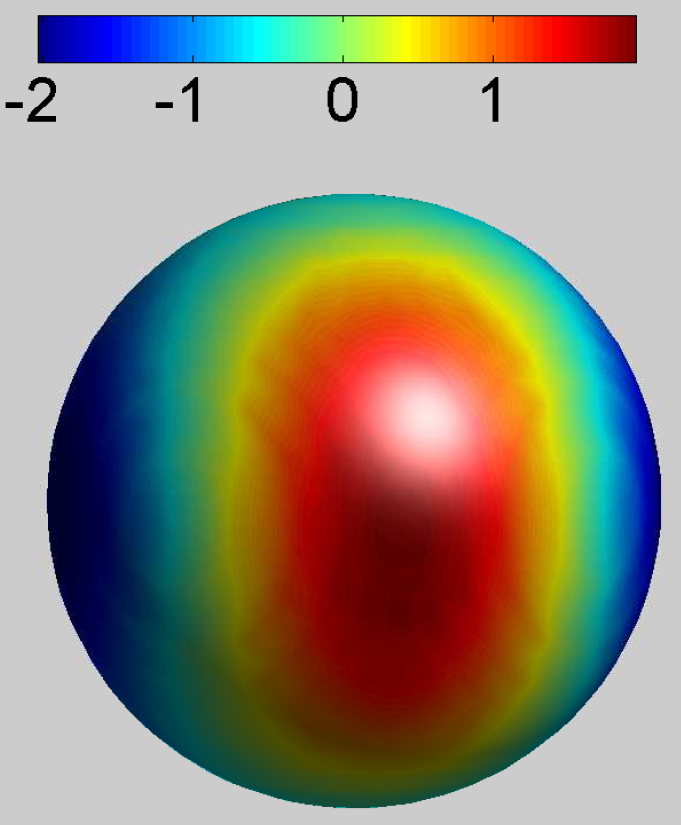}
\includegraphics[width=0.333\textwidth]{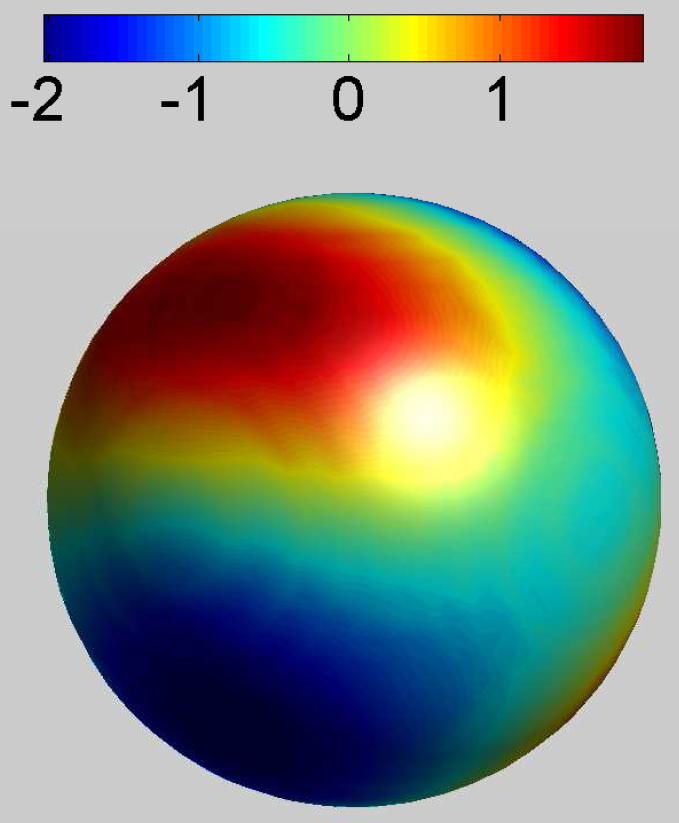}

\end{tabular}
\end{center}
\caption{Numerical solution to the sphere interface problem of Case 1 with 40 grid points along each direction. Left chart: $u_1$; Middle chart $u_2$; Right chart: $u_3$.
}
\label{shpere_sol}
\end{figure}

\begin{figure}
\begin{center}
\begin{tabular}{ccc}
\includegraphics[width=0.333\textwidth]{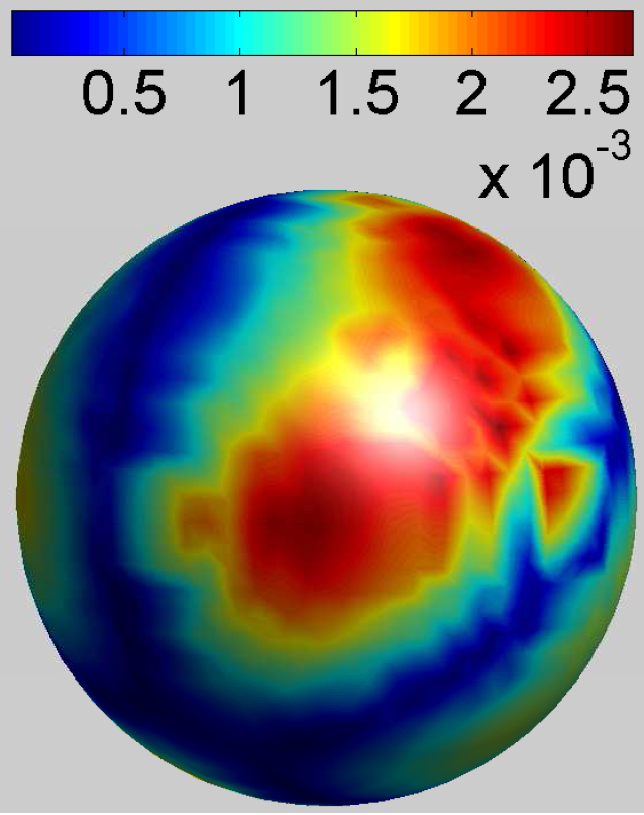}
\includegraphics[width=0.333\textwidth]{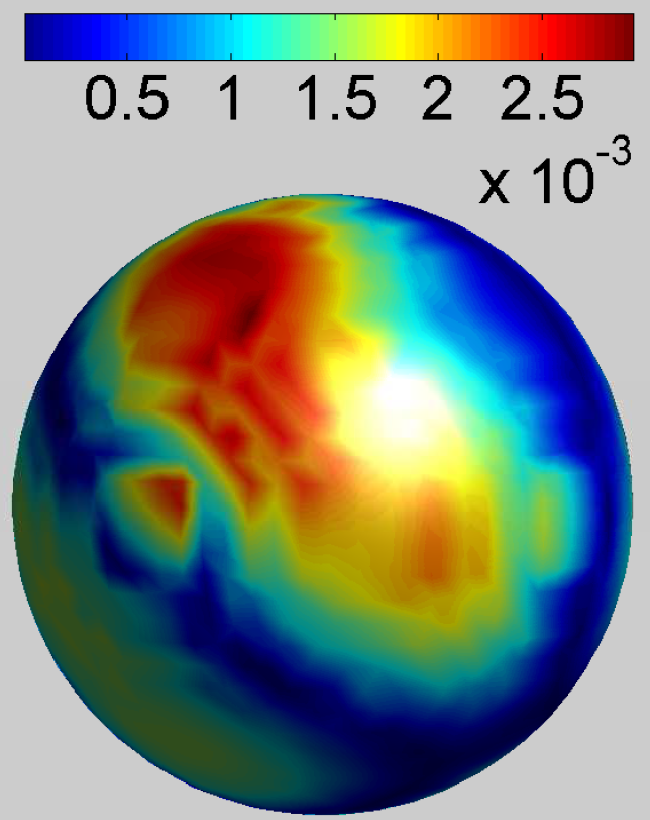}
\includegraphics[width=0.333\textwidth]{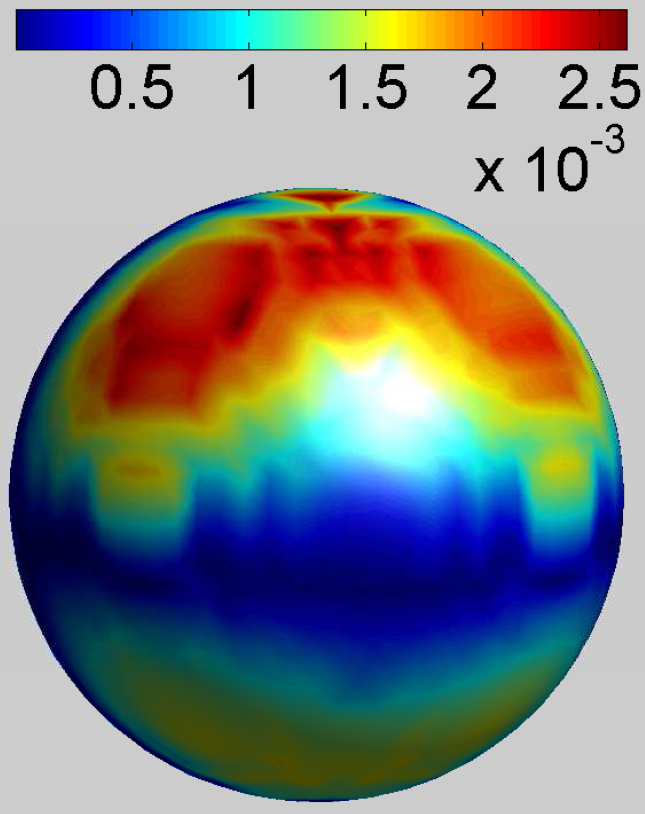}
\end{tabular}
\end{center}
\caption{Numerical error in solving the sphere interface problem of Case 1 with 40 grid points along each direction. Left chart: $u_1$; Middle chart $u_2$; Right chart: $u_3$.
}
\label{shpere_err}
\end{figure}

Figures \ref{shpere_sol} and \ref{shpere_err} illustrate the solution and error with 40  grid points along each direction. Apparently, the errors are quite small.

\textbf{Case~2.}
In this case, we test the proposed MIB method for large contrasts in material parameters across the interface.  We make the Poisson's ratio to be 1000 times in contrast
$$
\nu=
\left\{\begin{array}{ll}
\nu^+=0.00024, &\ \  \mbox{in}\ \ \Omega^+,\\
\nu^-=0.24, &\ \  \mbox{in}\ \ \Omega^-,
\end{array}\right.
$$
while the shear modulus remains unchanged,
$$
\mu=
\left\{\begin{array}{ll}
\mu^+=1500000, &\ \  \mbox{in}\ \ \Omega^+,\\
\mu^-=2000000, &\ \  \mbox{in}\ \ \Omega^-.
\end{array}\right.
$$

Table  \ref{Ex1_Case2_linf} lists the grid refinement analysis for the $L_\infty$ error. Similarly, Table \ref{Ex1_Case2_2} gives the grid refinement analysis for the $L_2$ error  of Case 2. It is seen that both the accuracy and convergence are not affected by the large contrast in the Poisson's ratio across the interface.

\begin{table}
\caption{The $L_\infty$ errors for the Case 2 of Example 1.}
\label{Ex1_Case2_linf}
\centering
\begin{tabular}{lllllllll}
\hline
%\multicolumn{5}{c}{The grid refinement analysis of $L_\infty$ error for case 2 of example1} \\
\cline{1-7}
$n_x \times n_y \times n_z$       &  $\  L_\infty(u_1) $ & Order              &$ L_\infty(u_2)$ & Order                &$  L_\infty(u_3)$   &Order    \\
\hline
$10\times 10\times10$ &$6.21\times 10^{-2}$     &           &$5.95\times 10^{-2}$      &     &$5.45\times 10^{-2}$     &        \\
$20\times 20\times20$ &$1.55\times 10^{-2}$     &2.00       &$1.55\times 10^{-2}$     &1.94  &$1.53\times 10^{-2}$     & 1.83   \\
$40\times40\times40$ &$3.13\times 10^{-3}$     &2.31       &$3.45\times 10^{-3}$      &2.17  &$3.28\times 10^{-3}$     &2.22  \\
$80\times80\times80$ &$8.19\times 10^{-4}$     &1.93       &$7.89\times 10^{-4}$      &2.13  &$7.88\times 10^{-4}$     &2.06  \\
\hline
\end{tabular}
\end{table}

\begin{table}
\caption{The $L_2$ errors for the Case 2 of Example 1.}
\label{Ex1_Case2_2}
\centering
\begin{tabular}{lllllllll}
\hline
%\multicolumn{5}{c}{The grid refinement analysis of $L_2$ error for case2 of example 1} \\
\cline{1-7}
$n_x \times n_y \times n_z$       & $L_2(u_1)$ &Order              &$L_2(u_2)$      &Order                &$L_2(u_3)$   &Order  \\
\hline
$10\times 10\times10$ &$1.29\times 10^{-2}$     &           &$1.28\times 10^{-2}$      &     &$1.28\times 10^{-2}$     &        \\
$20\times 20\times20$ &$3.29\times 10^{-3}$     &1.97        &$3.28\times 10^{-3}$     &1.96  &$3.22\times 10^{-3}$     &1.99        \\
$40\times40\times40$ &$8.49\times 10^{-4}$     &1.95       &$8.53\times 10^{-4}$      &1.94  &$8.41\times 10^{-4}$     &1.94  \\
$80\times80\times80$ &$2.29\times 10^{-4}$     &1.89       &$2.29\times 10^{-4}$      &1.90  &$2.28\times 10^{-4}$     &1.89  \\
\hline
\end{tabular}
\end{table}

\textbf{Case~3.}
Having tested the proposed MIB method for large contrast in the Poisson's ratio, let us enlarge the contrast of the shear modulus across the interface, while the Poisson's ratio is unchanged,
$$
\nu=
\left\{\begin{array}{ll}
\nu^+=0.20, &\ \  \mbox{in}\ \ \Omega^+,\\
\nu^-=0.24, &\ \  \mbox{in}\ \ \Omega^-,
\end{array}\right.
$$
and
$$
\mu=
\left\{\begin{array}{ll}
\mu^+=2000, &\ \  \mbox{in}\ \ \Omega^+,\\
\mu^-=2000000, &\ \  \mbox{in}\ \ \Omega^-.
\end{array}\right.
$$

\begin{table}
\caption{The $L_\infty$ errors for the Case 3 of Example 1.}
\label{Ex13inf}
\centering
\begin{tabular}{lllllllll}
\hline
%\multicolumn{5}{c}{The grid refinement analysis of $L_\infty$ error for case3 of example1} \\
\cline{1-7}
$n_x \times n_y \times n_z$       &  $\  L_\infty(u_1) $ & Order              &$ L_\infty(u_2)$ & Order                &$  L_\infty(u_3)$   &Order    \\
\hline
$10\times 10\times10$ &$6.70\times 10^{-2}$     &           &$6.31\times 10^{-2}$      &     &$5.68\times 10^{-2}$     &        \\
$20\times 20\times20$ &$1.40\times 10^{-2}$     &2.26       &$1.41\times 10^{-2}$     &1.94  &$1.31\times 10^{-2}$     &2.12   \\
$40\times40\times40$ &$2.72\times 10^{-3}$     &2.36       &$2.39\times 10^{-3}$      &2.16  &$2.69\times 10^{-3}$     &2.28  \\
$80\times80\times80$ &$7.58\times 10^{-4}$     &1.85       &$7.28\times 10^{-4}$      &1.72  &$7.17\times 10^{-4}$     &1.91  \\
\hline
\end{tabular}
\end{table}

\begin{table}
\caption{The $L_2$ errors for the Case 3 of Example 1.}
\label{Ex1_Case3_2}
\centering
\begin{tabular}{lllllllll}
\hline
%\multicolumn{5}{c}{The grid refinement analysis of $L_2$ error for case3 of example 1} \\
\cline{1-7}
$n_x \times n_y \times n_z$       & $L_2(u_1)$ &Order              &$L_2(u_2)$      &Order                &$L_2(u_3)$   &Order \\
\hline
$10\times 10\times10$ &$1.30\times 10^{-2}$     &           &$1.29\times 10^{-2}$      &     &$1.30\times 10^{-2}$     &        \\
$20\times 20\times20$ &$3.19\times 10^{-3}$     &2.03       &$3.20\times 10^{-3}$     &2.01  &$3.15\times 10^{-3}$     &2.05        \\
$40\times40\times40$ &$8.34\times 10^{-4}$     &1.94       &$8.39\times 10^{-4}$      &1.94  &$8.27\times 10^{-4}$     &1.93  \\
$80\times80\times80$ &$2.25\times 10^{-4}$     &1.89       &$2.24\times 10^{-4}$      &1.91  &$2.23\times 10^{-4}$     &1.89  \\
\hline
\end{tabular}
\end{table}

Table \ref{Ex13inf} gives the grid refinement analysis for the $L_\infty$ error  of Case 3. In
Table \ref{Ex1_Case3_2}, we  provide the grid refinement analysis for the $L_2$ error  of Case 3. Obviously, the accuracy and convergence are the same as those in Case 1. Therefore, the proposed method is very robust against large contrasts in material parameters. We have  obtained second order accuracy in both $L_\infty$ and $L_2$ error norms in all three cases in Example 1.

%********************************************************************************************************************************
\textbf{Example~2.}
In this example, we modify the interface geometry.  Let the computational domain be $[-3, 3]\times[-3, 3]\times[-3, 3]$ and the interface be given as a hemisphere
$$
\left\{\begin{array}{ll}
x^2+y^2+z^2=4,\\
z\geq 0,
\end{array}\right.
$$

To ensure the continuity of the solution across the interface, the analytic solution adopted in this example is the same as that in Example 1.
In this example, we also test the numerical scheme for three different cases of Poisson's ratio and shear modulus, in each case the
material parameters are inherited from the corresponded case in Example 1.

\textbf{Case~1.}
Table \ref{Ex2_Case1_linf} gives the grid refinement analysis of the $L_\infty$ error. Similarly
the grid refinement analysis of the $L_2$ error is presented in Table \ref{Ex2_case1_2}. It is seen that both the level of accuracy and the order of convergence are the same as those in the Case 1 of Example 1, which suggests that the proposed method is sensitive to the change in the geometry.

\begin{table}
\caption{The $L_\infty$ errors for the Case 1 of Example 2.}
\label{Ex2_Case1_linf}
\centering
\begin{tabular}{lllllllll}
\hline
%\multicolumn{5}{c}{The grid refinement analysis of $L_\infty$ error for example 2 case1} \\
\cline{1-7}
$n_x \times n_y \times n_z$       &  $\  L_\infty(u_1) $ & Order              &$ L_\infty(u_2)$ & Order                &$  L_\infty(u_3)$   &Order    \\
\hline
$10\times 10\times10$ &$6.38\times 10^{-2}$     &           &$5.93\times 10^{-2}$      &   &$6.17\times 10^{-2}$     &        \\
$20\times 20\times20$ &$1.35\times 10^{-2}$     &2.24       &$1.33\times 10^{-2}$     &2.16 &$1.35\times 10^{-2}$     & 2.19   \\
$40\times40\times40$ &$2.67\times 10^{-3}$     &2.34       &$2.97\times 10^{-3}$      &2.16  &$2.70\times 10^{-3}$     &2.32  \\
$80\times80\times80$ &$6.28\times 10^{-4}$     &2.09       &$6.52\times 10^{-4}$      &2.19  &$5.91\times 10^{-4}$     &2.19  \\
\hline
\end{tabular}
\end{table}

\begin{table}
\caption{The $L_2$ errors for the Case 1 of Example 2.}
\label{Ex2_case1_2}
\centering
\begin{tabular}{lllllllll}
\hline
%\multicolumn{5}{c}{The grid refinement analysis of $L_2$ error for example 2 case1} \\
\cline{1-7}
$n_x \times n_y \times n_z$       & $L_2(u_1)$ &Order              &$L_2(u_2)$      &Order                &$L_2(u_3)$   &Order \\
\hline
$10\times 10\times10$ &$1.33\times 10^{-2}$     &           &$1.32\times 10^{-2}$      &     &$1.57\times 10^{-2}$     &        \\
$20\times 20\times20$ &$3.35\times 10^{-3}$     &1.99        &$3.33\times 10^{-3}$     &1.99  &$3.44\times 10^{-3}$     &2.19        \\
$40\times40\times40$ &$8.61\times 10^{-4}$     &1.96       &$8.61\times 10^{-4}$      &1.95  &$8.65\times 10^{-4}$     &1.99  \\
$80\times80\times80$ &$2.01\times 10^{-4}$     &2.10       &$2.02\times 10^{-4}$      &2.09  &$2.01\times 10^{-4}$     &2.11  \\
\hline
\end{tabular}
\end{table}

\begin{figure}
\begin{center}
\begin{tabular}{ccc}
\includegraphics[width=0.333\textwidth]{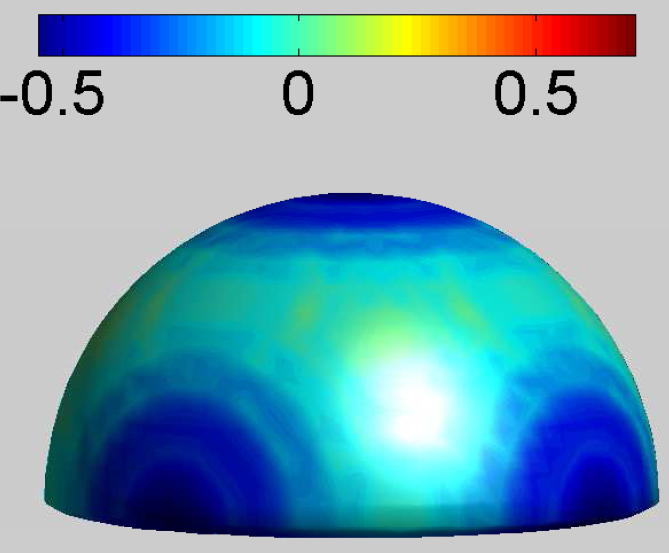}
\includegraphics[width=0.333\textwidth]{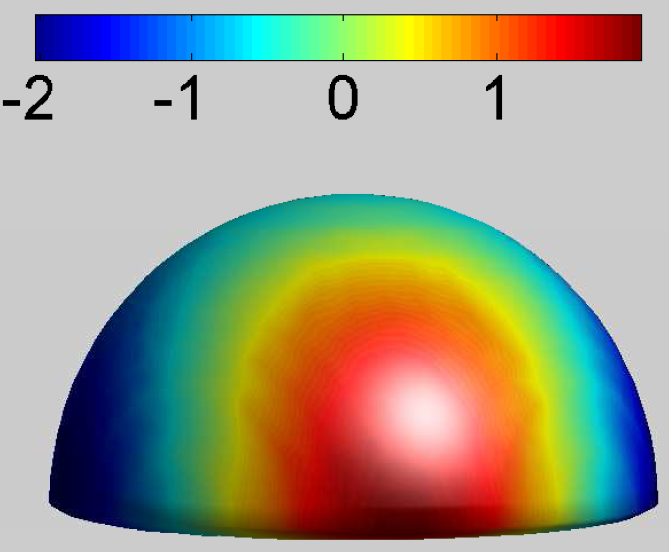}
\includegraphics[width=0.333\textwidth]{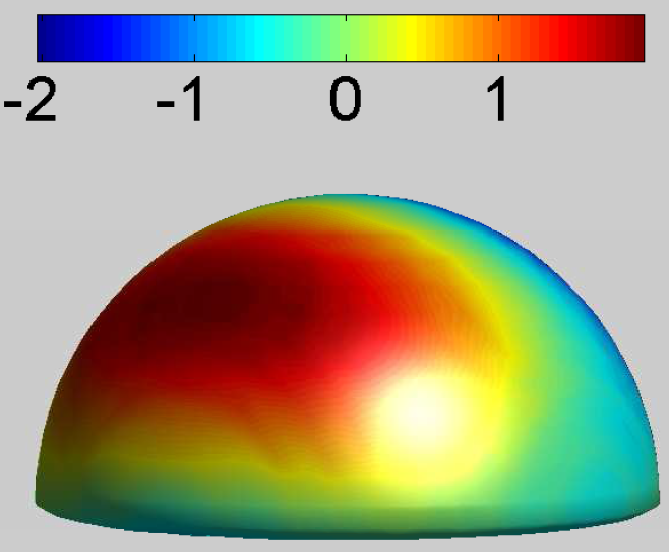}
\end{tabular}
\end{center}
\caption{Numerical solution to the Case 1 of the hemisphere interface problem with 40 grid points along each direction of the computational domain.
Left chart: $u_1$; Middle chart $u_2$; Right chart: $u_3$.
}
\label{hemi_sol}
\end{figure}

\begin{figure}
\begin{center}
\begin{tabular}{ccc}
\includegraphics[width=0.333\textwidth]{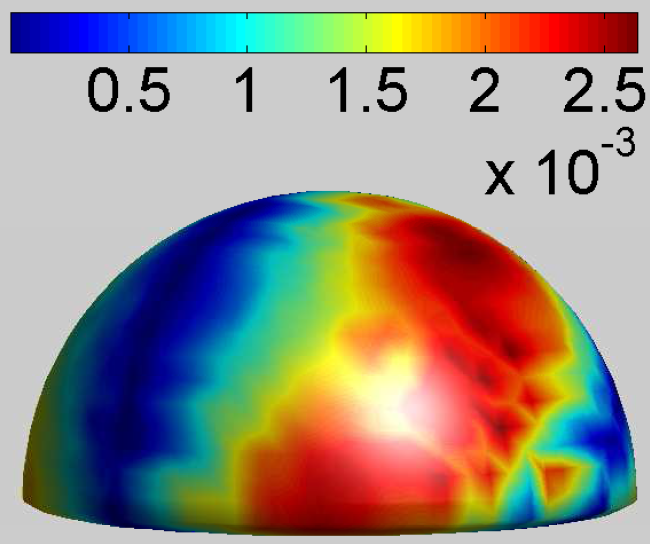}
\includegraphics[width=0.333\textwidth]{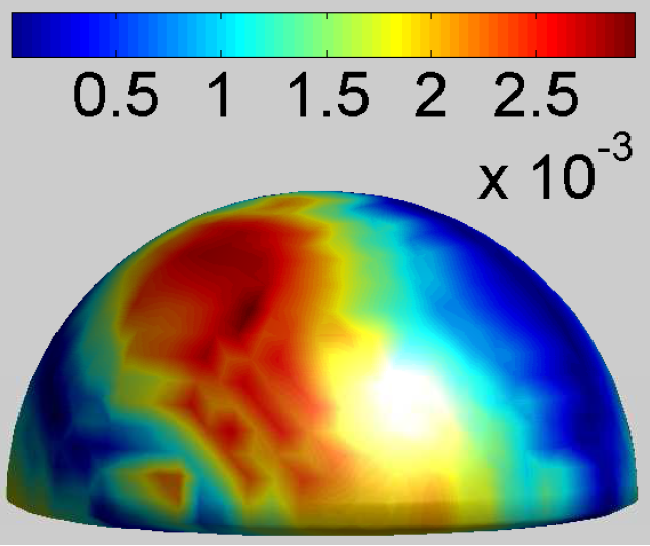}
\includegraphics[width=0.333\textwidth]{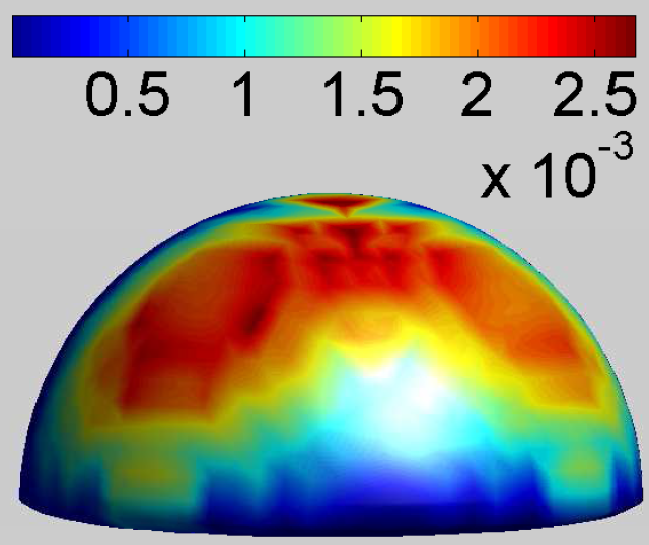}
\end{tabular}
\end{center}
\caption{Numerical error for the solving  the Case 1 of the hemisphere interface problem   with 40 grids along each direction of the computational domain.
Left chart: $u_1$; Middle chart $u_2$; Right chart: $u_3$.
}
\label{hemi_err}
\end{figure}

Figures \ref{hemi_sol} and \ref{hemi_err} show the numerical solution and error of the Case 1 of Example 2, respectively. The number of grids is 40 along each direction of the computational domain.

%\begin{figure}[!ht]
%\small
%\centering
%\includegraphics[width=12cm,height=4.0cm]{sol/hemi_sol_40.png}
%\caption{Solution to the hemisphere interface case 1 with 40 grids along each side of the computational domain, left $u_1$, middle $u_2$, right $u_3$}
%\label{hemi_sol}
%\end{figure}
%
%
%\begin{figure}[!ht]
%\small
%\centering
%\includegraphics[width=12cm,height=4.0cm]{sol/hemi_err_40.png}
%\caption{Error to the hemisphere interface case 1 with 40 grids along each side of the computational domain, left $u_1$, middle $u_2$, right $u_3$}
%\label{hemi_err}
%\end{figure}

\textbf{Case~2.}
Table \ref{Ex2_Case2_linf} gives the grid refinement analysis of the $L_\infty$ error for the large contrast between Poisson's ratio across the interface.
The numerical behavior is quite similar to that in the Case 2 of Example 1.
\begin{table}
\caption{The $L_\infty$ errors for the Case 2 of Example 2.}
\label{Ex2_Case2_linf}
\centering
\begin{tabular}{lllllllll}
\hline
%\multicolumn{5}{c}{The grid refinement analysis of $L_\infty$ error for example 2 case2} \\
\cline{1-7}
$n_x \times n_y \times n_z$       &  $\  L_\infty(u_1) $ & Order              &$ L_\infty(u_2)$ & Order                &$  L_\infty(u_3)$   &Order   \\
\hline
$10\times 10\times10$ &$5.87\times 10^{-2}$     &           &$5.59\times 10^{-2}$      &     &$5.93\times 10^{-2}$     &        \\
$20\times 20\times20$ &$1.48\times 10^{-2}$     &1.99       &$1.50\times 10^{-2}$     &1.90  &$1.58\times 10^{-2}$     & 1.91   \\
$40\times40\times40$ &$3.00\times 10^{-3}$     &2.30       &$3.47\times 10^{-3}$      &2.11  &$3.30\times 10^{-3}$     &2.26  \\
$80\times80\times80$ &$6.80\times 10^{-4}$     &2.14       &$7.00\times 10^{-4}$      &2.31  &$6.72\times 10^{-4}$     &2.30  \\
\hline
\end{tabular}
\end{table}

Table \ref{Ex2_case2_2} lists the grid refinement analysis of the $L_2$ error for the large contrast between Poisson's ratio across the interface.
We observe the second order convergence in the $L_2$ error norm.

\begin{table}
\caption{The $L_2$ errors for the Case 2 of Example 2.}
\label{Ex2_case2_2}
\centering
\begin{tabular}{lllllllll}
\hline
%\multicolumn{5}{c}{The grid refinement analysis of $L_2$ error for example 2 case2} \\
\cline{1-7}
$n_x \times n_y \times n_z$       & $L_2(u_1)$ &Order              &$L_2(u_2)$      &Order                &$L_2(u_3)$   &Order \\
\hline
$10\times 10\times10$ &$1.33\times 10^{-2}$     &           &$1.32\times 10^{-2}$      &     &$1.62\times 10^{-2}$     &        \\
$20\times 20\times20$ &$3.42\times 10^{-3}$     &1.96        &$3.40\times 10^{-3}$     &1.96  &$3.55\times 10^{-3}$     &2.19        \\
$40\times40\times40$ &$8.73\times 10^{-4}$     &1.97       &$8.75\times 10^{-4}$      &1.96  &$8.89\times 10^{-4}$     &2.00  \\
$80\times80\times80$ &$2.02\times 10^{-4}$     &2.11       &$2.02\times 10^{-4}$      &2.11  &$2.13\times 10^{-4}$     &2.06  \\
\hline
\end{tabular}
\end{table}

\textbf{Case~3.}
Table \ref{Ex2_Case3_linf} offers the grid refinement analysis of the $L_\infty$ error for the large contrast between shear modulus across the interface.
Table \ref{Ex2_case3_2} gives the grid refinement analysis of the $L_2$ error. In all these three cases in Example 2, the second order convergence in both $L_\infty$ and $L_2$ errors is essential reached. The level of accuracy is the same as that found in Example 1.

\begin{table}
\caption{The $L_\infty$ errors for Case 3 of Example 2.}
\label{Ex2_Case3_linf}
\centering
\begin{tabular}{lllllllll}
\hline
%\multicolumn{5}{c}{The grid refinement analysis of $L_\infty$ error for example 2 case3} \\
\cline{1-7}
$n_x \times n_y \times n_z$       &  $\  L_\infty(u_1) $ & Order              &$ L_\infty(u_2)$ & Order                &$  L_\infty(u_3)$   &Order   \\
\hline
$10\times 10\times10$ &$6.38\times 10^{-2}$     &           &$5.93\times 10^{-2}$      &     &$6.17\times 10^{-2}$     &        \\
$20\times 20\times20$ &$1.35\times 10^{-2}$     &2.24       &$1.33\times 10^{-2}$     &2.16  &$1.35\times 10^{-2}$     &2.19   \\
$40\times40\times40$ &$2.67\times 10^{-3}$     &2.34       &$2.97\times 10^{-3}$      &2.16  &$2.70\times 10^{-3}$     &2.32  \\
$80\times80\times80$ &$6.28\times 10^{-4}$     &2.09       &$6.52\times 10^{-4}$      &2.19  &$5.91\times 10^{-4}$     &2.19  \\
\hline
\end{tabular}
\end{table}

\begin{table}
\caption{The $L_2$ errors for the Case 3 of Example 2.}
\label{Ex2_case3_2}
\centering
\begin{tabular}{lllllllll}
\hline
%\multicolumn{5}{c}{The grid refinement analysis of $L_2$ error for example 2 case3} \\
\cline{1-7}
$n_x \times n_y \times n_z$       & $L_2(u_1)$ &Order              &$L_2(u_2)$      &Order                &$L_2(u_3)$   &Order \\
\hline
$10\times 10\times10$ &$1.33\times 10^{-2}$     &           &$1.32\times 10^{-2}$      &     &$1.57\times 10^{-2}$     &        \\
$20\times 20\times20$ &$3.35\times 10^{-3}$     &1.99       &$3.33\times 10^{-3}$     &1.99  &$3.44\times 10^{-3}$     &2.19        \\
$40\times40\times40$ &$8.60\times 10^{-4}$     &1.96       &$8.61\times 10^{-4}$      &1.95  &$8.65\times 10^{-4}$     &1.99  \\
$80\times80\times80$ &$2.01\times 10^{-4}$     &2.10       &$2.00\times 10^{-4}$      &2.10  &$2.00\times 10^{-4}$     &2.11  \\
\hline
\end{tabular}
\end{table}

%****************************************************************************************************************
\textbf{Example~3.}
In this example, the computational domain is set to be: $[-3, 3]\times[-4, 4]\times[-2, 2]$ with an
ellipsoid interface defined as $\frac{x^2}{4}+\frac{y^2}{9}+z^2=1$.

The Dirichlet boundary condition and interface jump conditions are determined from the following exact solution
$$
u_1(x, y)=
\left\{\begin{array}{ll}
\frac{x^2}{4}+\frac{y^2}{9}+z^2-1+\cos(x)\cos(y)\cos(z), &\ \  \mbox{in}\ \ \Omega^+,\\
\cos(x)\cos(y)\cos(z), &\ \  \mbox{in}\ \ \Omega^-.
\end{array}\right.
$$
$$
u_2(x, y)=
\left\{\begin{array}{ll}
\frac{x^2}{4}+\frac{y^2}{9}+z^2-1+xy+\cos(x)\cos(y)\cos(z), &\ \  \mbox{in}\ \ \Omega^+,\\
xy+\cos(x)\cos(y)\cos(z), &\ \  \mbox{in}\ \ \Omega^-.
\end{array}\right.
$$
and
$$
u_3(x, y)=
\left\{\begin{array}{ll}
\frac{x^2}{4}+\frac{y^2}{9}+z^2-1+yz+\cos(x)\cos(y)\cos(z), &\ \  \mbox{in}\ \ \Omega^+,\\
yz+\cos(x)\cos(y)\cos(z), &\ \  \mbox{in}\ \ \Omega^-.
\end{array}\right.
$$
Obviously, the property of solution continuity across the interface is also satisfied in the above solution. In this example, three different cases  of the material parameters used in the above two examples are adopted to examine the sensitivity of the proposed MIB method to the change in interface geometry.

\textbf{Case~1.}
Grid refinement  analysis for $L_\infty$ error is demonstrated in Table \ref{Ex3_case1_linf} for the ellipsoid interface. A similar
 analysis for $L_2$ error is listed in Table \ref{Ex3_case1_2} for the ellipsoid interface. Again, we see the same type of behavior in accuracy and convergence as that in last few examples.

\begin{table}
\caption{The $L_\infty$ errors for Case 1 of Example 3.}
\label{Ex3_case1_linf}
\centering
\begin{tabular}{lllllllll}
\hline
%\multicolumn{5}{c}{The grid refinement analysis of $L_\infty$ error for example case1} \\
\cline{1-7}
$n_x \times n_y \times n_z$       &  $\  L_\infty(u_1) $ & Order              &$ L_\infty(u_2)$ & Order                &$  L_\infty(u_3)$   &Order   \\
\hline
$10\times 10\times10$ &$3.96\times 10^{-2}$     &           &$5.23\times 10^{-2}$      &   &$3.16\times 10^{-2}$     &        \\
$20\times 20\times20$ &$1.52\times 10^{-2}$     &1.38       &$1.31\times 10^{-2}$     &2.00 &$8.07\times 10^{-2}$     & 1.97   \\
$40\times40\times40$ &$2.82\times 10^{-3}$     &2.43       &$3.45\times 10^{-3}$      &1.93  &$1.90\times 10^{-3}$     &2.09  \\
$80\times80\times80$ &$6.99\times 10^{-4}$     &2.01       &$8.81\times 10^{-4}$      &1.97  &$4.83\times 10^{-4}$     &1.98  \\
\hline
\end{tabular}
\end{table}

\begin{table}
\caption{The $L_2$ errors for the Case 1 of Example 3.}
\label{Ex3_case1_2}
\centering
\begin{tabular}{lllllllll}
\hline
%\multicolumn{5}{c}{The grid refinement analysis of $L_2$ error for example 3 case 1} \\
\cline{1-7}
$n_x \times n_y \times n_z$       & $L_2(u_1)$ &Order              &$L_2(u_2)$      &Order                &$L_2(u_3)$   &Order \\
\hline
$10\times 10\times10$ &$1.16\times 10^{-2}$     &           &$1.45\times 10^{-2}$      &     &$6.72\times 10^{-2}$     &        \\
$20\times 20\times20$ &$3.54\times 10^{-3}$     &1.71        &$3.96\times 10^{-3}$     &1.88  &$1.86\times 10^{-3}$     &1.85        \\
$40\times40\times40$ &$8.61\times 10^{-4}$     &2.04       &$1.08\times 10^{-4}$      &1.88  &$4.78\times 10^{-4}$     &1.96  \\
$80\times80\times80$ &$2.23\times 10^{-4}$     &1.95       &$2.86\times 10^{-4}$      &1.92  &$1.26\times 10^{-4}$     &1.93  \\
\hline
\end{tabular}
\end{table}

\begin{figure}
\begin{center}
\begin{tabular}{ccc}
\includegraphics[width=0.333\textwidth]{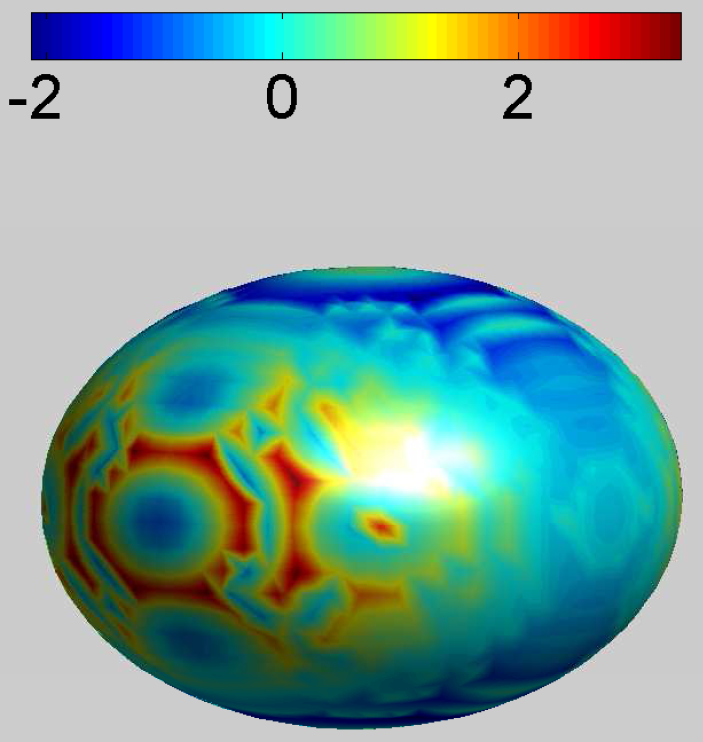}
\includegraphics[width=0.333\textwidth]{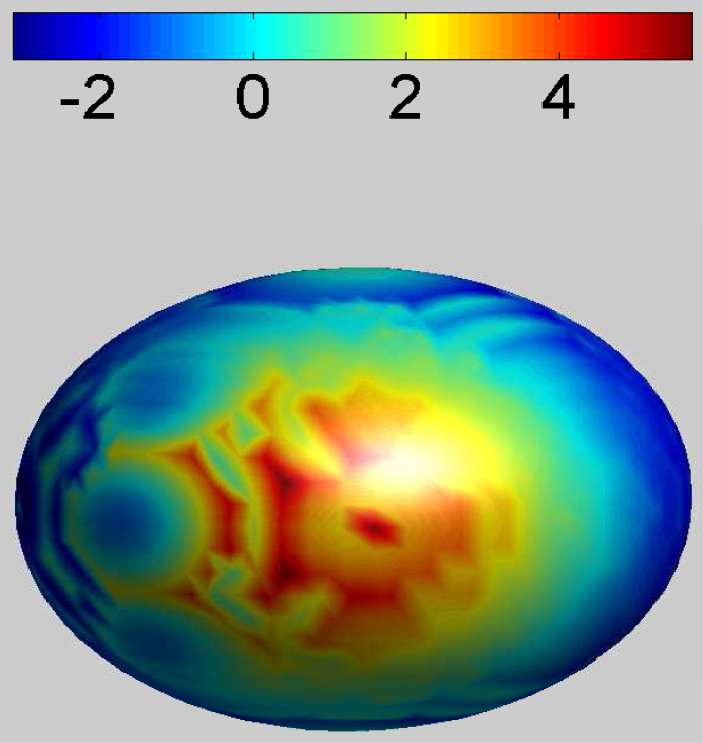}
\includegraphics[width=0.333\textwidth]{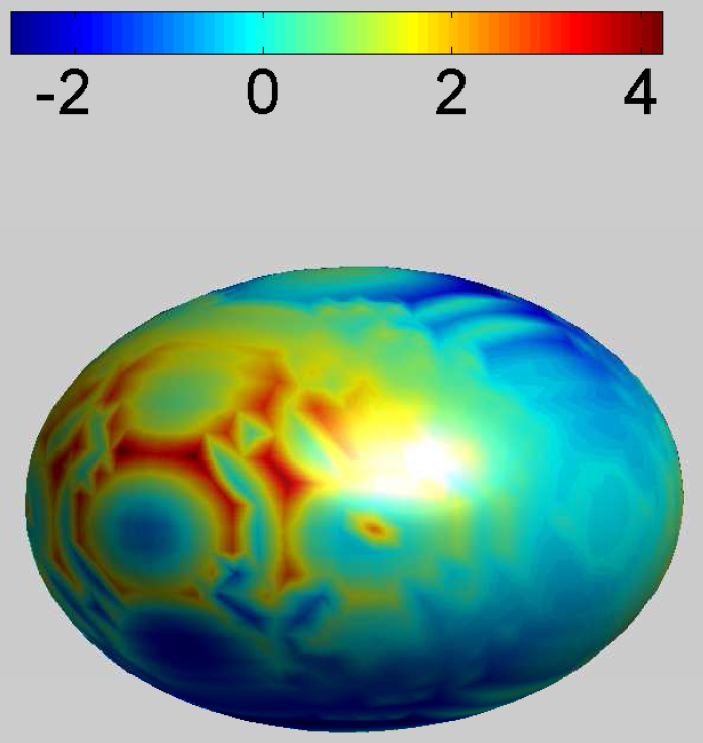}
\end{tabular}
\end{center}
\caption{Numerical solution to the Case 1 of the ellipsoid interface problem  with 40 grid points along each direction  of the computational domain.
Left chart: $u_1$; Middle chart $u_2$; Right chart: $u_3$.
}
\label{ell_sol}
\end{figure}

\begin{figure}
\begin{center}
\begin{tabular}{ccc}
\includegraphics[width=0.333\textwidth]{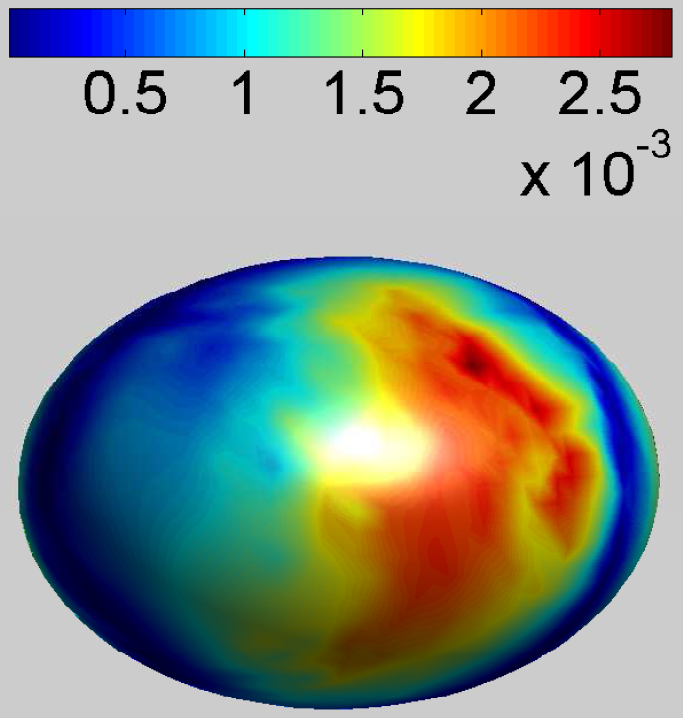}
\includegraphics[width=0.333\textwidth]{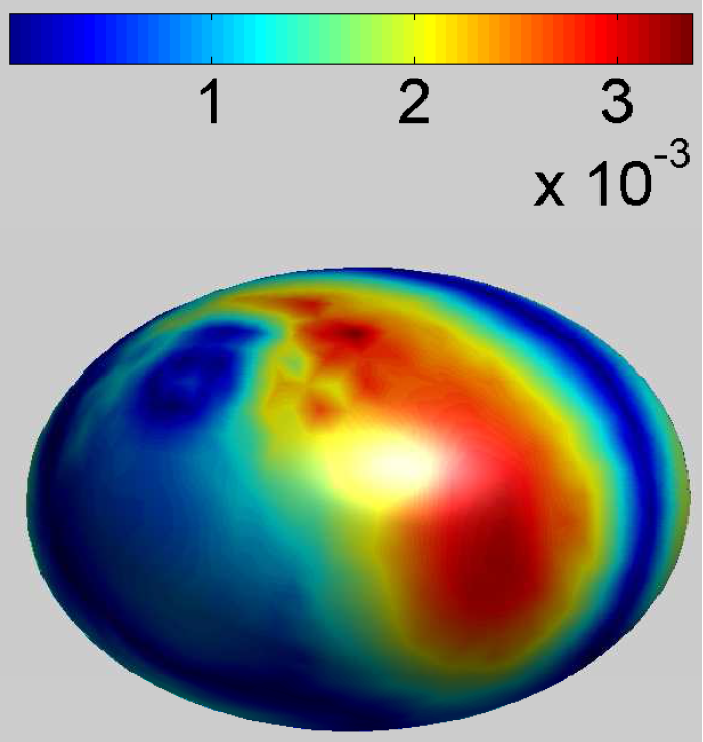}
\includegraphics[width=0.333\textwidth]{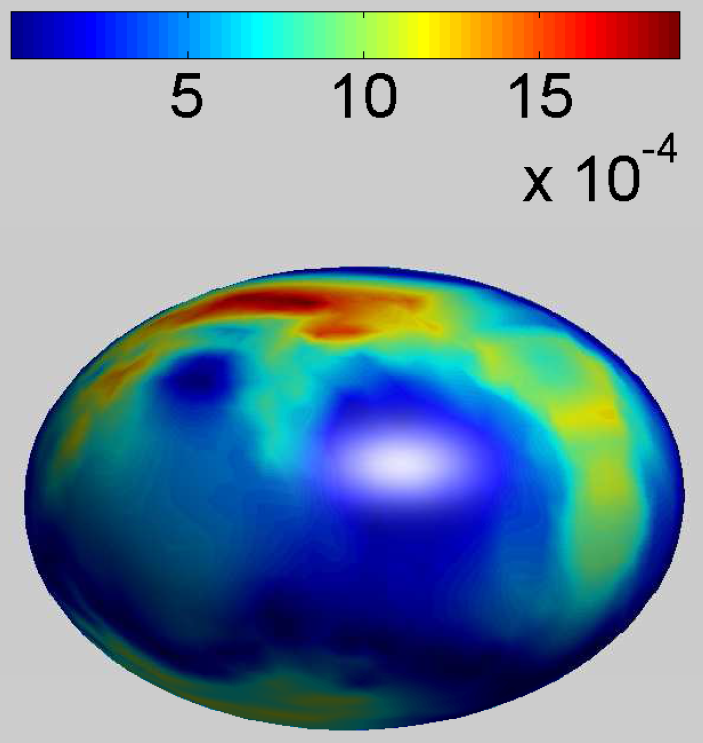}
\end{tabular}
\end{center}
\caption{Numerical error for solving  the Case 1 of the ellipsoid interface problem with 40 grid points along each direction of the computational domain.
Left chart: $u_1$; Middle chart $u_2$; Right chart: $u_3$.
}
\label{ell_err}
\end{figure}

The numerical solution and error of the ellipsoid interface problem are illustrated in Figs. \ref{ell_sol}-\ref{ell_err} with 40 grid points along each direction of the computational domain.

\textbf{Case~2.}
Grid refinement analysis for the $L_\infty$ error is demonstrated in Table \ref{Ex3_case2_linf}. A similar grid refinement  analysis for $L_2$ error is illustrated in Table \ref{Ex3_case2_2}.
\begin{table}
\caption{The $L_\infty$ errors for the Case 2 of Example 3.}
\label{Ex3_case2_linf}
\centering
\begin{tabular}{lllllllll}
\hline
%\multicolumn{5}{c}{The grid refinement analysis of $L_\infty$ error for example case2} \\
\cline{1-7}
$n_x \times n_y \times n_z$       &  $\  L_\infty(u_1) $ & Order              &$ L_\infty(u_2)$ & Order                &$  L_\infty(u_3)$   &Order   \\
\hline
$10\times 10\times10$ &$5.00\times 10^{-2}$     &           &$5.10\times 10^{-2}$      &     &$3.37\times 10^{-2}$     &        \\
$20\times 20\times20$ &$1.26\times 10^{-2}$     &1.99       &$1.37\times 10^{-2}$     &1.90  &$8.01\times 10^{-2}$     & 2.07   \\
$40\times40\times40$ &$3.24\times 10^{-3}$     &1.96       &$3.63\times 10^{-3}$      &1.92  &$2.00\times 10^{-3}$     &2.00  \\
$80\times80\times80$ &$7.73\times 10^{-4}$     &2.07       &$9.98\times 10^{-4}$      &1.87  &$5.15\times 10^{-4}$     &1.96  \\
\hline
\end{tabular}
\end{table}

\begin{table}
\caption{The $L_2$ errors for Case 2 of Example 3.}
\label{Ex3_case2_2}
\centering
\begin{tabular}{lllllllll}
\hline
%\multicolumn{5}{c}{The grid refinement analysis of $L_2$ error for example 3 case 2} \\
\cline{1-7}
$n_x \times n_y \times n_z$       & $L_2(u_1)$ &Order              &$L_2(u_2)$      &Order                &$L_2(u_3)$   &Order \\
\hline
$10\times 10\times10$ &$1.20\times 10^{-2}$     &           &$1.50\times 10^{-2}$      &     &$6.93\times 10^{-3}$     &        \\
$20\times 20\times20$ &$3.77\times 10^{-3}$     &1.67       &$4.13\times 10^{-3}$     &1.70  &$1.91\times 10^{-3}$     &1.86        \\
$40\times40\times40$ &$9.18\times 10^{-4}$     &2.04       &$1.12\times 10^{-4}$      &1.88  &$4.88\times 10^{-4}$     &1.97  \\
$80\times80\times80$ &$2.36\times 10^{-4}$     &1.96       &$2.94\times 10^{-4}$      &1.93  &$1.30\times 10^{-4}$     &1.91  \\
\hline
\end{tabular}
\end{table}

\textbf{Case~3.}
Grid refinement analysis for $L_\infty$ error is demonstrated in Table \ref{Ex3_case3_linf}. We also illustrate
the grid refinement analysis in terms of  $L_2$ error in Table \ref{Ex3_case1_3}.
 The second order convergence of the MIB algorithm is  essentially observed from all the  numerical tests in Example 3.

\begin{table}
\caption{The $L_\infty$ errors for the Case 3 of Example 3.}
\label{Ex3_case3_linf}
\centering
\begin{tabular}{lllllllll}
\hline
%\multicolumn{5}{c}{The grid refinement analysis of $L_\infty$ error for example case 3} \\
\cline{1-7}
$n_x \times n_y \times n_z$       &  $\  L_\infty(u_1) $ & Order              &$ L_\infty(u_2)$ & Order                &$  L_\infty(u_3)$   &Order    \\
\hline
$10\times 10\times10$ &$3.97\times 10^{-2}$     &           &$5.23\times 10^{-2}$      &     &$3.16\times 10^{-2}$     &        \\
$20\times 20\times20$ &$1.22\times 10^{-2}$     &1.70       &$1.31\times 10^{-2}$     &2.00  &$8.07\times 10^{-2}$     &1.97   \\
$40\times40\times40$ &$2.82\times 10^{-3}$     &2.11       &$3.45\times 10^{-3}$      &1.93  &$1.90\times 10^{-3}$     &2.07  \\
$80\times80\times80$ &$6.99\times 10^{-4}$     &2.01       &$8.81\times 10^{-4}$      &1.97  &$4.82\times 10^{-4}$     &1.98  \\
\hline
\end{tabular}
\end{table}

\begin{table}
\caption{The $L_2$ errors for the Case 3 of Example 3.}
\label{Ex3_case1_3}
\centering
\begin{tabular}{lllllllll}
\hline
%\multicolumn{5}{c}{The grid refinement analysis of $L_2$ error for example 3 case 3} \\
\cline{1-7}
$n_x \times n_y \times n_z$       & $L_2(u_1)$ &Order              &$L_2(u_2)$      &Order                &$L_2(u_3)$   &Order \\
\hline
$10\times 10\times10$ &$1.16\times 10^{-2}$     &           &$1.45\times 10^{-2}$      &     &$7.72\times 10^{-2}$     &        \\
$20\times 20\times20$ &$3.24\times 10^{-3}$     &1.84       &$3.96\times 10^{-3}$     &1.99  &$1.86\times 10^{-3}$     &2.05        \\
$40\times40\times40$ &$8.61\times 10^{-4}$     &1.91       &$1.08\times 10^{-3}$      &1.87  &$4.78\times 10^{-4}$     &1.96  \\
$80\times80\times80$ &$2.23\times 10^{-4}$     &1.95       &$2.86\times 10^{-4}$      &1.92  &$1.26\times 10^{-4}$     &1.92  \\
\hline
\end{tabular}
\end{table}

\begin{remark}
In all above examples, the continuity of the solution across the interface, i.e., the no fracture condition, has been carefully maintained in designing the analytical solutions. However, for real world problems, having fractures at the interface is very common. In the following two numerical experiments, the continuity condition of the function values across the interface is dropped.  We test our method for handling general jumps of the function values across the interface.  Numerically, this situation is slightly more difficult to deal with.
\end{remark}

%***************************************************************************************************************
\textbf{Example~4.}
The computational domain is set to be $[-2, 2]\times[-2, 2]\times[-2, 4.4]$ with a cylinder interface defined as
$$
\left\{\begin{array}{ll}
x^2+y^2=\frac{\pi^2}{4}, &\ \ ,\\
z\leq \pi, &\ \  ,\\
z\geq 0.
\end{array}\right.
$$

The Dirichlet boundary condition and interface conditions are determined from the following exact solutions.
$$
u_1(x, y)=
\left\{\begin{array}{ll}
x^2+y^2+z^2-4, &\ \  \mbox{in}\ \ \Omega^+,\\
\cos(x)\cos(y)\cos(z), &\ \  \mbox{in}\ \ \Omega^-.
\end{array}\right.
$$

$$
u_2(x, y)=
\left\{\begin{array}{ll}
x^2+y^2+z^2+xy-4, &\ \  \mbox{in}\ \ \Omega^+,\\
xy+\cos(x)\cos(y)\cos(z), &\ \  \mbox{in}\ \ \Omega^-.
\end{array}\right.
$$
and
$$
u_3(x, y)=
\left\{\begin{array}{ll}
x^2+y^2+z^2+yz-4, &\ \  \mbox{in}\ \ \Omega^+,\\
yz+\cos(x)\cos(y)\cos(z), &\ \  \mbox{in}\ \ \Omega^-.
\end{array}\right.
$$

The values of the Poisson's ratio and shear modulus are, respectively, set to
$$
\nu=
\left\{\begin{array}{ll}
\nu^+=0.20, &\ \  \mbox{in}\ \ \Omega^+,\\
\nu^-=0.24, &\ \  \mbox{in}\ \ \Omega^-.
\end{array}\right.
$$
and
$$
\mu=
\left\{\begin{array}{ll}
\mu^+=1500000, &\ \  \mbox{in}\ \ \Omega^+,\\
\mu^-=2000000, &\ \  \mbox{in}\ \ \Omega^-.
\end{array}\right.
$$

Table \ref{Ex4_linf} offers the grid refinement  analysis of the $L_\infty$ error for Example 4.
Similar grid refinement analysis of the $L_2$ error is given in Table \ref{Ex4_2} for Example 4.
A high level of accuracy and a robust order of convergence are observed from these tests.

\begin{table}
\caption{The $L_\infty$ errors of Example 4.}
\label{Ex4_linf}
\centering
\begin{tabular}{lllllllll}
\hline
%\multicolumn{5}{c}{The grid refinement analysis of $L_\infty$ error for Example 4} \\
\cline{1-7}
$n_x \times n_y \times n_z$       &  $\  L_\infty(u_1) $ & Order              &$ L_\infty(u_2)$ & Order                &$  L_\infty(u_3)$   &Order    \\
\hline
$20\times 20\times20$ &$4.68\times 10^{-3}$     &           &$4.68\times 10^{-3}$      &   &$7.07\times 10^{-3}$     &        \\
$40\times 40\times40$ &$1.16\times 10^{-3}$     &2.01       &$1.17\times 10^{-3}$     &2.00 &$1.74\times 10^{-3}$     & 2.02   \\
$80\times80\times80$ &$2.87\times 10^{-4}$     &2.02       &$2.91\times 10^{-4}$      &2.00  &$4.23\times 10^{-4}$     &2.04  \\
\hline
\end{tabular}
\end{table}

\begin{table}
\caption{The $L_2$ errors of Example 4.}
\label{Ex4_2}
\centering
\begin{tabular}{lllllllll}
\hline
%\multicolumn{5}{c}{The grid refinement analysis of $L_2$ error for Example 4} \\
\cline{1-7}
$n_x \times n_y \times n_z$       & $L_2(u_1)$ &Order              &$L_2(u_2)$      &Order                &$L_2(u_3)$   &Order \\
\hline
$20\times 20\times20$ &$1.04\times 10^{-3}$     &           &$1.04\times 10^{-3}$      &     &$1.58\times 10^{-3}$      &        \\
$40\times 40\times40$ &$2.61\times 10^{-4}$     &1.99       &$2.62\times 10^{-4}$      &1.99  &$3.86\times 10^{-4}$     &2.03        \\
$80\times80\times80$ &$6.69\times 10^{-5}$      &1.96       &$6.77\times 10^{-5}$      &1.95  &$9.77\times 10^{-5}$     &1.98  \\
\hline
\end{tabular}
\end{table}

\begin{figure}
\begin{center}
\begin{tabular}{ccc}
\includegraphics[width=0.333\textwidth]{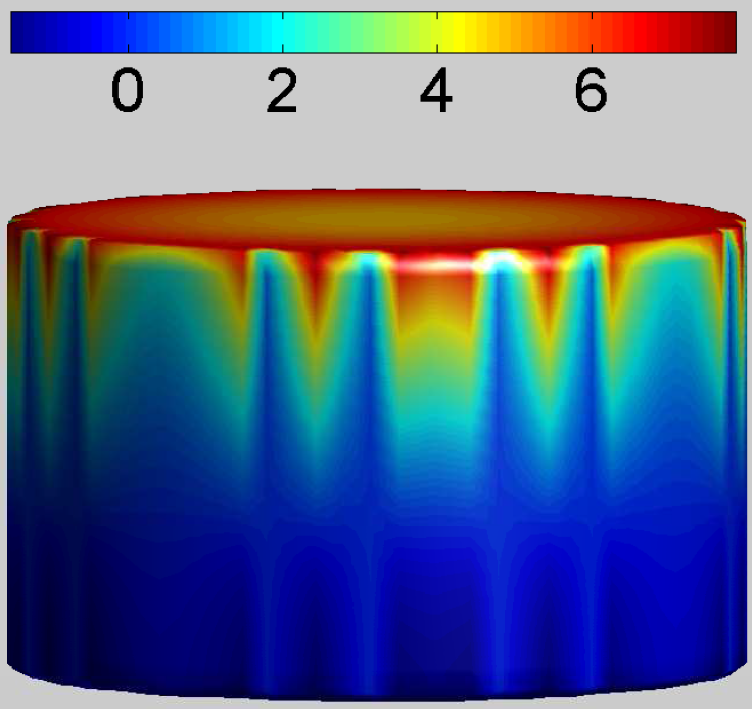}
\includegraphics[width=0.333\textwidth]{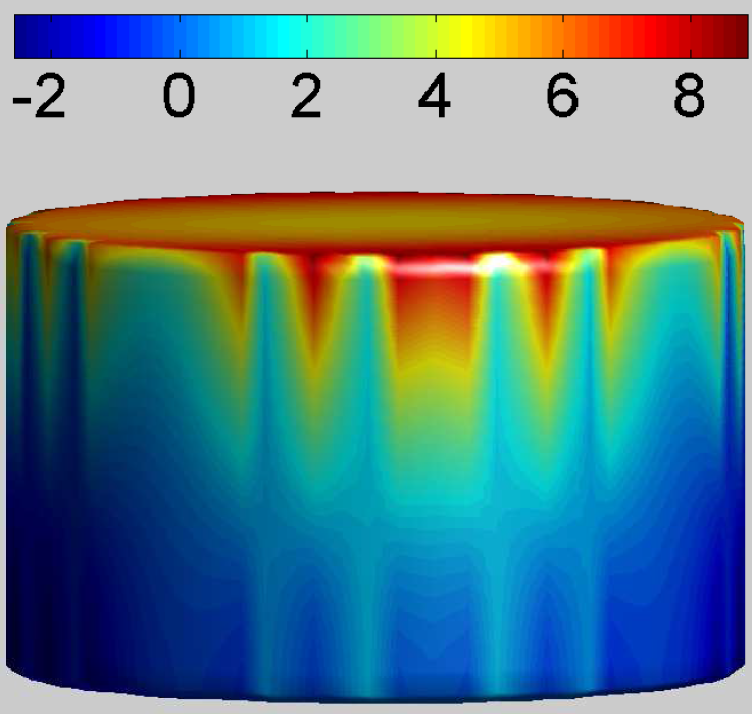}
\includegraphics[width=0.333\textwidth]{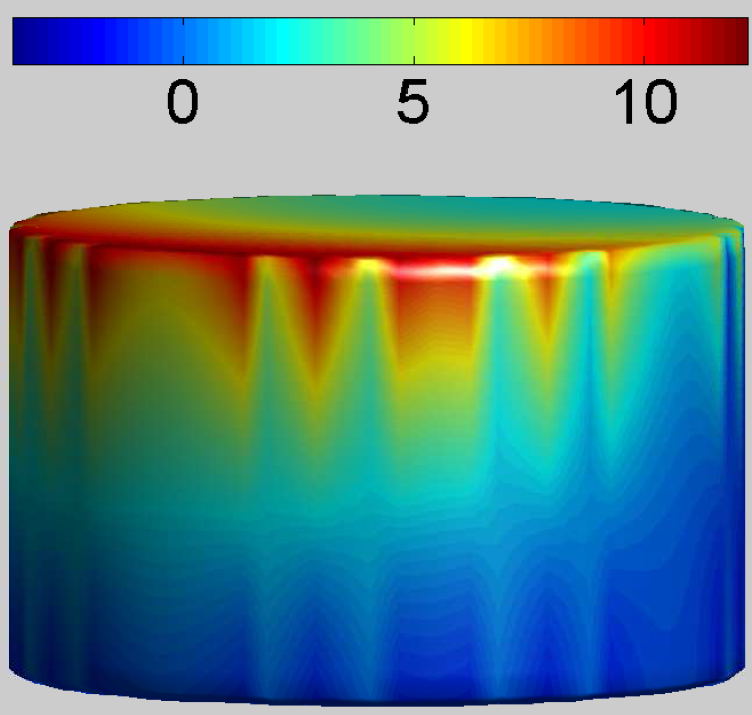}
\end{tabular}
\end{center}
\caption{Numerical solution to the cylinder interface problem with 40 grid points along each direction of the computational domain.
Left chart: $u_1$; Middle chart $u_2$; Right chart: $u_3$.
}
\label{cyl_sol}
\end{figure}

\begin{figure}
\begin{center}
\begin{tabular}{ccc}
\includegraphics[width=0.333\textwidth]{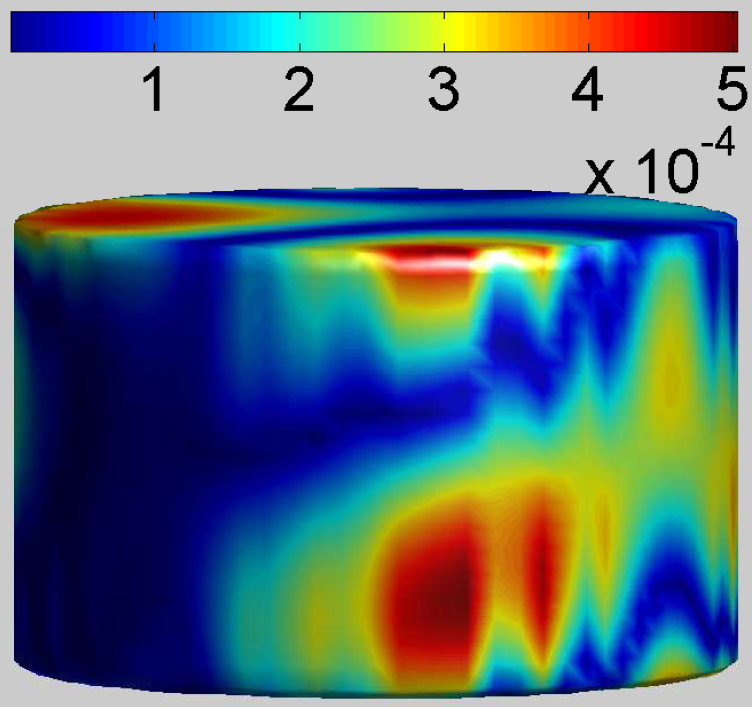}
\includegraphics[width=0.333\textwidth]{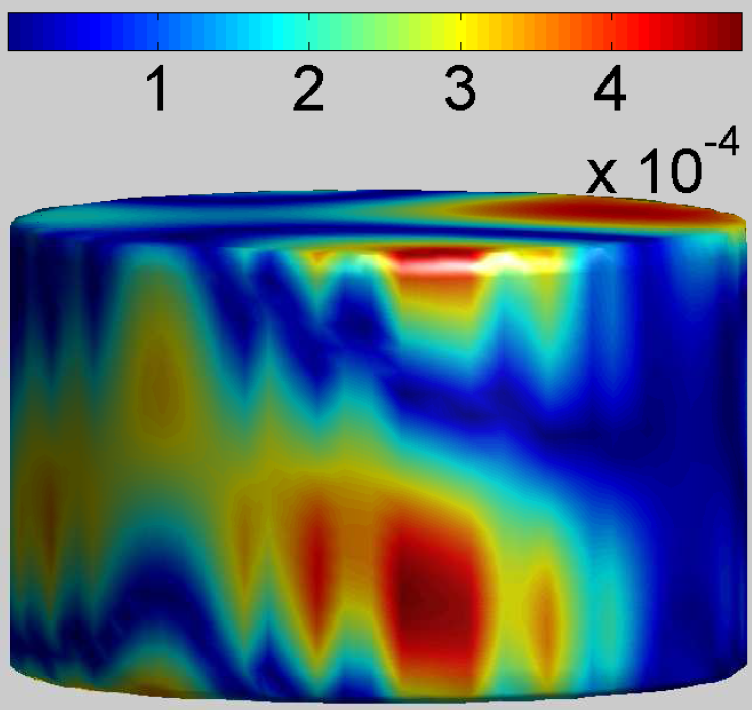}
\includegraphics[width=0.333\textwidth]{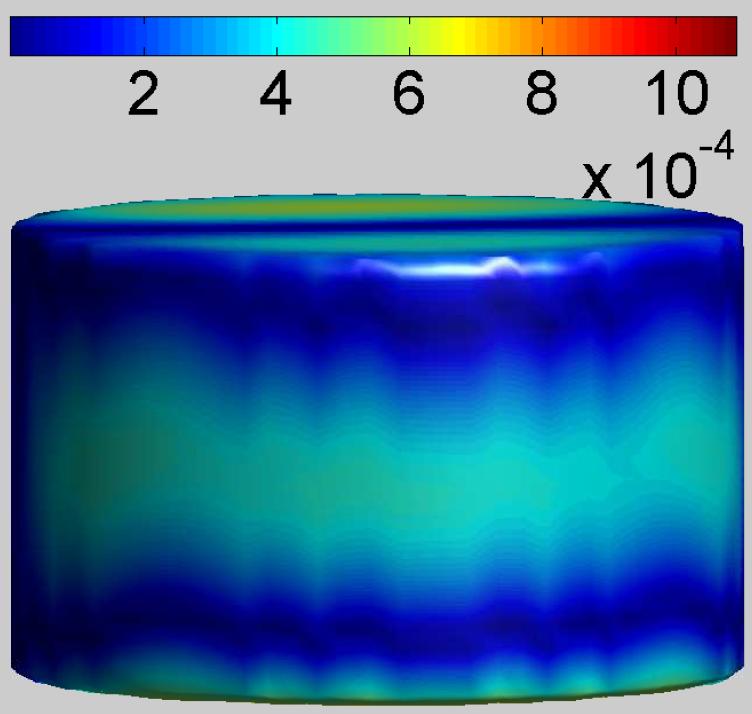}
\end{tabular}
\end{center}
\caption{Numerical error of solving  the cylinder interface problem with 40 grid points along each direction of the computational domain.
Left chart: $u_1$; Middle chart $u_2$; Right chart: $u_3$.
}
\label{cyl_err}
\end{figure}

Numerical solution and error are depicted in Figs. \ref{cyl_sol} and \ref{cyl_err}, respectively, where 40 grid points are used along each direction  of
the computational domain.  Obviously, the error is very small in all of three solutions.

%***************************************************************************************************************
\textbf{Example~5.}
Geometric complexity is a major issue in interface problems. It is often the case that numerical methods designed for simple interface geometries do not work  for  complex interface geometries. In this example, we  consider a  more complicated interface, which is defined to be a torus
$$
\left\{\begin{array}{ll}
x(u, v)=(R+r\cos v)\cos u, &\ \ \\
y(u, v)=(R+r\sin v)\sin u, &\ \  \\
z(u, v)=r\sin v,
\end{array}\right.
$$
where $u, v\in [0,2\pi]$ are two parameters. The computational domain is set to be $[-10, 10]\times[-10, 10]\times[-5, 5]$

The above torus can also be represented as
$$
(R-\sqrt{x^2+y^2})^2+z^2=r^2.
$$
We set  $R=4, r=2$ in our numerical experiments.

The Poisson's ratio, shear modulus and  designed analytic solution in Example 4 are adopted for this example.
Grid refinement analysis in terms of $L_\infty$ error is given in  Table \ref{Ex5_linf}. A similar grid refinement analysis in terms of $L_2$ error is given Table \ref{Ex5_2}.
Although there is some small fluctuation in the convergent order, the second order convergence is essentially obtained in this test.

\begin{table}
\caption{The $L_\infty$ errors of Example 5.}
\label{Ex5_linf}
\centering
\begin{tabular}{lllllllll}
\hline
%\multicolumn{5}{c}{The grid refinement analysis of $L_\infty$ error for Example 5} \\
\cline{1-7}
$n_x \times n_y \times n_z$       &  $\  L_\infty(u_1) $ & Order              &$ L_\infty(u_2)$ & Order                &$  L_\infty(u_3)$   &Order   \\
\hline
$20\times 20\times20$ &$2.04\times 10^{-1}$     &           &$2.04\times 10^{-1}$      &   &$1.12\times 10^{-1}$     &        \\
$40\times 40\times40$ &$4.14\times 10^{-2}$     &2.30       &$4.05\times 10^{-2}$     &2.33 &$2.34\times 10^{-2}$      & 2.09   \\
$80\times80\times80$ &$1.24\times 10^{-2}$      &1.74       &$1.09\times 10^{-2}$      &1.89  &$4.66\times 10^{-3}$     &1.97  \\
\hline
\end{tabular}
%\caption{The grid refinement analysis of $L_\infty$ error for Example 5}
\end{table}

\begin{table}
\caption{The $L_2$ errors of Example 5.}
\label{Ex5_2}
\centering
\begin{tabular}{lllllllll}
\hline
%\multicolumn{5}{c}{The grid refinement analysis of $L_2$ error for Example 5} \\
\cline{1-7}
$n_x \times n_y \times n_z$       & $L_2(u_1)$ &Order              &$L_2(u_2)$      &Order                &$L_2(u_3)$   &Order \\
\hline
$20\times 20\times20$ &$4.54\times 10^{-2}$     &           &$4.52\times 10^{-2}$      &     &$1.87\times 10^{-2}$     &        \\
$40\times 40\times40$ &$1.12\times 10^{-2}$     &1.71        &$1.10\times 10^{-2}$     &2.04  &$4.40\times 10^{-3}$     &2.09        \\
$80\times80\times80$ &$2.92\times 10^{-3}$      &2.04       &$2.90\times 10^{-3}$      &1.92  &$1.12\times 10^{-3}$     &1.97  \\
\hline
\end{tabular}
%\caption{The grid refinement analysis of $L_2$ error for Example 5}
\end{table}

\begin{figure}
\begin{center}
\begin{tabular}{ccc}
\includegraphics[width=0.333\textwidth]{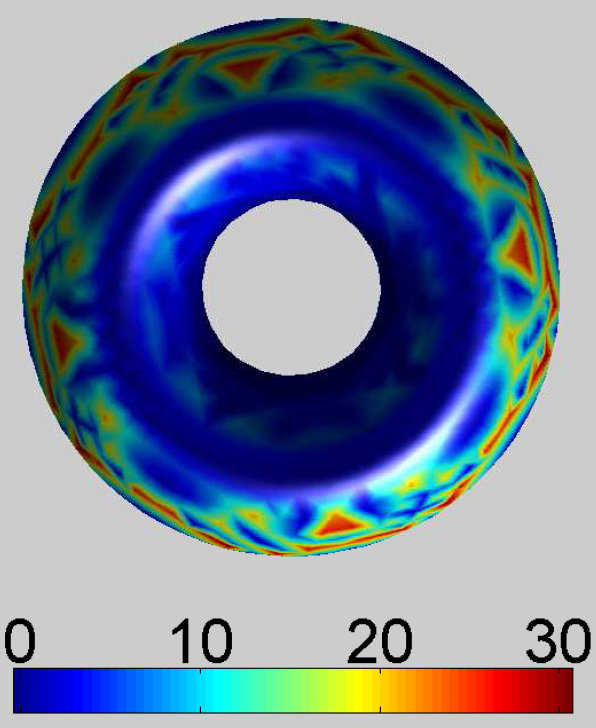}
\includegraphics[width=0.333\textwidth]{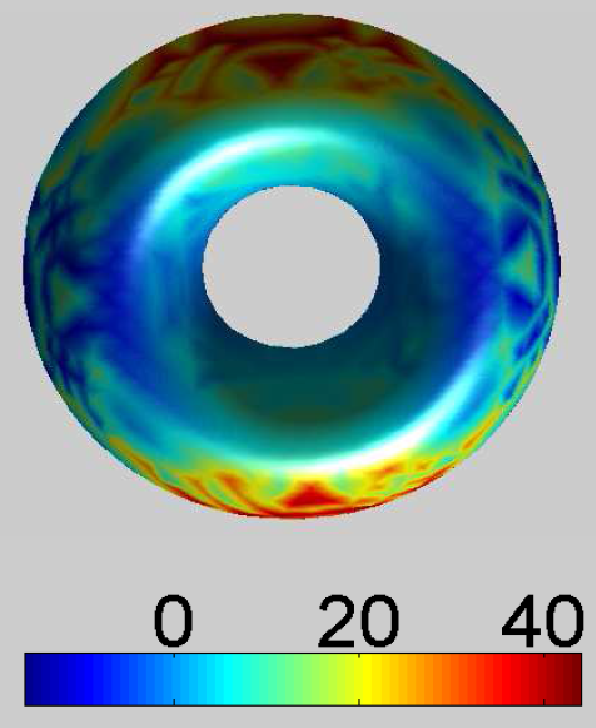}
\includegraphics[width=0.333\textwidth]{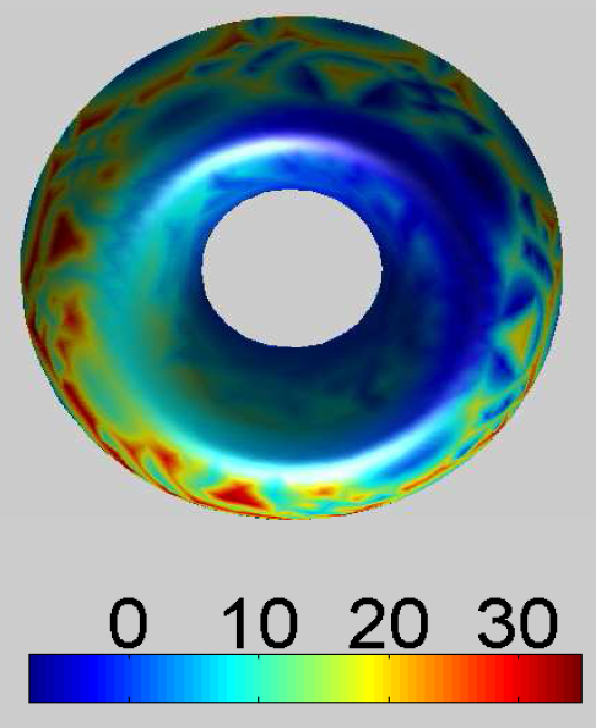}
\end{tabular}
\end{center}
\caption{Numerical solution to the torus interface problem with 40 grid points along each direction  of the computational domain. Left chart: $u_1$; Middle chart: $u_2$;  Right chart: $u_3$.
}
\label{torus_sol}
\end{figure}

\begin{figure}
\begin{center}
\begin{tabular}{ccc}
\includegraphics[width=0.333\textwidth]{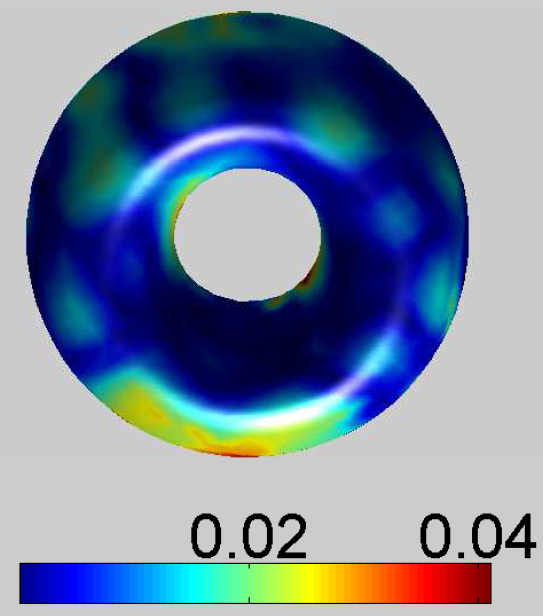}
\includegraphics[width=0.333\textwidth]{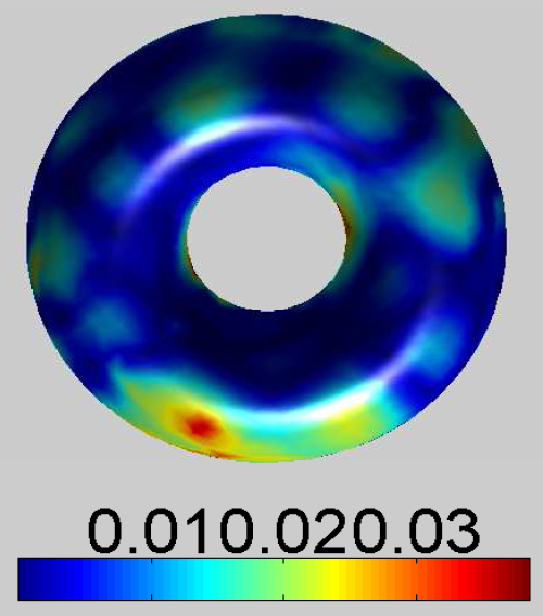}
\includegraphics[width=0.333\textwidth]{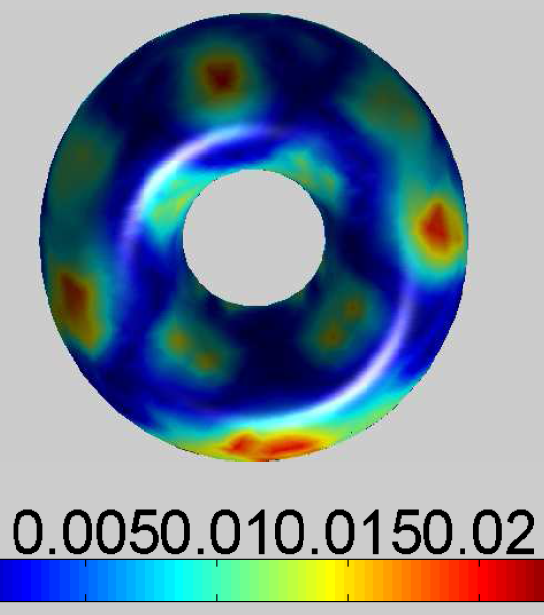}
\end{tabular}
\end{center}
\caption{Numerical error to the torus interface problem with 40 grid points along each direction of the computational domain, Left chart: $u_1$; Middle chart: $u_2$;  Right chart: $u_3$.
}
\label{torus_err}
\end{figure}

%\begin{figure}[!ht]
%\small
%\centering
%\includegraphics[width=12.5cm,height=5.5cm]{sol/torus_sol_40.png}
%\caption{Solution to the torus interface with 40 grids along each side of the computational domain, left $u_1$, middle $u_2$, right $u_3$}
%\label{torus_sol}
%\end{figure}
%
%
%\begin{figure}[!ht]
%\small
%\centering
%\includegraphics[width=12.5cm,height=5.5cm]{sol/torus_err_40.png}
%\caption{Error to the torus interface with 40 grids along each side of the computational domain, left $u_1$, middle $u_2$, right $u_3$}
%\label{torus_err}
%\end{figure}

Figures \ref{torus_sol} and \ref{torus_err} illustrate  the numerical solution and the error in a $40\times40\times40$ mesh. Note that errors appear large in this test example. However, the amplitude of the solution is much larger too, due to a much larger computational domain.

\textbf{Example~6.}
For the last example of the smooth interface with piecewise constant material parameters, we consider a   more complicated interface geometry, i.e., a flower-like cylinder interface. The interface can be represented as
$$
\left\{\begin{array}{ll}
r=\frac{5}{2}+\frac{5}{7}\sin{5\theta},\\
-\frac{2}{3}\leq z \leq \frac{2}{3},
\end{array}\right.
$$
where $r=\sqrt{x^2+y^2}$ and $\theta=\arctan{\frac{y}{x}}$. The computational domain is set to $[-5, 5]\times[-5, 5]\times[-2, 2]$.

Material parameters and exact solutions designed in Example 4 are utilized in this example.
Grid refinement analysis in terms of $L_\infty$ error is given in Table \ref{Ex66_linf} and a similar analysis  in terms of $L_2$ error is given in Table \ref{Ex66_2}. It is quite interesting to see that good convergent orders are attained.

\begin{table}
\caption{The $L_\infty$ errors of Example 6.}
\label{Ex66_linf}
\centering
\begin{tabular}{lllllllll}
\hline
%\multicolumn{5}{c}{The grid refinement analysis of $L_\infty$ error for Example 6} \\
\cline{1-7}
Grid Size       &$L_\infty(u_1)$ &Order              &$L_\infty(u_2)   $ &Order                &$L_\infty(u_3)  $   &Order  \\
\hline
$0.5$   &$4.29\times 10^{-2}$     &           &$4.49\times 10^{-2}$      &      &$1.95\times 10^{-2}$     &        \\
$0.25$  &$9.04\times 10^{-3}$     &2.25       &$9.46\times 10^{-3}$      &2.25   &$4.97\times 10^{-3}$    & 1.97   \\
$0.125$ &$1.96\times 10^{-3}$     &2.21       &$2.17\times 10^{-3}$      &2.03  &$7.02\times 10^{-3}$     & 2.82  \\
\hline
\end{tabular}
%\caption{The grid refinement analysis of $L_\infty$ error for Example 6}
\end{table}

\begin{table}
\caption{The $L_2$ errors of Example 6.}
\label{Ex66_2}
\centering
\begin{tabular}{lllllllll}
\hline
%\multicolumn{5}{c}{The grid refinement analysis of $L_2$ error for Example 6} \\
\cline{1-7}
Grid Size       & $L_2(u_1)  $ &Order              &Error  $L_2(u_2)    $ &Order                &    $L_2(u_3)  $   &Order  \\
\hline
$0.5$   &$4.12\times 10^{-3}$     &           &$4.96\times 10^{-3}$      &      &$2.35\times 10^{-3}$     &        \\
$0.25$  &$9.90\times 10^{-4}$     &2.06       &$1.11\times 10^{-3}$     &2.16   &$4.10\times 10^{-4}$     &2.52   \\
$0.125$ &$2.11\times 10^{-4}$     &2.23       &$2.38\times 10^{-4}$      &2.22  &$7.68\times 10^{-5}$     &2.42  \\
\hline
\end{tabular}
%\caption{The grid refinement analysis of $L_2$ error for Example 6}
\end{table}

\begin{figure}
\begin{center}
\begin{tabular}{ccc}
\includegraphics[width=0.333\textwidth]{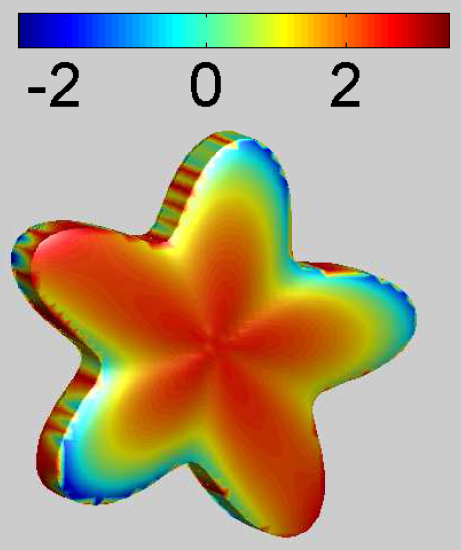}
\includegraphics[width=0.333\textwidth]{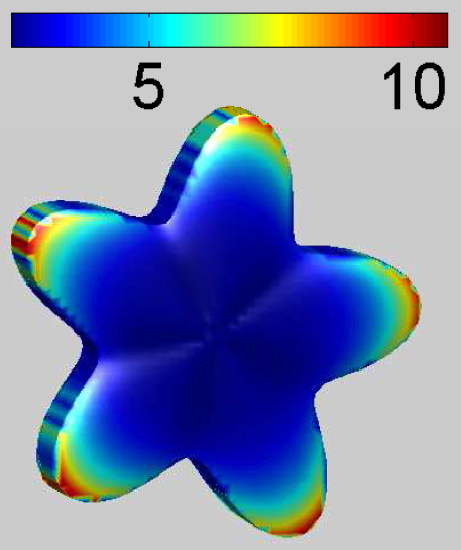}
\includegraphics[width=0.333\textwidth]{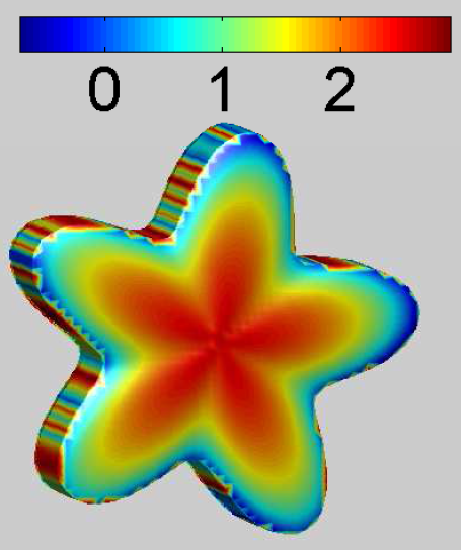}
\end{tabular}
\end{center}
\caption{Numerical solution to the flower interface problem with grid size 0.125. Left chart: $u_1$; Middle chart: $u_2$;  Right chart: $u_3$.}
\label{flower_sol}
\end{figure}

\begin{figure}
\begin{center}
\begin{tabular}{ccc}
\includegraphics[width=0.333\textwidth]{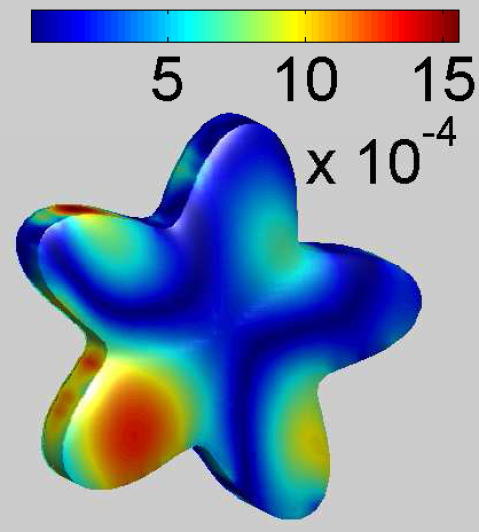}
\includegraphics[width=0.333\textwidth]{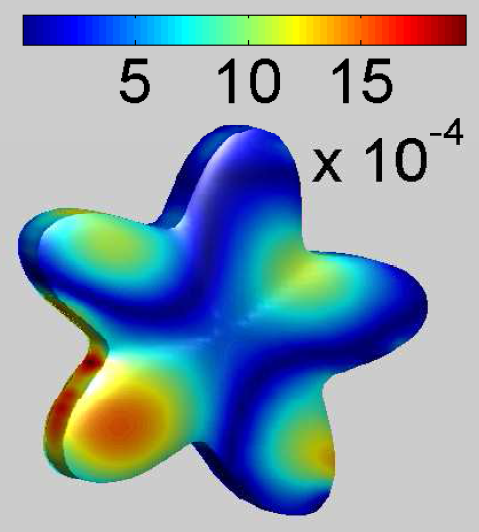}
\includegraphics[width=0.333\textwidth]{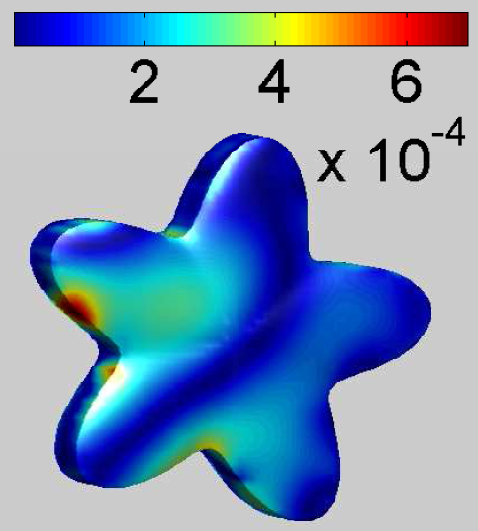}
\end{tabular}
\end{center}
\caption{Numerical error for solving  the flower interface problem with grid size 0.125. Left chart: $u_1$; Middle chart: $u_2$;  Right chart: $u_3$.}
\label{flower_err}
\end{figure}

%\begin{figure}[!ht]
%\small
%\centering
%\includegraphics[width=12.5cm,height=5.5cm]{sol/flower_sol_125.png}
%\caption{Solution to the flower interface with grid size 0.125, left $u_1$, middle $u_2$, right $u_3$}
%\label{flower_sol}
%\end{figure}
%
%
%\begin{figure}[!ht]
%\small
%\centering
%\includegraphics[width=12.5cm,height=5.5cm]{sol/flower_err_125.png}
%\caption{Error to the flower interface with grid size 0.125, left $u_1$, middle $u_2$, right $u_3$}
%\label{flower_err}
%\end{figure}

Figures \ref{flower_sol} and \ref{flower_err} demonstrate  the solution and error of the flower-liked interface problem with grid size $0.125$.
Note that the error amplitude depends on the mesh size.

\subsubsection{Position dependent shear modulus}

Spatially varying shear modulus occurs frequently in natural and man-made materials and devices. The ability to deal with position-dependent material parameters cannot be overemphasized for practical applications.  For example, the  protein molecules can have variable shear modulus \cite{KLXia:2013d}.
 In this subsection, we consider that the shear modulus is given as a  position-dependent function.

%*****************************************************************************************************************
\textbf{Example~7.}
In this example, we consider the problem defined in Example 1, while  replace the shear modulus in Example 1 by the following position-dependent function
$$
\mu=
\left\{\begin{array}{ll}
\mu^+=1500000+(x+y+z), &\ \  \mbox{in}\ \ \Omega^+,\\
\mu^-=2000000+xyz, &\ \  \mbox{in}\ \ \Omega^-.
\end{array}\right.
$$

The error analysis is given in Tables \ref{Ex6_linf} and \ref{Ex6_2} for $L_\infty$ and $L_2$, respectively. Essentially, second order convergence is obtained. The level of accuracy is the same as that obtained for Example 1, which indicates that the proposed MIB method is not sensitive to position-dependent material parameters.

\begin{table}
\caption{The $L_\infty$ errors of Example 7.}
\label{Ex6_linf}
\centering
\begin{tabular}{lllllllll}
\hline
%\multicolumn{5}{c}{The grid refinement analysis of $L_\infty$ error for Example 7} \\
\cline{1-7}
$n_x\times n_y\times n_z$      &$L_\infty(u_1)$ &Order              &$L_\infty(u_2)   $ &Order                &$L_\infty(u_3)  $   &Order  \\
\hline
$10\times 10\times10$ &$6.61\times 10^{-2}$     &           &$6.27\times 10^{-2}$      &    &$5.67\times 10^{-2}$     &        \\
$20\times 20\times20$ &$1.37\times 10^{-2}$     &2.27       &$1.34\times 10^{-2}$     &2.27 &$1.31\times 10^{-2}$     & 2.11   \\
$40\times 40\times40$ &$2.66\times 10^{-3}$     &2.36       &$2.84\times 10^{-3}$     &2.24 &$2.67\times 10^{-3}$     & 2.29   \\
$80\times80\times80$ &$7.41\times 10^{-4}$     &1.90       &$7.15\times 10^{-4}$      &1.99 &$7.26\times 10^{-4}$     &1.89  \\
\hline
\end{tabular}
%\caption{The grid refinement analysis of $L_\infty$ error for Example 7}
\end{table}

\begin{table}
\caption{The $L_2$ errors of Example 7.}
\label{Ex6_2}
\centering
\begin{tabular}{lllllllll}
\hline
%\multicolumn{5}{c}{The grid refinement analysis of $L_2$ error for Example 7} \\
\cline{1-7}
$n_x\times n_y\times n_z$     & $L_2(u_1)  $ &Order              &Error  $L_2(u_2)    $ &Order                &    $L_2(u_3)  $   &Order  \\
\hline
$10\times 10\times10$ &$1.28\times 10^{-2}$     &           &$1.28\times 10^{-2}$      &     &$1.29\times 10^{-2}$     &        \\
$20\times 20\times20$ &$3.18\times 10^{-3}$     &2.01       &$3.19\times 10^{-3}$     &2.00  &$3.14\times 10^{-3}$     & 2.04   \\
$40\times 40\times40$ &$8.30\times 10^{-4}$     &1.94        &$8.36\times 10^{-4}$     &1.93  &$8.26\times 10^{-4}$     &1.93        \\
$80\times80\times80$ &$2.24\times 10^{-4}$     &1.89       &$2.24\times 10^{-4}$      &1.90  &$2.23\times 10^{-4}$     &1.89  \\
\hline
\end{tabular}
%\caption{The grid refinement analysis of $L_2$ error for Example 7}
\end{table}

%******************************************************************************************************************
\textbf{Example~8.}
To further test the proposed method for its performance in dealing with variable material parameters, we consider
an example by setting the shear modulus in Example 4 to the following spatially dependent functions
$$
\mu=
\left\{\begin{array}{ll}
\mu^+=1500000+2000(x+y+z), &\ \  \mbox{in}\ \ \Omega^+,\\
\mu^-=2000000+1500xyz, &\ \  \mbox{in}\ \ \Omega^-.
\end{array}\right.
$$

\begin{table}
\caption{The $L_\infty$ errors of Example 8.}
\label{Ex7_linf}
\centering
\begin{tabular}{lllllllll}
\hline
%\multicolumn{5}{c}{The grid refinement analysis of $L_\infty$ error for Example 8} \\
\cline{1-7}
$n_x\times n_y\times n_z$       &$L_\infty(u_1)$ &Order              &$L_\infty(u_2)   $ &Order                &$L_\infty(u_3)  $   &Order\\
\hline
$10\times 10\times10$ &$1.85\times 10^{-2}$     &           &$1.85\times 10^{-2}$      &   &$3.14\times 10^{-2}$     &        \\
$20\times 20\times20$ &$4.68\times 10^{-3}$     &1.98       &$4.68\times 10^{-3}$     &1.98 &$7.07\times 10^{-3}$     & 2.15   \\
$40\times 40\times40$ &$1.15\times 10^{-3}$     &2.02       &$1.17\times 10^{-3}$     &2.00 &$1.74\times 10^{-3}$     & 2.02   \\
$80\times80\times80$ &$2.99\times 10^{-4}$     &1.94       &$3.19\times 10^{-4}$      &1.87  &$4.23\times 10^{-4}$     &2.01  \\
\hline
\end{tabular}
%\caption{The grid refinement analysis of $L_\infty$ error for Example 8}
\end{table}

\begin{table}
\caption{The $L_2$ errors of Example 8.}
\label{Ex7_2}
\centering
\begin{tabular}{lllllllll}
\hline
%\multicolumn{5}{c}{The grid refinement analysis of $L_2$ error for Example 8} \\
\cline{1-7}
$n_x\times n_y\times n_z$      & $L_2(u_1)  $ &Order              &Error  $L_2(u_2)    $ &Order                &    $L_2(u_3)  $   &Order   \\
\hline
$10\times 10\times10$ &$4.16\times 10^{-3}$     &           &$4.16\times 10^{-3}$      &     &$7.47\times 10^{-3}$     &        \\
$20\times 20\times20$ &$1.04\times 10^{-3}$     &2.00        &$1.04\times 10^{-3}$     &2.00  &$1.58\times 10^{-3}$     &2.24        \\
$40\times 40\times40$ &$2.63\times 10^{-4}$     &1.98        &$2.64\times 10^{-4}$     &1.98  &$3.92\times 10^{-4}$     &2.01        \\
$80\times80\times80$ &$6.82\times 10^{-5}$     &1.95       &$7.07\times 10^{-5}$      &1.90  &$1.00\times 10^{-4}$     &1.97  \\
\hline
\end{tabular}
%\caption{The grid refinement analysis of $L_2$ error for Example 8}
\end{table}

The $L_\infty$ and $L_2$ errors are analyzed in Tables \ref{Ex7_linf} and \ref{Ex7_2}, respectively. %Again, we observe the second order accuracy in both  $L_\infty$ and $L_2$  norms.

%*******************************************************************************************************************
\textbf{Example~9.}

%What was the problem???
In this numerical experiment, we further investigate the robustness of the proposed MIB algorithm to the position dependent shear modulus, now we redo example 8, however, change the shear modulus in example 8 to be the following functions.
$$
\mu=
\left\{\begin{array}{ll}
\mu^+=1500000+2000(x^2+y^2+z^2), &\ \  \mbox{in}\ \ \Omega^+,\\
\mu^-=2000000+1500x^2y^2z^2, &\ \  \mbox{in}\ \ \Omega^-.
\end{array}\right.
$$

The $L_\infty$ and $L_2$ errors are analyzed in Tables \ref{Ex8_linf} and \ref{Ex8_2}, respectively.

\begin{table}
\caption{The $L_\infty$ errors of Example 9.}
\label{Ex8_linf}
\centering
\begin{tabular}{lllllllll}
\hline
%\multicolumn{5}{c}{The grid refinement analysis of $L_\infty$ error for Example 9} \\
\cline{1-7}
Grid Size       &$L_\infty(u_1)$ &Order              &$L_\infty(u_2)   $ &Order                &$L_\infty(u_3)  $   &Order \\
\hline
$20\times 20\times20$ &$1.67\times 10^{-1}$     &           &$1.65\times 10^{-1}$      &   &$9.39\times 10^{-2}$     &        \\
$40\times 40\times40$ &$4.20\times 10^{-2}$     &1.99       &$5.36\times 10^{-2}$     &1.62 &$2.66\times 10^{-2}$     & 1.82   \\
$80\times80\times80$ &$9.97\times 10^{-3}$     &2.07       &$9.71\times 10^{-3}$      &2.48  &$5.96\times 10^{-3}$     &2.43  \\
\hline
\end{tabular}
%\caption{The grid refinement analysis of $L_\infty$ error for Example 9}
\end{table}

\begin{table}
\caption{The $L_2$ errors of Example 9.}
\label{Ex8_2}
\centering
\begin{tabular}{lllllllll}
\hline
%\multicolumn{5}{c}{The grid refinement analysis of $L_2$ error for Example 9} \\
\cline{1-7}
Grid Size       & $L_2(u_1)  $ &Order              &Error  $L_2(u_2)    $ &Order                &    $L_2(u_3)  $   &Order \\
\hline
$20\times 20\times20$ &$4.27\times 10^{-2}$     &           &$4.26\times 10^{-2}$      &     &$1.78\times 10^{-2}$     &        \\
$40\times 40\times40$ &$1.10\times 10^{-2}$     &1.96        &$1.09\times 10^{-2}$     &1.97  &$4.48\times 10^{-3}$     &1.99        \\
$80\times80\times80$ &$2.80\times 10^{-3}$     &1.97       &$2.78\times 10^{-3}$      &1.97  &$1.08\times 10^{-3}$     &2.05  \\
\hline
\end{tabular}
%\caption{The grid refinement analysis of $L_2$ error for Example 9}
\end{table}

\begin{remark}
From the above three examples, we observe the second order accuracy in both $L_\infty$ and $L_2$  norms for elasticity interfaces with  position-dependent material parameters. Additionally the level of accuracy is not affected by the spatially varying material parameters.
\end{remark}

\subsection{Nonsmooth  interfaces}

Nonsmooth interfaces are omnipresent in practical applications and give rise to challenges for numerical algorithm design. In this section, we consider  a few elasticity  interface problems with geometric singularities.

\textbf{Example~10.}
In this example, let us consider an apple-like interface \cite{Yu:2007a}
$$
\rho=1.9\left(1-\cos{\phi}\right),
$$
where $\rho=\sqrt{x^2+y^2+z^2}$ and $\phi=\arccos{\frac{z}{\rho}}$. The computational domain is set to $[-5, 4.6]\times[-5, 4.6]\times[-8, 4]$.

The values of the Poisson's ratio and shear modulus are, respectively
$$
\nu=
\left\{\begin{array}{ll}
\nu^+=0.24, &\ \  \mbox{in}\ \ \Omega^+,\\
\nu^-=0.20, &\ \  \mbox{in}\ \ \Omega^-.
\end{array}\right.
$$
and
$$
\mu=
\left\{\begin{array}{ll}
\mu^+=2000000, &\ \  \mbox{in}\ \ \Omega^+,\\
\mu^-=1500000, &\ \  \mbox{in}\ \ \Omega^-.
\end{array}\right.
$$

The Dirichlet boundary condition and  interface jump conditions can be determined by the following exact solution
$$
u_1(x, y)=
\left\{\begin{array}{ll}
\cos{x}\cos{y}\cos{z}+xyz, &\ \  \mbox{in}\ \ \Omega^+,\\
3, &\ \  \mbox{in}\ \ \Omega^-.
\end{array}\right.
$$

$$
u_2(x, y)=
\left\{\begin{array}{ll}
\cos{x}\cos{y}\cos{z}+x^2+y^2+z^2, &\ \  \mbox{in}\ \ \Omega^+,\\
3, &\ \  \mbox{in}\ \ \Omega^-.
\end{array}\right.
$$
and
$$
u_3(x, y)=
\left\{\begin{array}{ll}
\cos{x}\cos{y}\cos{z}, &\ \  \mbox{in}\ \ \Omega^+,\\
3, &\ \  \mbox{in}\ \ \Omega^-.
\end{array}\right.
$$

Grid refinement analysis in terms of $L_\infty$ error is given in Table \ref{Ex10_linf}. A similar  $L_2$ error analysis is given in Table \ref{Ex10_2}.
The level of accuracy and the order of convergence are similar to those observed in earlier cases.

\begin{table}
\caption{The $L_\infty$ errors of Example 10.}
\label{Ex10_linf}
\centering
\begin{tabular}{lllllllll}
\hline
%\multicolumn{5}{c}{The grid refinement analysis of $L_\infty$ error for Example 10} \\
\cline{1-7}
Grid Size       &$L_\infty(u1)$ &Order              &$L_\infty(u2)   $ &Order                &$L_\infty(u3)  $   &Order \\
\hline
$0.6$   &$5.08\times 10^{-2}$     &           &$5.18\times 10^{-2}$      &      &$6.60\times 10^{-2}$     &        \\
$0.3$  &$1.39\times 10^{-2}$     &1.87       &$1.41\times 10^{-2}$     &2.06   &$1.74\times 10^{-2}$     & 1.92   \\
$0.15$ &$3.07\times 10^{-3}$     &2.18       &$3.77\times 10^{-3}$      &2.18  &$4.09\times 10^{-3}$     &2.09  \\
\hline
\end{tabular}
%\caption{The grid refinement analysis of $L_\infty$ error for Example 10}
\end{table}

\begin{table}
\caption{The $L_2$ errors of Example 10.}
\label{Ex10_2}
\centering
\begin{tabular}{lllllllll}
\hline
%\multicolumn{5}{c}{The grid refinement analysis of $L_2$ error for Example 10} \\
\cline{1-7}
Grid Size       & $L_2(u1)  $ &Order              &Error  $L_2(u2)    $ &Order                &    $L_2(u3)  $   &Order \\
\hline
$0.6$   &$5.86\times 10^{-3}$     &           &$6.11\times 10^{-3}$      &      &$9.44\times 10^{-3}$     &        \\
$0.3$  &$1.51\times 10^{-3}$     &2.17       &$1.61\times 10^{-3}$     &1.92   &$2.55\times 10^{-3}$     &1.89   \\
$0.15$ &$3.77\times 10^{-4}$     &2.07       &$3.96\times 10^{-4}$      &2.02  &$6.69\times 10^{-4}$     &1.93  \\
\hline
\end{tabular}
%\caption{The grid refinement analysis of $L_2$ error for Example 10}
\end{table}

\begin{figure}
\begin{center}
\begin{tabular}{ccc}
\includegraphics[width=0.333\textwidth]{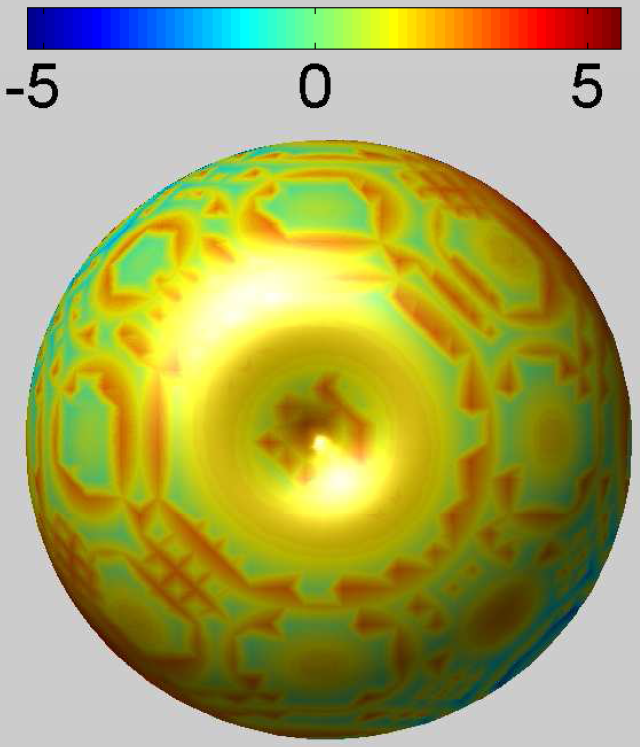}
\includegraphics[width=0.333\textwidth]{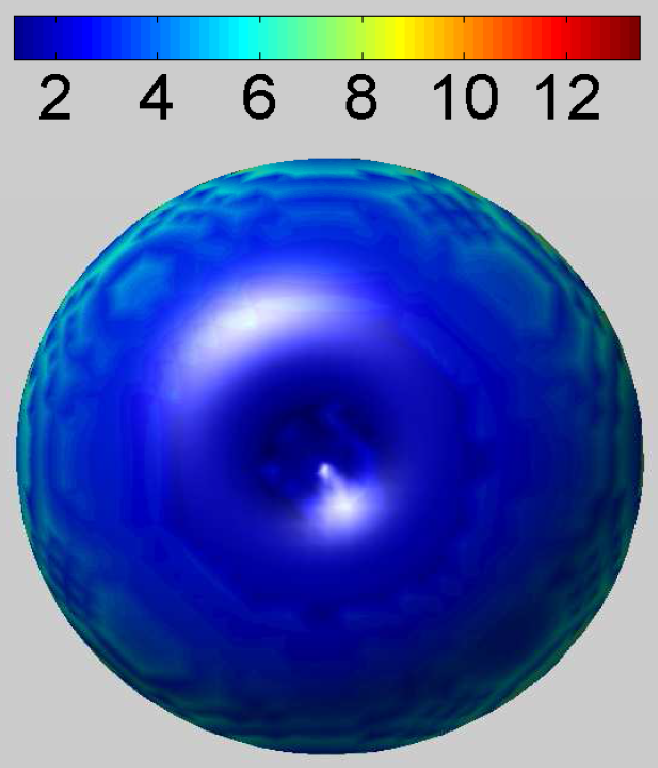}
\includegraphics[width=0.333\textwidth]{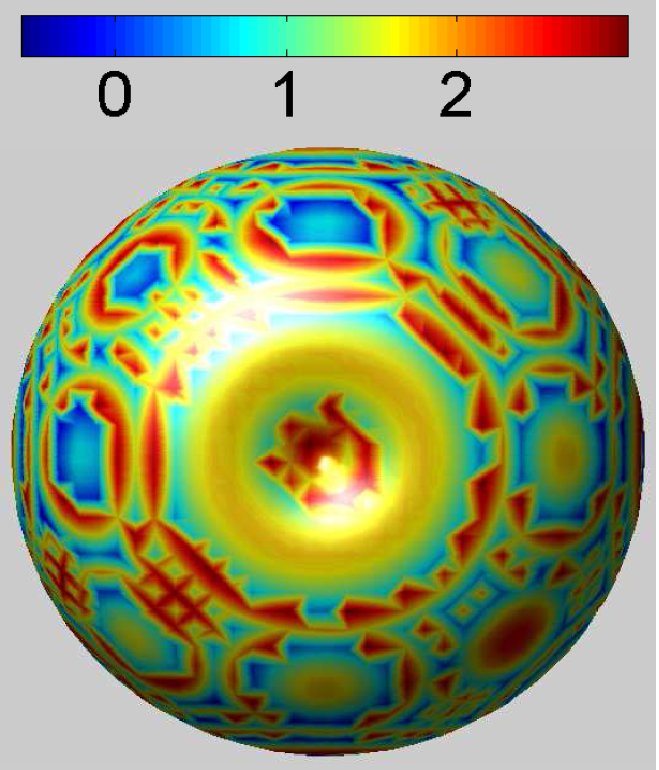}
\end{tabular}
\end{center}
\caption{Solution to the apple-like  interface problem with grid size 0.15. Left chart: $u_1$; Middle chart: $u_2$;  Right chart: $u_3$.}
\label{app_sol}
\end{figure}

\begin{figure}
\begin{center}
\begin{tabular}{ccc}
\includegraphics[width=0.333\textwidth]{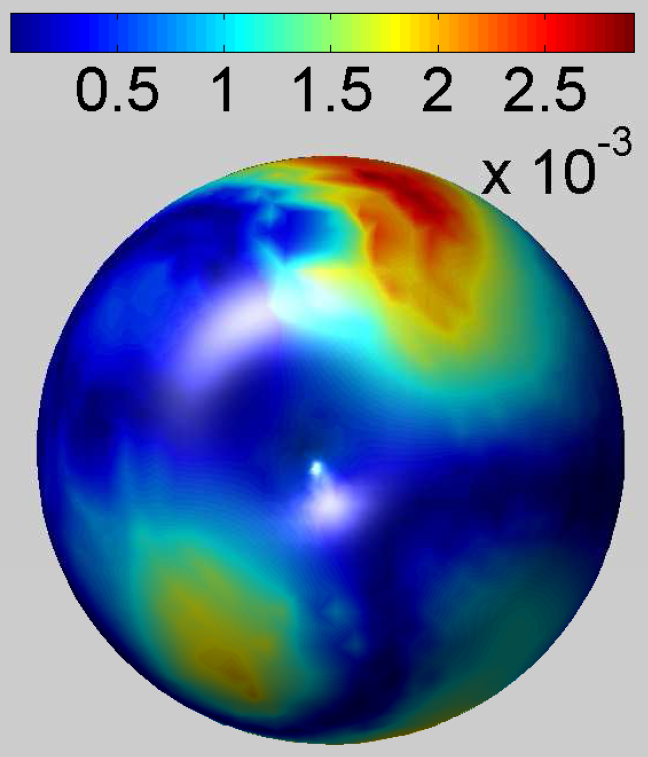}
\includegraphics[width=0.333\textwidth]{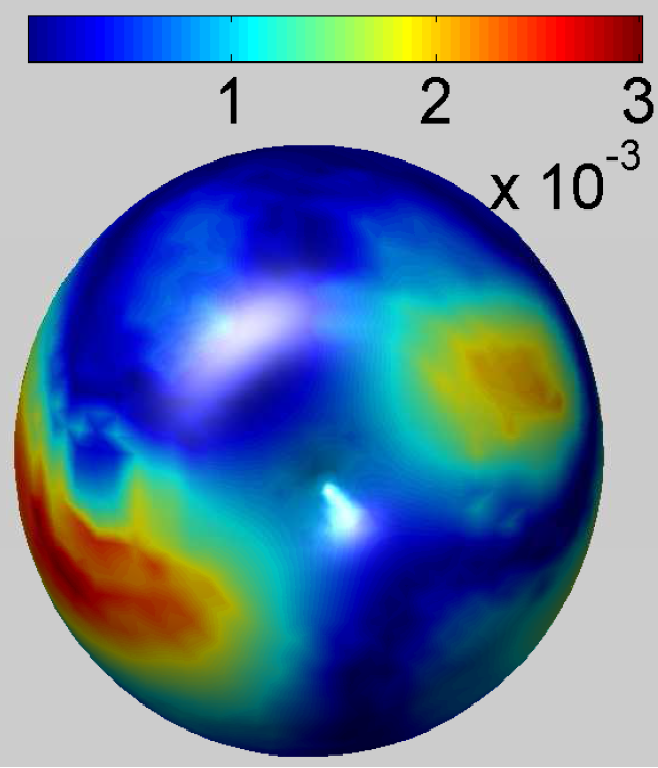}
\includegraphics[width=0.333\textwidth]{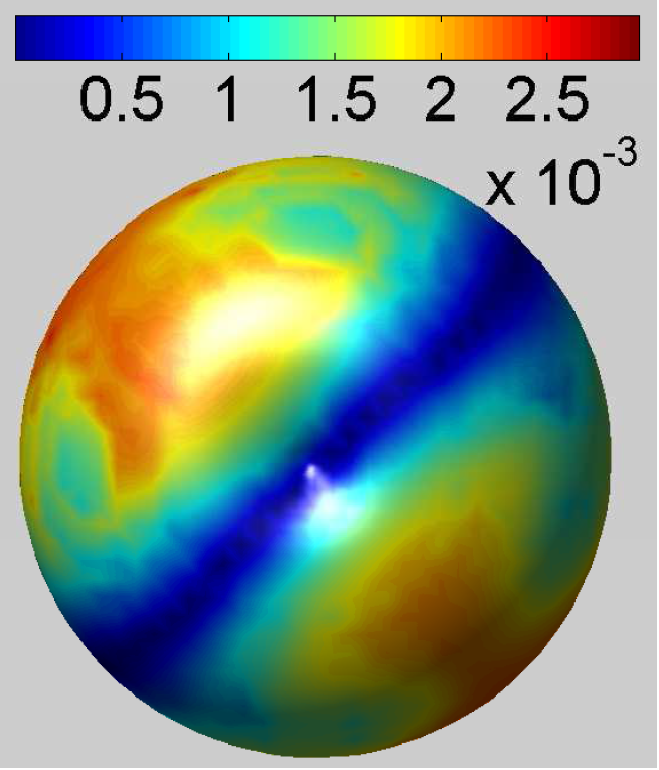}
\end{tabular}
\end{center}
\caption{Numerical error for solving the apple-life interface problem with grid size 0.15. Left chart: $u_1$; Middle chart: $u_2$;  Right chart: $u_3$.}
\label{app_err}
\end{figure}

%\begin{figure}[!ht]
%\small
%\centering
%\includegraphics[width=12.5cm,height=5.5cm]{sol/apple_sol_15.png}
%\caption{Solution to the apple interface with grid size 0.15, left $u_1$, middle $u_2$, right $u_3$}
%\label{app_sol}
%\end{figure}
%
%
%\begin{figure}[!ht]
%\small
%\centering
%\includegraphics[width=12.5cm,height=5.5cm]{sol/apple_err_15.png}
%\caption{Error to the apple interface with grid size 0.15, left $u_1$, middle $u_2$, right $u_3$}
%\label{app_err}
%\end{figure}

Figures \ref{app_sol} and \ref{app_err} illustrate   the numerical solution and error of solving the apple-liked interface problem with grid size $0.15$. The grid refinement analysis  in Table \ref{Ex10_linf} and Table \ref{Ex10_2} indicate the second order convergence in both $L_2$ and $L_\infty$ error norms. Note that largest error did not occur at the geometric singularity, which indicates that  the proposed MIB method works well for geometric singularities.

\textbf{Example~11.}
Next, we consider an oak-acorn interface geometry \cite{Yu:2007a}
$$
\left\{\begin{array}{ll}
\left(\frac{x}{d}\right)^2+\left(\frac{y}{d}\right)^2=(z-q)^2, &\ \  \mbox{if} z>0,\\
x^2+y^2+(z-g)^2=R^2, &\ \  \mbox{if} z\leq 0,
\end{array}\right.
$$
where $q=-\frac{6}{7}$, $g=\frac{1}{2}$, $R=\frac{15}{7}$ and $d=\sqrt{\frac{R^2-g^2}{q^2}}$. The computational domain is set to $[-5, 4.6]\times[-5, 4.6]\times[-5, 4.6]$. Note that this interface has a   tip.

The material parameters and exact solutions in Example 10 are adopted. Grid refinement analysis in terms of $L_\infty$ error is given in Table \ref{Ex11_linf}. In Table \ref{Ex11_2}, similar  analysis in terms of $L_2$ error  is also given. These results show that the second order convergence in both $L_2$ and $L_\infty$ error norms is achieved.

\begin{table}
\caption{The $L_\infty$ errors of Example 11.}
\label{Ex11_linf}
\centering
\begin{tabular}{lllllllll}
\hline
%\multicolumn{5}{c}{The grid refinement analysis of $L_\infty$ error for Example 11} \\
\cline{1-7}
Grid Size       &$L_\infty(u1)$ &Order              &$L_\infty(u2)   $ &Order                &$L_\infty(u3)  $   &Order \\
\hline
$0.48$   &$3.90\times 10^{-2}$     &           &$4.28\times 10^{-2}$      &      &$6.18\times 10^{-2}$     &        \\
$0.24$  &$9.92\times 10^{-3}$     &1.98       &$1.01\times 10^{-2}$     &2.08   &$1.19\times 10^{-2}$     &2.38   \\
$0.12$ &$2.29\times 10^{-3}$     &2.12       &$2.54\times 10^{-3}$      &1.99  &$2.60\times 10^{-3}$     &2.19  \\
\hline
\end{tabular}
%\caption{The grid refinement analysis of $L_\infty$ error for Example 11}
\end{table}

\begin{table}
\caption{The $L_2$ errors of Example 11.}
\label{Ex11_2}
\centering
\begin{tabular}{lllllllll}
\hline
%\multicolumn{5}{c}{The grid refinement analysis of $L_2$ error for Example 11} \\
\cline{1-7}
Grid Size       & $L_2(u1)  $ &Order              &Error  $L_2(u2)    $ &Order                &    $L_2(u3)  $   &Order \\
\hline
$0.48$   &$5.91\times 10^{-3}$     &           &$6.37\times 10^{-3}$      &      &$7.44\times 10^{-3}$     &        \\
$0.24$  &$1.36\times 10^{-3}$     &2.17       &$1.48\times 10^{-3}$     &2.11   &$1.88\times 10^{-3}$     &1.98   \\
$0.12$ &$3.25\times 10^{-4}$     &2.07       &$3.60\times 10^{-4}$      &2.04  &$4.06\times 10^{-4}$     &2.21  \\
\hline
\end{tabular}
%\caption{The grid refinement analysis of $L_2$ error for Example 11}
\end{table}

\begin{figure}
\begin{center}
\begin{tabular}{ccc}
\includegraphics[width=0.333\textwidth]{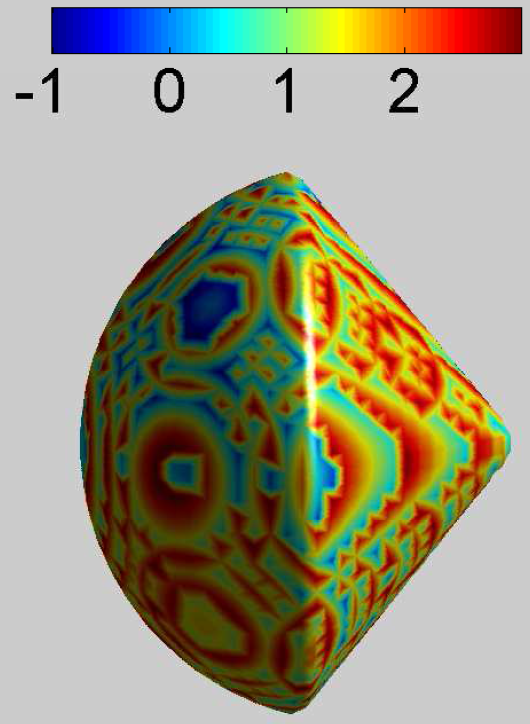}
\includegraphics[width=0.333\textwidth]{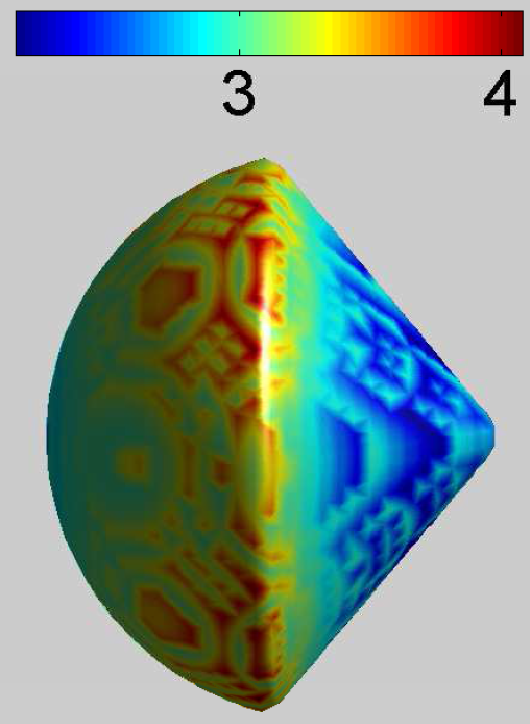}
\includegraphics[width=0.333\textwidth]{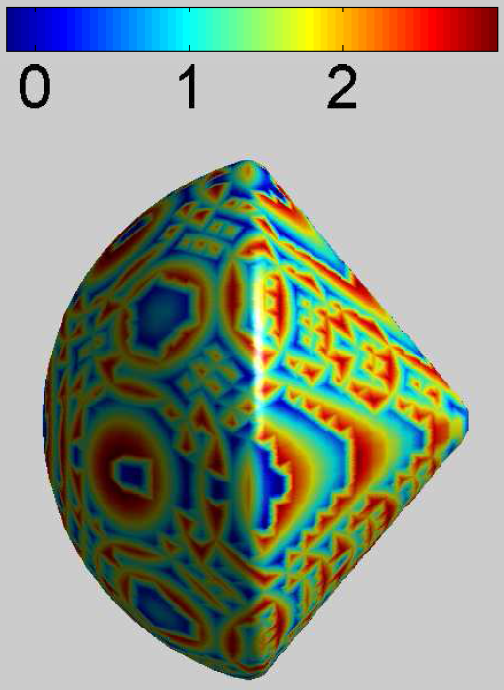}
\end{tabular}
\end{center}
\caption{Numerical solution to the acorn interface problem with grid size 0.12. Left chart: $u_1$; Middle chart: $u_2$;  Right chart: $u_3$.}
\label{acorn_sol}
\end{figure}

\begin{figure}
\begin{center}
\begin{tabular}{ccc}
\includegraphics[width=0.333\textwidth]{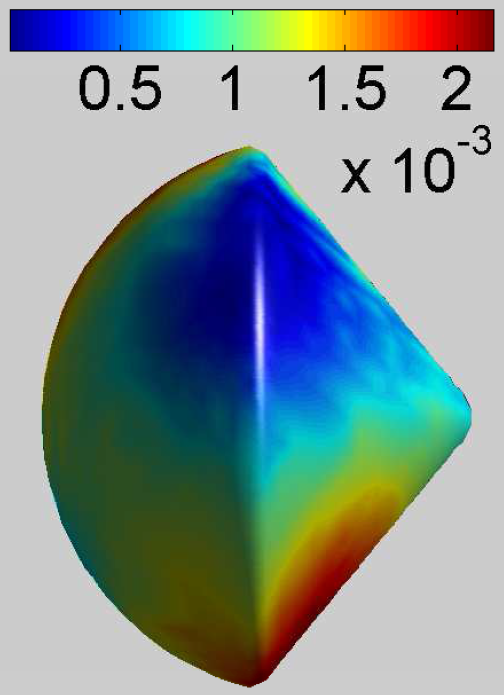}
\includegraphics[width=0.333\textwidth]{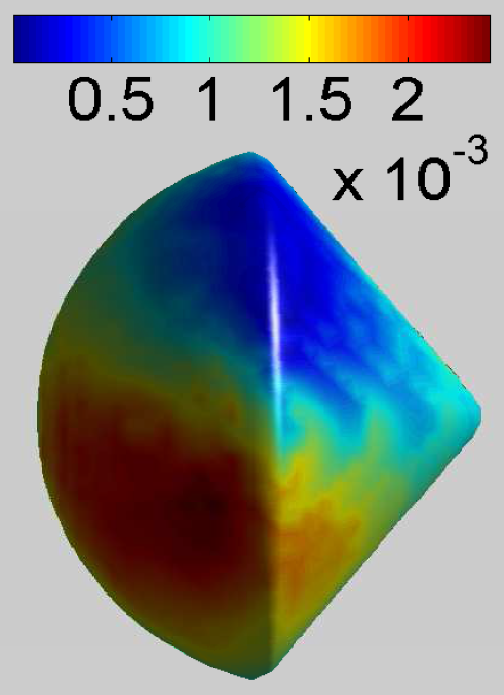}
\includegraphics[width=0.333\textwidth]{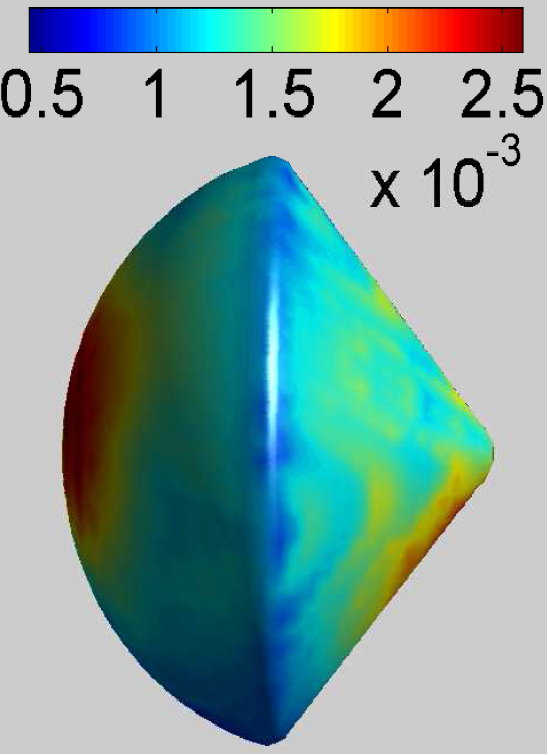}
\end{tabular}
\end{center}
\caption{Numerical error for solving the acorn interface problem with grid size 0.12. Left chart: $u_1$; Middle chart: $u_2$;  Right chart: $u_3$.}
\label{acorn_err}
\end{figure}

%\begin{figure}[!ht]
%\small
%\centering
%\includegraphics[width=12.5cm,height=5.5cm]{sol/acorn_sol_12.png}
%\caption{Solution to the acorn interface with grid size 0.12, left $u_1$, middle $u_2$, right $u_3$}
%\label{acorn_sol}
%\end{figure}
%
%
%\begin{figure}[!ht]
%\small
%\centering
%\includegraphics[width=12.5cm,height=5.5cm]{sol/acorn_err_12.png}
%\caption{Error to the acorn interface with grid size 0.12, left $u_1$, middle $u_2$, right $u_3$}
%\label{acorn_err}
%\end{figure}

The geometry, numerical solution and error distribution are provided in  Figures \ref{acorn_sol} and \ref{acorn_err}, which are computed with grid size $0.15$ in all directions.  Again, the largest error is away from the tip, which indicates the robustness of the present MIB method for dealing with geometric singularity.

\textbf{Example~12.}
Finally, let us extend the benchmark pentagon-star interface test used  in 2D to a 3D one, which is a more complicated interface with very a sharp edge. We set the interface as
$$
\phi(r, \theta)=\left\{\begin{array}{ll}
\frac{R\sin{(\theta_t/2)}}{\sin{(\theta_t/2+\theta-\theta_r-2\pi(i-1)/5})}-r, &\ \  \theta_r+\pi(2i-2)/5\leq \theta <\theta_r+\pi(2i-1)/5,\\
\frac{R\sin{(\theta_t/2)}}{\sin{(\theta_t/2-\theta+\theta_r+2\pi(i-1)/5})}-r, &\ \  \theta_r+\pi(2i-3)/5\leq \theta <\theta_r+\pi(2i-2)/5,
\end{array}\right.
$$
where $\theta_t=\frac{\pi}{5}$, $\theta_t=\frac{\pi}{7}$, $R=\frac{6}{7}$ and $i=1, 2, 3, 4, 5$. Furthermore, we have $r=\sqrt{x^2+y^2}$ and $\theta=\arctan{\frac{y}{x}}$.
The $z$-direction of the interface is constrained by
$$
-\frac{\sqrt{3}}{2}\leq z \leq \frac{\sqrt{3}}{2}.
$$
The computational domain is set to $[-1.3, 1.1]\times[-1.3, 1.1]\times[-1.3, 1.1]$.
The material parameters and exact solutions in Example 10 are utilized for this problem.

Grid refinement $L_\infty$ error analysis is given in  Table \ref{Ex12_linf}. We also shown the grid refinement analysis in terms of  $L_2$ error in Table \ref{Ex12_2}. The grid refinement analysis  shows that the second order convergence in both $L_2$ and $L_\infty$ error norms is obtained.

\begin{table}
\caption{The $L_\infty$ errors of Example 12.}
\label{Ex12_linf}
\centering
\begin{tabular}{lllllllll}
\hline
%\multicolumn{5}{c}{The grid refinement analysis of $L_\infty$ error for Example 12} \\
\cline{1-7}
Grid Size       &$L_\infty(u1)$ &Order              &$L_\infty(u2)   $ &Order                &$L_\infty(u3)  $   &Order \\
\hline
$0.12$   &$1.27\times 10^{-3}$     &           &$1.57\times 10^{-3}$      &       &$1.75\times 10^{-2}$     &        \\
$0.06$   &$2.08\times 10^{-4}$     &2.61       &$3.98\times 10^{-4}$      &1.98   &$1.76\times 10^{-4}$     &3.31   \\
$0.03$   &$3.80\times 10^{-5}$     &2.45       &$6.56\times 10^{-5}$      &2.60   &$2.79\times 10^{-5}$     &2.66  \\
\hline
\end{tabular}
%\caption{The grid refinement analysis of $L_\infty$ error for Example 12}
\end{table}

\begin{table}
\caption{The $L_2$ errors of Example 12.}
\label{Ex12_2}
\centering
\begin{tabular}{lllllllll}
\hline
%\multicolumn{5}{c}{The grid refinement analysis of $L_2$ error for Example 12} \\
\cline{1-7}
Grid Size       & $L_2(u1)  $ &Order              &Error  $L_2(u2)    $ &Order                &    $L_2(u3)  $   &Order \\
\hline
$0.12$   &$3.71\times 10^{-4}$     &           &$1.67\times 10^{-4}$      &       &$3.57\times 10^{-4}$     &        \\
$0.06$   &$3.79\times 10^{-5}$     &3.29       &$4.47\times 10^{-5}$      &1.90   &$3.76\times 10^{-5}$     &3.25   \\
$0.03$   &$7.67\times 10^{-6}$     &2.30       &$1.07\times 10^{-5}$      &2.60   &$6.08\times 10^{-6}$     &2.63  \\
\hline
\end{tabular}
%\caption{The grid refinement analysis of $L_2$ error for Example 12}
\end{table}

\begin{figure}[!ht]
\begin{center}
\begin{tabular}{ccc}
\includegraphics[width=0.333\textwidth]{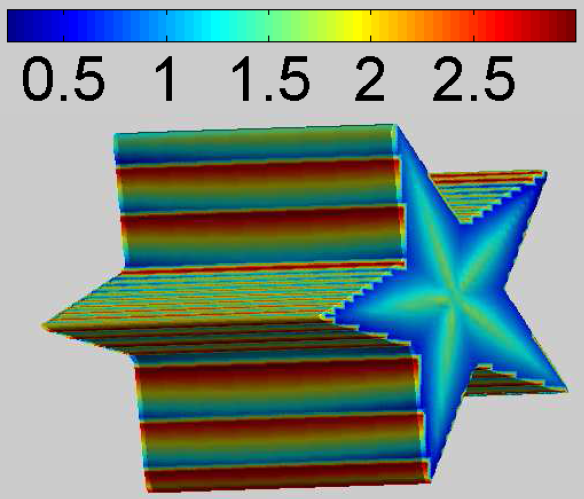}
\includegraphics[width=0.333\textwidth]{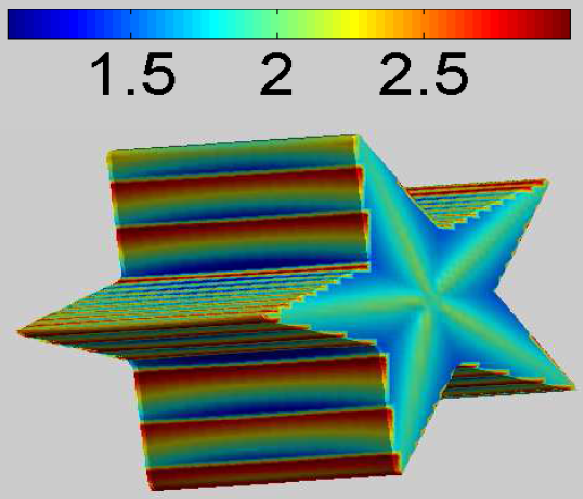}
\includegraphics[width=0.333\textwidth]{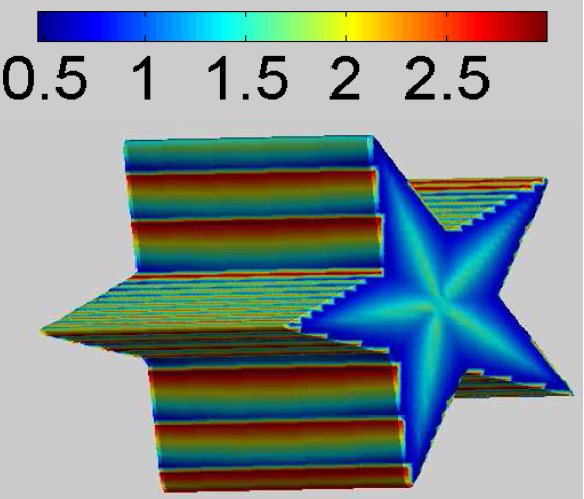}
\end{tabular}
\end{center}
\caption{Numerical solution to the pentagon star interface problem with grid size 0.03. Left chart: $u_1$; Middle chart: $u_2$;  Right chart: $u_3$.  }
\label{pen_sol}
\end{figure}

\begin{figure}[!ht]
\begin{center}
\begin{tabular}{ccc}
\includegraphics[width=0.333\textwidth]{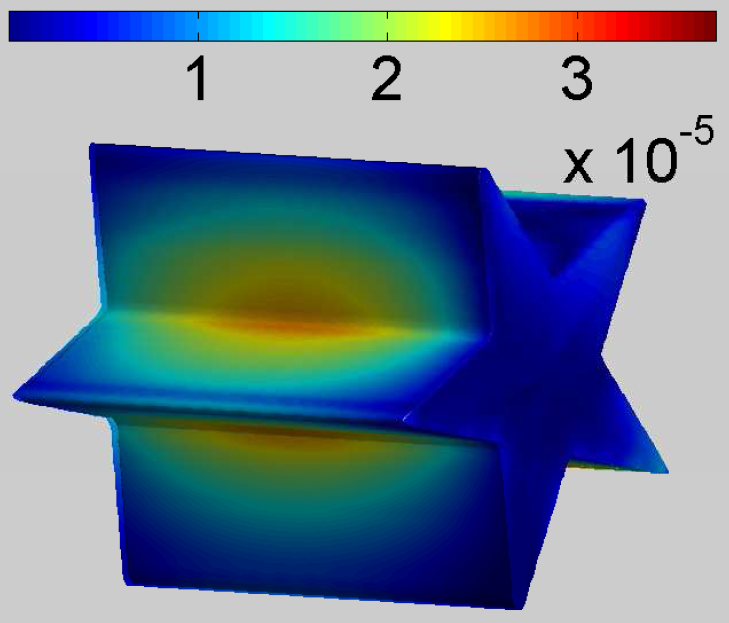}
\includegraphics[width=0.333\textwidth]{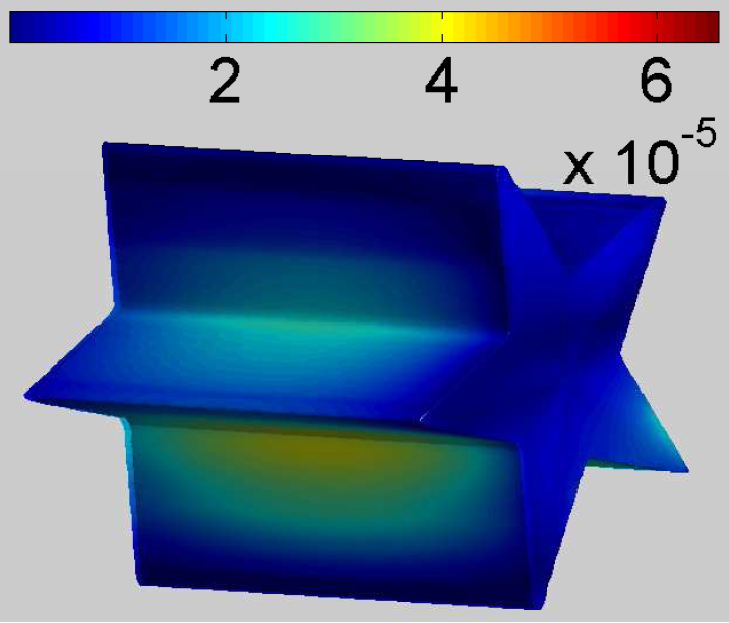}
\includegraphics[width=0.333\textwidth]{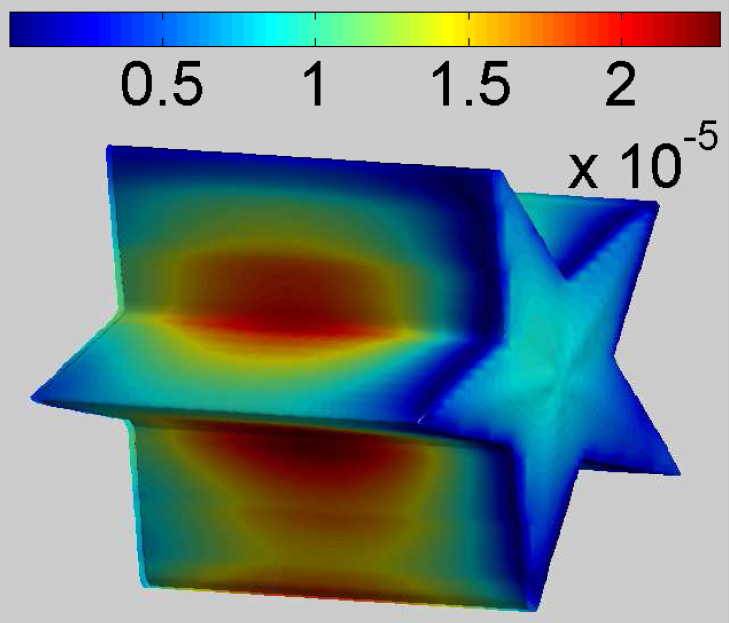}
\end{tabular}
\end{center}
\caption{Numerical error for solving  the pentagon star interface problem with grid size 0.03. Left chart: $u_1$; Middle chart: $u_2$;  Right chart: $u_3$.}
\label{pen_err}
\end{figure}

To give a visualization of the numerical solution, error and interface geometry of   the pentagon star liked interface problem, we provide
Figs. \ref{pen_sol} and \ref{pen_err}, which are plotted with grid size $0.03$. In general, error is very small. Additionally, the largest error does not occur at the sharp edge of the interface. Therefore, from the above three test examples,  we can conclude that the proposed MIB method is very robust in handling geometric singularities.

\section{Conclusion}

In this work, we develop the  matched interface and boundary (MIB) method for solving  three dimensional (3D) elasticity interface problems. Both isotropic homogeneous material and  isotropic inhomogeneous material are considered in the theoretical modeling and numerical computation.   In particular, the isotropic inhomogeneous material  is described by a strain-stress constitutive law with a position-dependent modulus  function.

Most previous effort in the MIB method has been for elliptic interface problems.  Its essential idea is to replace function values   on irregular grid points fictitious values in the discretization so that the standard finite difference schemes can be systematically employed as if there were no interface. Interface jump conditions are enforced on the intersecting points between the interface and the mesh lines, which in turn determines  fictitious values. In principle, the MIB method developed for one interface problem can be utilized for solving another interface problem because the MIB procedure does depend on the form of the partial differential equation. However,  elasticity interface equations are exceptional because they involve both central derivatives and cross derivatives, which lead to new difficulties in determining fictitious values. Additionally, the elasticity interface equation  is a vector equation with three deformation components in a 3D setting, which results in more demanding in efficient numerical schemes in terms of computer memory storage and convergent speed in solving the linear algebraic system.   Consequently, a new MIB method has been  developed in this work to address these issues.   To make the MIB scheme of second order convergence, a number of techniques for central derivatives and cross derivatives is proposed in this work. For central derivatives, techniques such as  local coordinate transformation, disassociation, two sets of jump conditions are utilized, while  for cross derivatives, disassociation, extrapolation and neighbor combination techniques are proposed to determining fictitious values.  The resulting large sparse linear systems for the coupled vector equations are solved efficiently by using the bi-conjugated gradient method.

The proposed MIB method has been  validated by using a variety of benchmark examples. In terms of interface complexity, we considered both smooth interfaces and nonsmooth interfaces.  Smooth interface geometries include sphere, hemisphere, genus-1 torus, flower and cylinder. In the category of nonsmooth interface geometries,  apple-shaped,  oak-acorn-shaped and pentagon star interfaces are considered. It is well-known that in order to achieve second order convergence, nonsmooth interface geometries require special considerations in the interface algorithm design. The robustness of the MIB method proved by showing that the largest error occurs away from the geometric singularities.

The proposed MIB method has also been tested for  elasticity interface problems with both weak discontinuity and strong discontinuity in solutions. Another standard test is the stability of the numerical schemes for large contrasts in material parameters across the interface. These aspects are
investigated with numerous examples. We have demonstrated that the proposed MIB method is not sensitive to change in discontinuity and material contrast.

Finally, two classes of material parameters, namely, piecewise constants and spatially varying Poisson's ratio and  shear modulus are considered in our numerical experiments. We have demonstrated with extensive numerical examples that proposed MIB method achieve second order convergence in both $L_\infty$ and $L_2$ error  norms for all the tests described above. Additionally, the level of MIB accuracy is not affected the above mentioned test issues. We therefore believe that the present MIB is ready for applications to the real world problems. In fact, application to complex biomolecular systems is under our consideration.

\section*{Acknowledgments}

This work was supported in part by NSF grants   IIS-1302285 and DMS-1160352,   NIH grant R01GM-090208 and MSU Center for Mathematical Molecular Biosciences Initiative.

\clearpage
\section*{Literature cited}
\renewcommand\refname{}

% Use alpha to check for repeated references
\bibliographystyle{abbrv}
\bibliography{refs}

\end{document}